\documentclass{article}

\usepackage{arxiv}

\usepackage[utf8]{inputenc} 
\usepackage[T1]{fontenc}    
\usepackage{hyperref}       
\usepackage{url}            
\usepackage{booktabs}       
\usepackage{amsmath,amsfonts,amssymb}       
\usepackage{nicefrac}       
\usepackage{microtype}      
\usepackage{cleveref}       
\usepackage{lipsum}         
\usepackage{graphicx}
\usepackage[square,numbers]{natbib}
\usepackage{doi}
\usepackage[dvipsnames]{xcolor}
\usepackage{bm}
\usepackage{derivative}
\usepackage{caption}
\usepackage{cancel}
\usepackage{comment}
\usepackage{algorithm}
\usepackage{float}
\usepackage{tikz}
\usetikzlibrary{positioning, shapes, decorations.markings, decorations.pathreplacing, decorations.pathmorphing, arrows, calc, matrix, fit, backgrounds}
\usepackage{subcaption}
\usepackage{mathtools}
\usepackage{makecell}
\usepackage[graphicx]{realboxes}
\usepackage{pdflscape}
\usepackage{svg}
\usepackage{tcolorbox}
\usepackage{enumitem}

\newcommand{\frulargs}[2]{f_{R U L_{#1}}\left(#2\right)}

\newcommand{\dt}{\mathrm{d} t}
\newcommand{\dr}{\mathrm{d} r}
\newcommand{\EE}[1]{\mathbb{E}\left[#1\right]}
\newcommand{\hEE}[1]{\widehat{\mathbb{E}}\left[#1\right]}
\newcommand{\Var}[1]{\mathbb{V}\mathrm{ar}\left[#1\right]}
\newcommand{\hVar}[1]{\widehat{\mathbb{V}\mathrm{ar}}\left[#1\right]}
\newcommand{\trep}{t_{\mathrm{rep}}}

\newcommand{\Pro}[1]{\mathrm{Pr}\left(#1\right)}

\newcommand{\Tfi}{T_\mathrm{F}^{(i)}}
\newcommand{\Tf}{T_{\mathrm{F}}}

\title{Branching Out:\\
Prognostics-Based Replacement Policies for Series Systems
}

\author{Daniel Koutas \\
	ERA Group\\
	TU Munich\\
	Munich, Germany \\
	\texttt{daniel.koutas@tum.de} \\
	\And
	Daniel Straub \\
	ERA Group\\
	TU Munich\\
	Munich, Germany \\
	\texttt{straub@tum.de} \\
}

\renewcommand{\shorttitle}{Branching Out}

\begin{document}

\newcommand{\workflow}{
\begin{figure}[H]
\def\ldist{3}
\def\adist{3}
\def\awidth{2.5}
\def\apicfac{0.27}
\def\apicheight{0.88}
\def\sydist{1.7}
\def\swidth{2.3}
\def\sheight{1.46}
\def\boxtextyshift{0.5}
\def\xpad{0.2}
\def\ypad{0.2}
\pgfmathsetmacro\radius{1*\sheight}
\pgfmathsetmacro\modelwidth{0.7*\radius}
\pgfmathsetmacro\sensorwidth{0.95*\swidth}

\def\termradius{0.2}
\def\boundline{very thick}
\def\arrowthick{very thick}
\def\dashthick{thick}
\tikzset{circle_node/.style={circle, draw, \boundline, minimum size=\radius cm, inner sep=0pt}}
\tikzset{square_node/.style={rectangle, draw, align=center, rounded corners, \boundline, minimum width=\swidth cm, minimum height=\sheight cm, inner sep=0pt}}
\tikzset{ghost_node/.style={circle, minimum size=0 cm, inner sep=0pt}}
    \centering
    \begin{tikzpicture}
    \coordinate (s1cord) at (0,0);
    \node[square_node, minimum width=\swidth cm, minimum height=\sheight cm] (s1) at (s1cord) {\includegraphics[width=\sensorwidth cm]{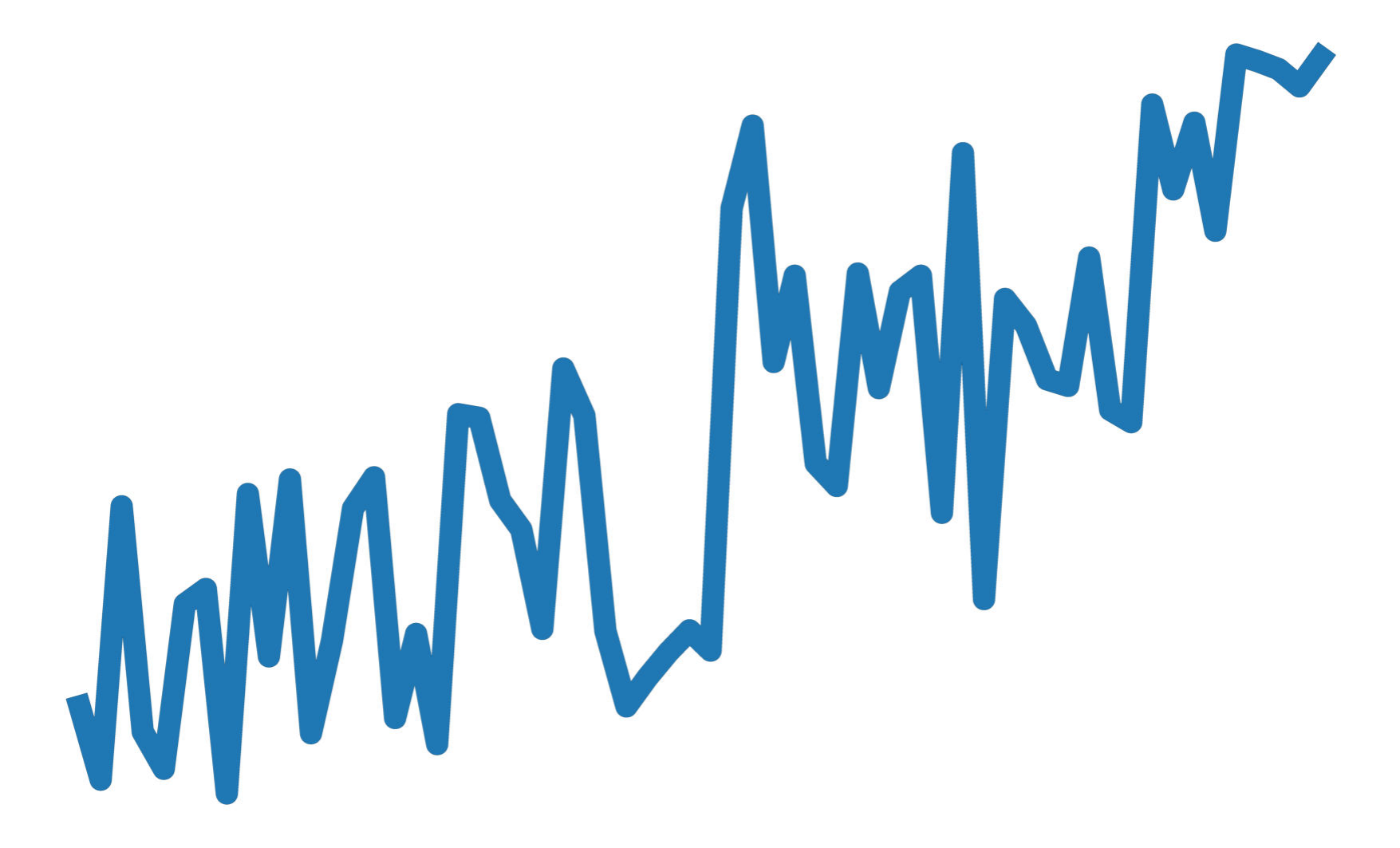}};
    \coordinate[right=\ldist cm of s1cord] (m1cord);
    \node[circle_node] (m1) at (m1cord) {\includegraphics[width=\modelwidth cm]{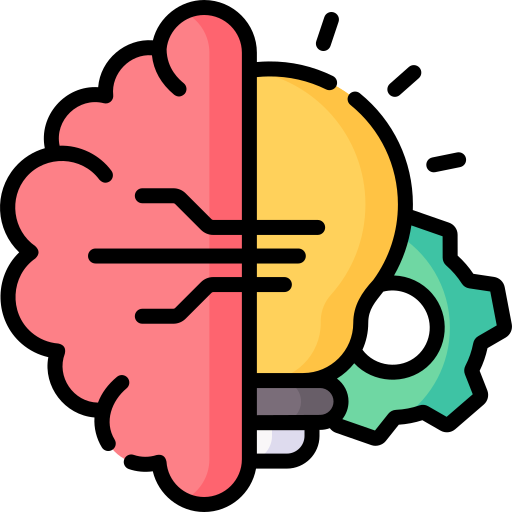}};
    \coordinate[right=\ldist cm of m1cord] (r1cord);
    \node[square_node] (r1) at (r1cord) {\includegraphics[width=\sensorwidth cm]{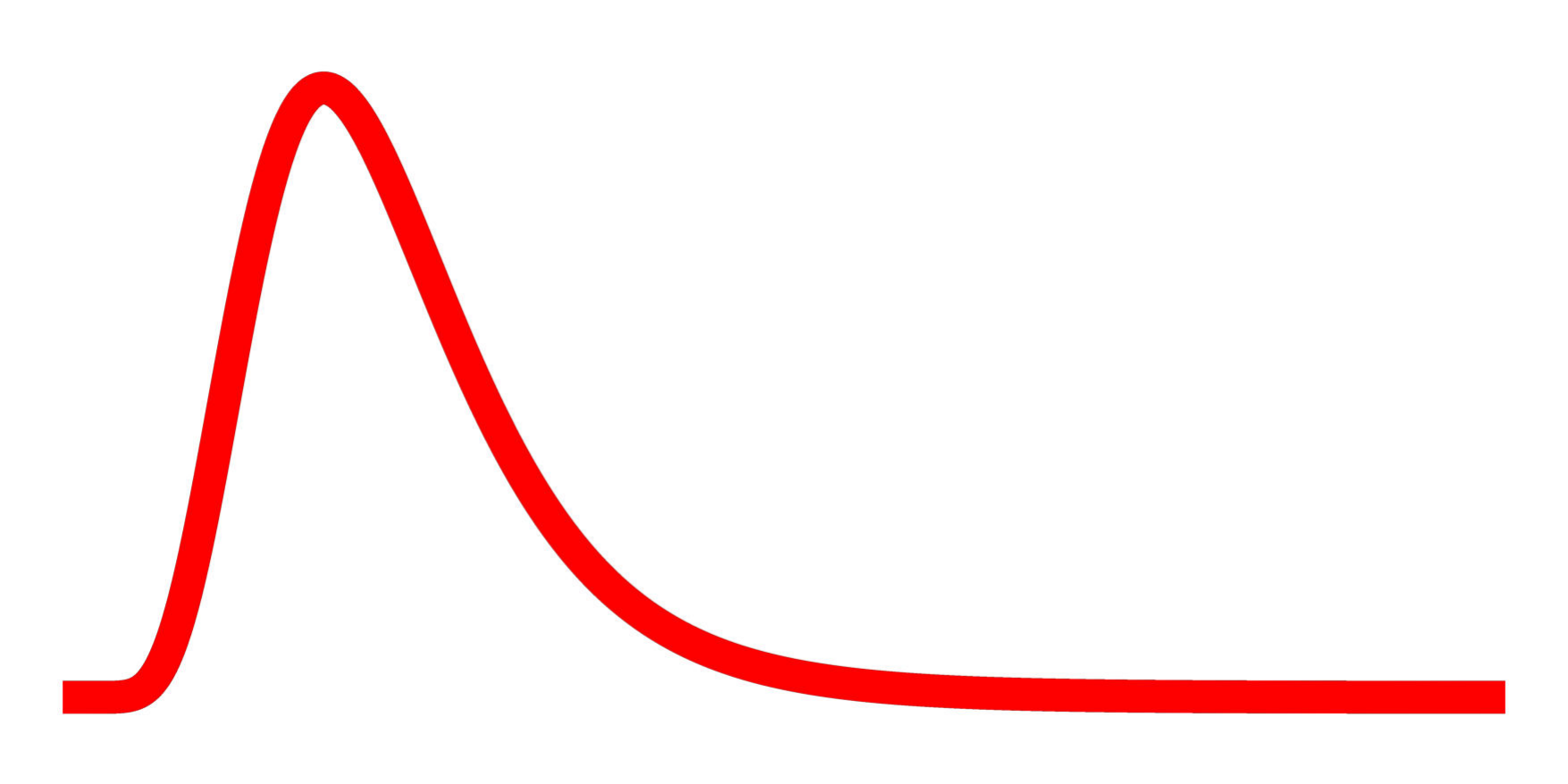}};
    \coordinate[below=\sydist cm of s1cord] (s2cord);
    \node[ghost_node] (s2) at (s2cord) {\Large $\vdots$};
    \coordinate[below=\sydist cm of m1cord] (m2cord);
    \node[ghost_node] (m2) at (m2cord) {\Large $\vdots$};
    \coordinate[right=\ldist cm of m2cord] (r2cord);
    \node[ghost_node] (r2) at (r2cord) {\Large $\vdots$};
    \coordinate[below=\sydist cm of s2cord] (sMcord);
    \node[square_node,] (sM) at (sMcord) {\includegraphics[width=\sensorwidth cm]{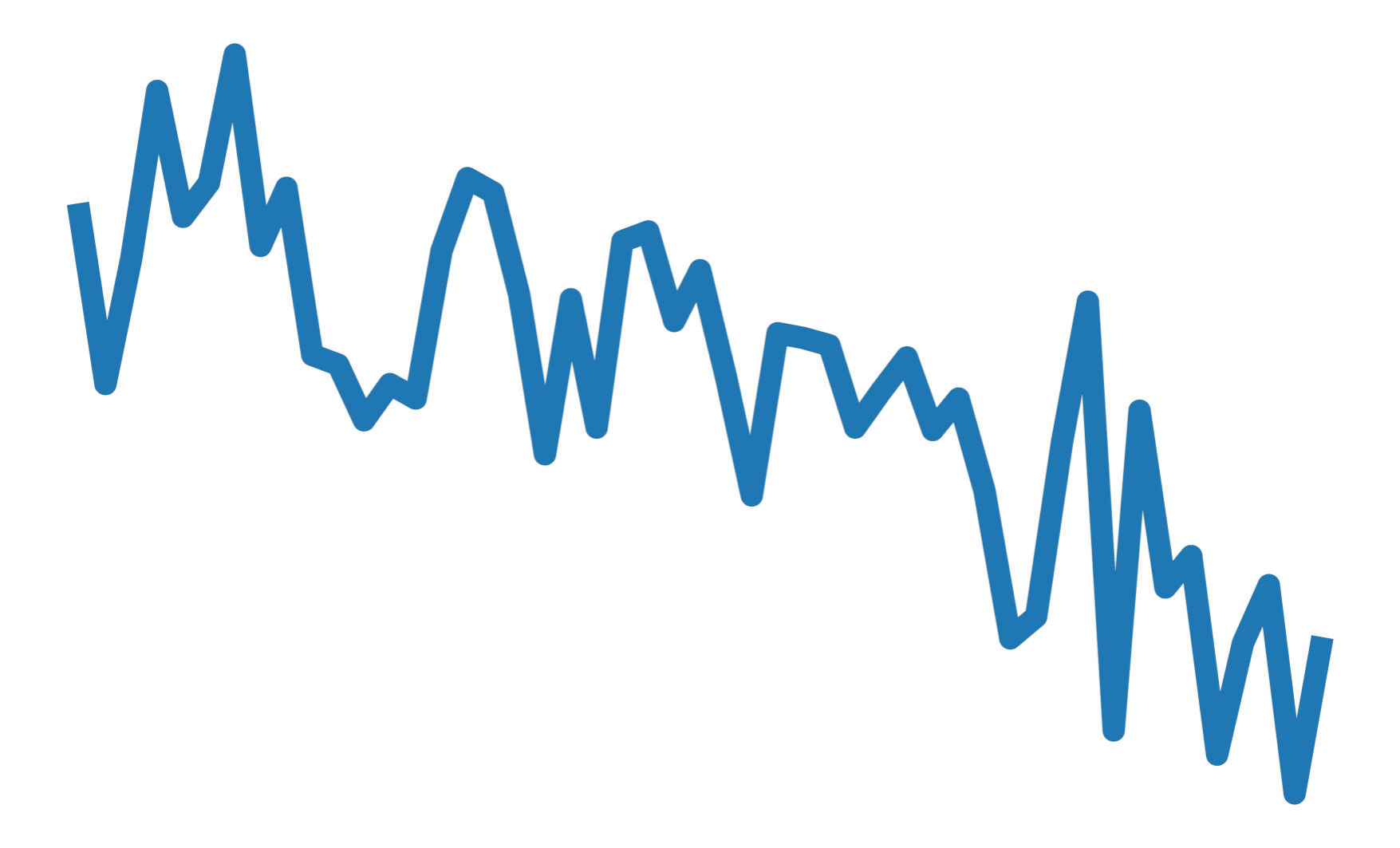}};
    \coordinate[right=\ldist cm of sMcord] (mMcord);
    \node[circle_node] (mM) at (mMcord) {\includegraphics[width=\modelwidth cm]{./ML_model_icon.png}};
    \coordinate[right=\ldist cm of mMcord] (rMcord);
    \node[square_node] (rM) at (rMcord) {\includegraphics[width=\sensorwidth cm]{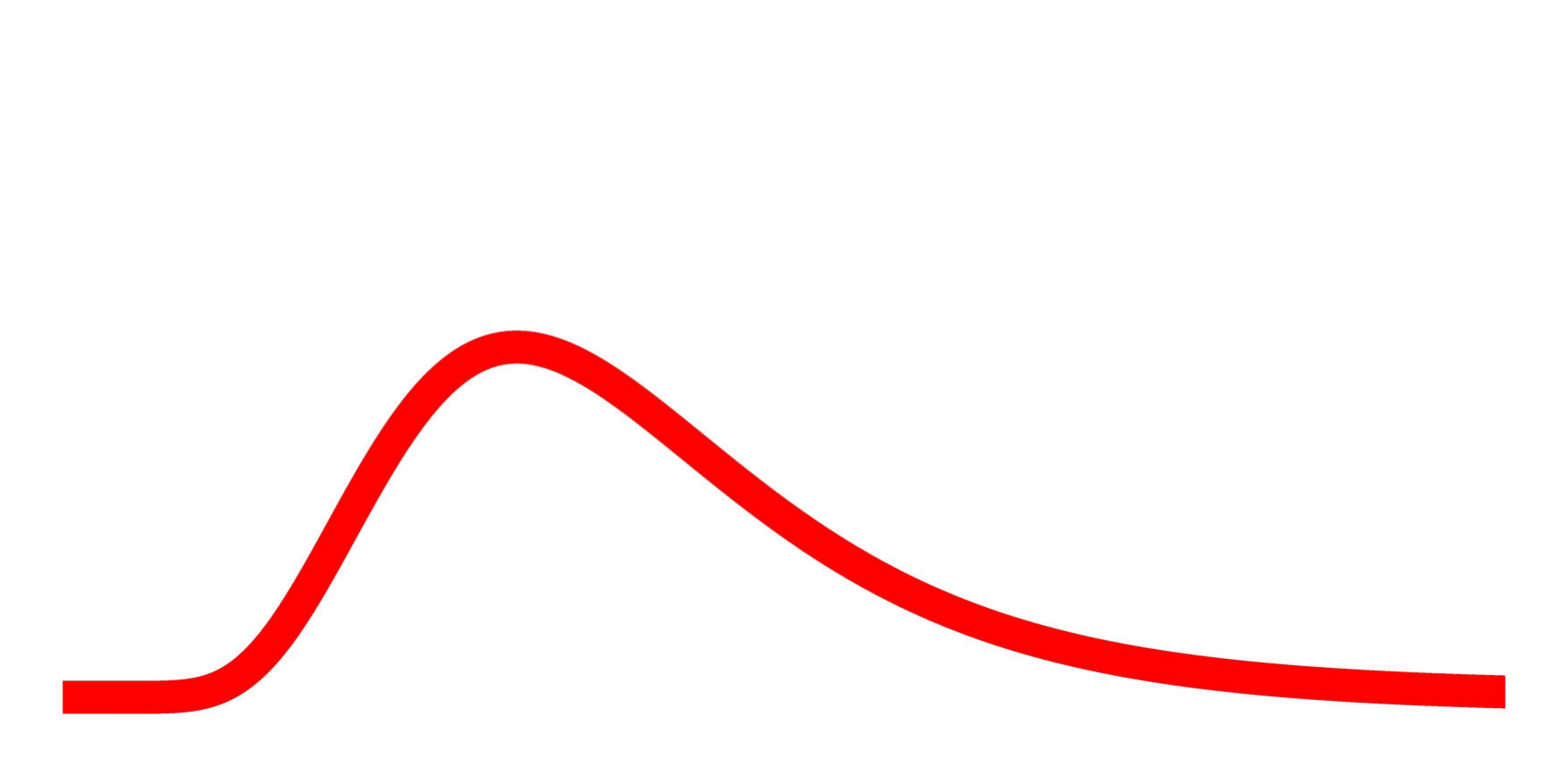}};
    
    \coordinate[right=\ldist cm of r2cord] (picord);
    \node[circle_node, minimum size=\sheight cm] (pi) at (picord) {\includegraphics[width=\modelwidth cm]{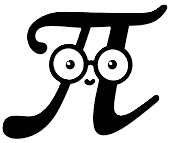}};

    \coordinate[right=\adist cm of picord] (a2cord);
    \node[ghost_node, inner sep=3pt] (a2) at (a2cord) {\Large $\vdots$};
    \coordinate[above=\sydist cm of a2cord] (a1cord);
    \node[square_node, inner sep=3pt, minimum width=\awidth cm, minimum height=\sheight cm] (a1) at (a1cord) {\Huge /};
    \coordinate[left=\apicfac*\awidth of a1cord] (a1dncord);
    \node[ghost_node] (a1dn) at (a1dncord){\includegraphics[height=\apicheight cm]{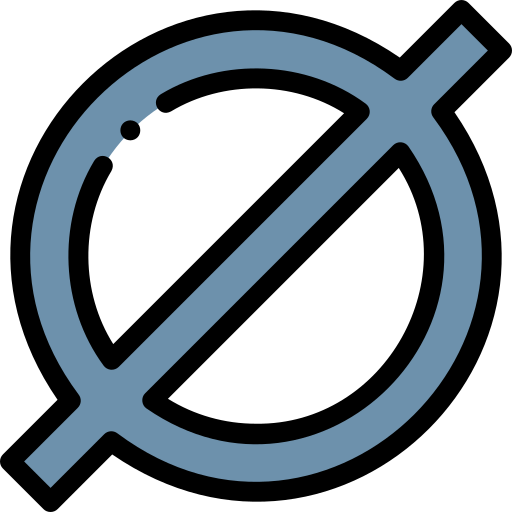}};
    \coordinate[right=\apicfac*\awidth of a1cord] (a1prcord);
    \node[ghost_node] (a1pr) at (a1prcord){\includegraphics[height=\apicheight cm]{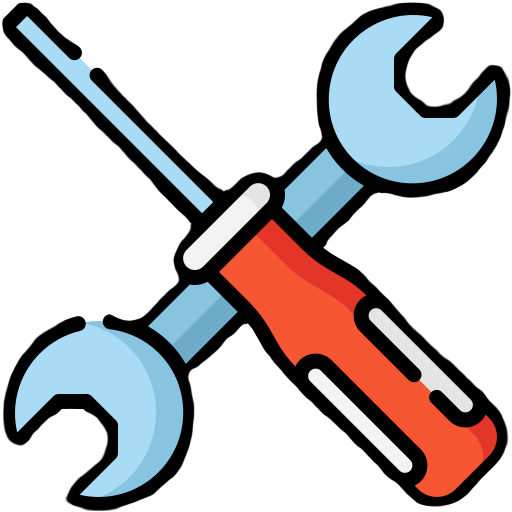}};
    \coordinate[below=\sydist cm of a2cord] (aMcord);
    \node[square_node, inner sep=3pt, minimum width=\awidth cm, minimum height=\sheight cm] (aM) at (aMcord) {\Huge /};
    \coordinate[left=\apicfac*\awidth of aMcord] (aMdncord);
    \node[ghost_node] (aMdn) at (aMdncord){\includegraphics[height=\apicheight cm]{./nothing.png}};
    \coordinate[right=\apicfac*\awidth of aMcord] (aMprcord);
    \node[ghost_node] (aMpr) at (aMprcord){\includegraphics[height=\apicheight cm]{./maint_nb.png}};

    \draw[-latex, \arrowthick] (s1) -- (m1);
    \draw[-latex, \arrowthick] (m1) -- (r1);
    \draw[-latex, \arrowthick] (sM) -- (mM);
    \draw[-latex, \arrowthick] (mM) -- (rM);
    \draw[-latex, \arrowthick] (r1) -- (pi);
    \draw[-latex, \arrowthick] (rM) -- (pi);
    \draw[-latex, \arrowthick] (pi) -- (a1dn);
    \draw[-latex, \arrowthick] (pi) -- (aMdn);

    \def\padtop{1.4}
    \def\padleft{1}
    \def\padright{0.4}
    \def\padbot{1.4}
    \begin{pgfonlayer}{background}

        \fill[teal!18!white, rounded corners]
        ($(s1.north west) + 
        (-\padleft,\padtop)$) --
        ($(a1.north east) + 
        (\padright,\padtop)$) --
        ($(aM.south east) + 
        (\padright,-\padbot)$) --
        ($(sM.south west) + 
        (-\padleft,-\padbot)$) -- 
        ($(s1.north west) + 
        (-\padleft,\padtop)$);
    \end{pgfonlayer}
    \node[anchor=south west, yshift=0pt] at ($(sM.south west) + (-\padleft,-\padtop)$)
    {{\large \textcolor{teal!90!black}{Repeat at each time step $t_k$}}};

    \def\off{1}
    \pgfmathsetmacro\dashxdist{0.5*\swidth+\xpad}
    \pgfmathsetmacro\dashydist{0.5*\sheight+\ypad}
    
    \coordinate[left=\dashxdist cm of s1cord] (bslcord);
    \coordinate[above=\dashydist cm of s1cord] (bstcord);
    \coordinate[right=\dashxdist cm of s1cord] (bsrcord);
    \coordinate[below=\dashydist cm of sMcord] (bsbcord);
    \coordinate[above=\off cm of bstcord] (bst2cord);
    \coordinate[right=\dashxdist cm of bstcord] (bsr2cord);
    \coordinate[left=\dashxdist cm of bstcord] (bsl2cord);
    \coordinate[above=\boxtextyshift cm of bstcord] (bstextcord);
    \node[ghost_node] (bstext) at (bstextcord) {SHM data};
    \def\comptextpad{1.7}
    \coordinate[left=\comptextpad cm of s1cord] (comp1textcord);
    \node[ghost_node, rotate=90] (comp1text) at (comp1textcord) {Comp. $1$};
    \coordinate[left=\comptextpad cm of sMcord] (compMtextcord);
    \node[ghost_node, rotate=90] (compMtext) at (compMtextcord) {Comp. $M$};

    \pgfmathsetmacro\dashxdistc{0.5*\radius+\xpad}
    \coordinate[left=\dashxdistc cm of m1cord] (bmlcord);
    \coordinate[above=\dashydist cm of m1cord] (bmtcord);
    \coordinate[right=\dashxdistc cm of m1cord] (bmrcord);
    \coordinate[below=\dashydist cm of mMcord] (bmbcord);
    \coordinate[above=\off cm of bmtcord] (bmt2cord);
    \coordinate[right=\dashxdistc cm of bmtcord] (bmr2cord);
    \coordinate[left=\dashxdistc cm of bmtcord] (bml2cord);
    \coordinate[above=\boxtextyshift cm of bmtcord] (bmtextcord);
    \node[ghost_node, align=center] (bmtext) at (bmtextcord) {Prognostic \\ models};

    \coordinate[left=\dashxdist cm of r1cord] (brlcord);
    \coordinate[above=\dashydist cm of r1cord] (brtcord);
    \coordinate[right=\dashxdist cm of r1cord] (brrcord);
    \coordinate[below=\dashydist cm of rMcord] (brbcord);
    \coordinate[above=\off cm of brtcord] (brt2cord);
    \coordinate[right=\dashxdist cm of brtcord] (brr2cord);
    \coordinate[left=\dashxdist cm of brtcord] (brl2cord);
    \coordinate[above=\boxtextyshift cm of brtcord] (brtextcord);
    \node[ghost_node, align=center] (brtext) at (brtextcord) {Predicted \\RUL distributions};

    \pgfmathsetmacro\dashxdista{0.5*\awidth+\xpad}
    \coordinate[left=\dashxdista cm of a1cord] (balcord);
    \coordinate[above=\dashydist cm of a1cord] (batcord);
    \coordinate[right=\dashxdista cm of a1cord] (barcord);
    \coordinate[below=\dashydist cm of aMdncord] (babcord);
    \coordinate[above=\off cm of batcord] (bat2cord);
    \coordinate[right=\dashxdista cm of batcord] (bar2cord);
    \coordinate[left=\dashxdista cm of batcord] (bal2cord);
    \coordinate[above=\boxtextyshift cm of batcord] (batextcord);
    \node[ghost_node, align=center] (batext) at (batextcord) {Do nothing / \\ Maintain};

    \pgfmathsetmacro\dashxdistpi{0.5*\sheight+\xpad}
    \coordinate[left=\dashxdistpi cm of picord] (bpilcord);
    \coordinate[above=\dashydist cm of picord] (bpitcord);
    \coordinate[right=\dashxdistpi cm of picord] (bpircord);
    \coordinate[below=\dashydist cm of picord] (bpibcord);
    \coordinate[above=\off cm of bpitcord] (bpit2cord);
    \coordinate[right=\dashxdistpi cm of bpitcord] (bpir2cord);
    \coordinate[left=\dashxdistpi cm of bpitcord] (bpil2cord);
    \coordinate[above=\boxtextyshift cm of bpitcord] (pitextcord);
    \node[ghost_node, align=center] (pitext) at (pitextcord) {Policy};

    \def\boxcolor{white}

    \begin{pgfonlayer}{background}
        \node[
            fit=(bslcord) (bst2cord) (bsrcord) (bsbcord), 
            fill=\boxcolor,
            draw=black,
            very thick,
            inner sep=0pt,
            rounded corners
        ] (bluebox) {};
        \draw[very thick] (bsl2cord) -- (bsr2cord);

        \node[
            fit=(bmlcord) (bmt2cord) (bmrcord) (bmbcord), 
            fill=\boxcolor,
            draw=black,
            very thick,
            inner sep=0pt,
            rounded corners
        ] (bluebox) {};
        \draw[very thick] (bml2cord) -- (bmr2cord);
        
        \node[
            fit=(brlcord) (brt2cord) (brrcord) (brbcord), 
            fill=\boxcolor,
            draw=black,
            very thick,
            inner sep=0pt,
            rounded corners
        ] (bluebox) {};
        \draw[very thick] (brl2cord) -- (brr2cord);
        
        \node[
            fit=(bpilcord) (bpit2cord) (bpircord) (bpibcord), 
            fill=\boxcolor,
            draw=black,
            very thick,
            inner sep=0pt,
            rounded corners
        ] (bluebox) {};
        \draw[very thick] (bpil2cord) -- (bpir2cord);
        
        \node[
            fit=(balcord) (bat2cord) (barcord) (babcord), 
            fill=\boxcolor,
            draw=black,
            very thick,
            inner sep=0pt,
            rounded corners
        ] (bluebox) {};
        \draw[very thick] (bal2cord) -- (bar2cord);
    \end{pgfonlayer}

    \end{tikzpicture}
    \captionsetup{width=\textwidth}
    \caption{General workflow\protect\footnotemark of multi-component predictive maintenance policies. At every time step $t_k$, new SHM data is gathered for each component, fed into the component-specific prognostic models, which output corresponding RUL distributions. These are then collectively fed into a policy ($\pi$), which outputs a joint action vector, assigning each component either a “do nothing” or “maintain” decision.}
    \label{fig:decision_setting}
\end{figure}
\footnotetext{Prognostic model, do nothing and maintain icons are taken from Flaticon.com, and the policy icon from Vecteezy.com}
}

\newcommand{\decisiontree}{
\begin{figure}[H]
\def\xdist{5}
\def\sydist{2.4}
\def\aydist{0.65}
\def\xtextdist{0.4}
\def\ytextdist{1.0}
\def\tdist{9.4}
\def\width{0.4}
\def\radius{0.4}
\def\termradius{0.2}
\def\capsize{\normalsize}
\def\boundline{very thick}
\def\arrowthick{very thick}
\def\tshift{6}
\def\fixshift{0.4}
\def\safecolor{green!60!black}
\def\failcolor{red!80!black}
\def\actfac{3}
\def\condpos{0.75}
\def\pfpos{0.65}
\def\xenumshift{0.6}
\def\yenumshift{0.25}
\tikzset{circle_node/.style={circle, draw, \boundline, minimum size=\radius cm, inner sep=0pt}}
\tikzset{failed_node/.style={circle, draw, fill=black, \boundline, minimum size=\termradius cm, inner sep=0pt}}
\tikzset{semi_replaced_node/.style={circle, draw, fill=orange!40!white, \boundline, minimum size=\radius cm, inner sep=0pt}}
\tikzset{replaced_node/.style={circle, draw, fill=orange!70!black, \boundline, minimum size=\radius cm, inner sep=0pt}}
\tikzset{act_node/.style={rectangle, draw, fill=cyan, \boundline, minimum size=\width cm, inner sep=0pt}}
\tikzset{cost_node/.style={rectangle, draw, fill=lime, \boundline, minimum size=1.2*\width cm, inner xsep=2pt, inner ysep=0pt}}

    \centering
    \vspace{-1.2cm}
    \begin{tikzpicture}
    \coordinate (a0cord) at (0,0);
    \node[act_node] (a0) at (a0cord) {};
    
    \coordinate[right=\xdist cm of a0cord, yshift=\sydist cm] (s1dnprcord) {};
    \node[semi_replaced_node] (s1dnpr) at (s1dnprcord) {};
    \coordinate[right=\xdist cm of a0cord, yshift=-\sydist cm] (s1prdncord) {};
    \node[semi_replaced_node] (s1prdn) at (s1prdncord) {};
    \coordinate[above=2*\sydist cm of s1dnprcord] (s1dndncord) {};
    \node[circle_node] (s1dndn) at (s1dndncord) {};
    \coordinate[below=2*\sydist cm of s1prdncord] (s1prprcord) {};
    \node[replaced_node] (s1prpr) at (s1prprcord) {};

    \coordinate[right=\xdist cm of s1dndncord, yshift=\actfac*\aydist cm] (dndna1cord) {};
    \node[act_node] (dndna1) at (dndna1cord) {};
    \coordinate[right=\xdist cm of s1dndncord, yshift=\aydist cm] (dndnx1fcord) {};
    \node[failed_node] (dndnx1f) at (dndnx1fcord) {};
    \coordinate[right=\xdist cm of s1dndncord, yshift=-\aydist cm] (dndnx2fcord) {};
    \node[failed_node] (dndnx2f) at (dndnx2fcord) {};
    \coordinate[right=\xdist cm of s1dndncord, yshift=-\actfac*\aydist cm] (dndnx1fx2fcord) {};
    \node[failed_node] (dndnx1fx2f) at (dndnx1fx2fcord) {};

    \coordinate[right=\xdist cm of s1dnprcord, yshift=\actfac*\aydist cm] (dnpra1cord) {};
    \node[act_node] (dnpra1) at (dnpra1cord) {};
    \coordinate[right=\xdist cm of s1dnprcord, yshift=\aydist cm] (dnprx1fcord) {};
    \node[failed_node] (dnprx1f) at (dnprx1fcord) {};
    \coordinate[right=\xdist cm of s1dnprcord, yshift=-\aydist cm] (dnprx2fcord) {};
    \node[failed_node] (dnprx2f) at (dnprx2fcord) {};
    \coordinate[right=\xdist cm of s1dnprcord, yshift=-\actfac*\aydist cm] (dnprx1fx2fcord) {};
    \node[failed_node] (dnprx1fx2f) at (dnprx1fx2fcord) {};

    \coordinate[right=\xdist cm of s1prdncord, yshift=\actfac*\aydist cm] (prdna1cord) {};
    \node[act_node] (prdna1) at (prdna1cord) {};
    \coordinate[right=\xdist cm of s1prdncord, yshift=\aydist cm] (prdnx1fcord) {};
    \node[failed_node] (prdnx1f) at (prdnx1fcord) {};
    \coordinate[right=\xdist cm of s1prdncord, yshift=-\aydist cm] (prdnx2fcord) {};
    \node[failed_node] (prdnx2f) at (prdnx2fcord) {};
    \coordinate[right=\xdist cm of s1prdncord, yshift=-\actfac*\aydist cm] (prdnx1fx2fcord) {};
    \node[failed_node] (prdnx1fx2f) at (prdnx1fx2fcord) {};

    \coordinate[right=\xdist cm of s1prprcord, yshift=\actfac*\aydist cm] (prpra1cord) {};
    \node[act_node] (prpra1) at (prpra1cord) {};
    \coordinate[right=\xdist cm of s1prprcord, yshift=\aydist cm] (prprx1fcord) {};
    \node[failed_node] (prprx1f) at (prprx1fcord) {};
    \coordinate[right=\xdist cm of s1prprcord, yshift=-\aydist cm] (prprx2fcord) {};
    \node[failed_node] (prprx2f) at (prprx2fcord) {};
    \coordinate[right=\xdist cm of s1prprcord, yshift=-\actfac*\aydist cm] (prprx1fx2fcord) {};
    \node[failed_node] (prprx1fx2f) at (prprx1fx2fcord) {};

    \coordinate[right=\xtextdist cm of dndna1cord] (dndna1costcord) {};
    \node[cost_node, anchor=west] (dndna1cost) at (dndna1costcord) {$c_I^1$};
    \coordinate[right=\xtextdist cm of dndnx1fcord] (dndnx1fcostcord) {}; 
    \node[cost_node, anchor=west] (dndnx1fcost) at (dndnx1fcostcord) {$c_{F,1}+c_I^2$};
    \coordinate[right=\xtextdist cm of dndnx2fcord] (dndnx2fcostcord) {};
    \node[cost_node, anchor=west] (dndnx2fcost) at (dndnx2fcostcord) {$c_{F,2}+c_I^3$};
    \coordinate[right=\xtextdist cm of dndnx1fx2fcord] (dndnx1fx2fcostcord) {};
    \node[cost_node, anchor=west] (dndnx1fx2fcost) at (dndnx1fx2fcostcord) {$c_{F,1}+c_{F,2}+c_I^4$};

    \coordinate[right=\xtextdist cm of dnpra1cord] (dnpra1costcord) {};
    \node[cost_node, anchor=west] (dnpra1cost) at (dnpra1costcord) {$c_I^1$};
    \coordinate[right=\xtextdist cm of dnprx1fcord] (dnprx1fcostcord) {}; 
    \node[cost_node, anchor=west] (dnprx1fcost) at (dnprx1fcostcord) {$c_{F,1}+c_I^2$};
    \coordinate[right=\xtextdist cm of dnprx2fcord] (dnprx2fcostcord) {};
    \node[cost_node, anchor=west] (dnprx2fcost) at (dnprx2fcostcord) {$c_{F,2}+c_I^3$};
    \coordinate[right=\xtextdist cm of dnprx1fx2fcord] (dnprx1fx2fcostcord) {};
    \node[cost_node, anchor=west] (dnprx1fx2fcost) at (dnprx1fx2fcostcord) {$c_{F,1}+c_{F,2}+c_I^4$};
    \coordinate[below=\fixshift cm of s1dnprcord] (dnprfixcord) {};
    \node[cost_node, anchor=north] (dnprfix) at (dnprfixcord) {$c_f+c_{v,2}$};

    \coordinate[right=\xtextdist cm of prdna1cord] (prdna1costcord) {};
    \node[cost_node, anchor=west] (prdna1cost) at (prdna1costcord) {$c_I^1$};
    \coordinate[right=\xtextdist cm of prdnx1fcord] (prdnx1fcostcord) {}; 
    \node[cost_node, anchor=west] (prdnx1fcost) at (prdnx1fcostcord) {$c_{F,1}+c_I^2$};
    \coordinate[right=\xtextdist cm of prdnx2fcord] (prdnx2fcostcord) {};
    \node[cost_node, anchor=west] (prdnx2fcost) at (prdnx2fcostcord) {$c_{F,2}+c_I^3$};
    \coordinate[right=\xtextdist cm of prdnx1fx2fcord] (prdnx1fx2fcostcord) {};
    \node[cost_node, anchor=west] (prdnx1fx2fcost) at (prdnx1fx2fcostcord) {$c_{F,1}+c_{F,2}+c_I^4$};
    \coordinate[below=\fixshift cm of s1prdncord] (prdnfixcord) {};
    \node[cost_node, anchor=north] (prdnfix) at (prdnfixcord) {$c_f+c_{v,1}$};

    \coordinate[right=\xtextdist cm of prpra1cord] (prpra1costcord) {};
    \node[cost_node, anchor=west] (prpra1cost) at (prpra1costcord) {$c_I^1$};
    \coordinate[right=\xtextdist cm of prprx1fcord] (prprx1fcostcord) {}; 
    \node[cost_node, anchor=west] (prprx1fcost) at (prprx1fcostcord) {$c_{F,1}+c_I^2$};
    \coordinate[right=\xtextdist cm of prprx2fcord] (prprx2fcostcord) {};
    \node[cost_node, anchor=west] (prprx2fcost) at (prprx2fcostcord) {$c_{F,2}+c_I^3$};
    \coordinate[right=\xtextdist cm of prprx1fx2fcord] (prprx1fx2fcostcord) {};
    \node[cost_node, anchor=west] (prprx1fx2fcost) at (prprx1fx2fcostcord) {$c_{F,1}+c_{F,2}+c_I^4$};
    \coordinate[below=\fixshift cm of s1prprcord] (prprfixcord) {};
    \node[cost_node, anchor=north] (prprfix) at (prprfixcord) {$c_f + c_{v,1}+c_{v,2}$};


    \draw[\arrowthick] (a0) -- (s1dndn) node[midway, sloped, above] {DN,DN};
    \draw[\arrowthick] (a0) -- (s1dnpr) node[midway, sloped, above] {DN,PR};
    \draw[\arrowthick] (a0) -- (s1prdn) node[midway, sloped, above] {PR,DN};
    \draw[\arrowthick] (a0) -- (s1prpr) node[midway, sloped, above] {PR,PR};

    \draw[\arrowthick] (s1dndn) -- (dndna1) 
        node[pos=\condpos, sloped, yshift=\tshift pt] {\textcolor{\safecolor}{Safe}}
        node[pos=\pfpos, sloped, yshift=-\tshift pt] {$(1-p_{F,1})\cdot(1-p_{F,2})$};
    \draw[\arrowthick] (s1dndn) -- (dndnx1f) 
        node[pos=\condpos, sloped, yshift=\tshift pt] {\textcolor{\failcolor}{$X_1$ fail}}
        node[pos=\pfpos, sloped, yshift=-\tshift pt] {$p_{F,1}\cdot(1-p_{F,2})$};
    \draw[\arrowthick] (s1dndn) -- (dndnx2f) 
        node[pos=\condpos, sloped, yshift=\tshift pt] {\textcolor{\failcolor}{$X_2$ fail}}
        node[pos=\pfpos, sloped, yshift=-\tshift pt] {$(1-p_{F,1})\cdot p_{F,2}$};
    \draw[\arrowthick] (s1dndn) -- (dndnx1fx2f) 
        node[pos=\condpos, sloped, yshift=\tshift pt] {\textcolor{\failcolor}{$X_1~\&~X_2$ fail}}
        node[pos=\pfpos, sloped, yshift=-\tshift pt] {$p_{F,1}\cdot p_{F,2}$};

    \draw[\arrowthick] (s1dnpr) -- (dnpra1) 
        node[pos=\condpos, sloped, yshift=\tshift pt] {\textcolor{\safecolor}{Safe}}
        node[pos=\pfpos, sloped, yshift=-\tshift pt] {$1-p_{F,1}$};
    \draw[\arrowthick] (s1dnpr) -- (dnprx1f) 
        node[pos=\condpos, sloped, yshift=\tshift pt] {\textcolor{\failcolor}{$X_1$ fail}}
        node[pos=\pfpos, sloped, yshift=-\tshift pt] {$p_{F,1}$};
    \draw[\arrowthick, dashed] (s1dnpr) -- (dnprx2f) 
        node[pos=\condpos, sloped, yshift=\tshift pt] {\textcolor{\failcolor}{$X_2$ fail}}
        node[pos=\pfpos, sloped, yshift=-\tshift pt] {$\approx 0$};
    \draw[\arrowthick, dashed] (s1dnpr) -- (dnprx1fx2f) 
        node[pos=\condpos, sloped, yshift=\tshift pt] {\textcolor{\failcolor}{$X_1~\&~X_2$ fail}}
        node[pos=\pfpos, sloped, yshift=-\tshift pt] {$\approx 0$};
        
    \draw[\arrowthick] (s1prdn) -- (prdna1) 
        node[pos=\condpos, sloped, yshift=\tshift pt] {\textcolor{\safecolor}{Safe}}
        node[pos=\pfpos, sloped, yshift=-\tshift pt] {$1-p_{F,2}$};
    \draw[\arrowthick, dashed] (s1prdn) -- (prdnx1f) 
        node[pos=\condpos, sloped, yshift=\tshift pt] {\textcolor{\failcolor}{$X_1$ fail}}
        node[pos=\pfpos, sloped, yshift=-\tshift pt] {$\approx 0$};
    \draw[\arrowthick] (s1prdn) -- (prdnx2f) 
        node[pos=\condpos, sloped, yshift=\tshift pt] {\textcolor{\failcolor}{$X_2$ fail}}
        node[pos=\pfpos, sloped, yshift=-\tshift pt] {$p_{F,2}$};
    \draw[\arrowthick, dashed] (s1prdn) -- (prdnx1fx2f) 
        node[pos=\condpos, sloped, yshift=\tshift pt] {\textcolor{\failcolor}{$X_1~\&~X_2$ fail}}
        node[pos=\pfpos, sloped, yshift=-\tshift pt] {$\approx 0$};
        
    \draw[\arrowthick] (s1prpr) -- (prpra1) 
        node[pos=\condpos, sloped, yshift=\tshift pt] {\textcolor{\safecolor}{Safe}}
        node[pos=\pfpos, sloped, yshift=-\tshift pt] {$1-p_{F,2}$};
    \draw[\arrowthick, dashed] (s1prpr) -- (prprx1f) 
        node[pos=\condpos, sloped, yshift=\tshift pt] {\textcolor{\failcolor}{$X_1$ fail}}
        node[pos=\pfpos, sloped, yshift=-\tshift pt] {$\approx 0$};
    \draw[\arrowthick, dashed] (s1prpr) -- (prprx2f) 
        node[pos=\condpos, sloped, yshift=\tshift pt] {\textcolor{\failcolor}{$X_2$ fail}}
        node[pos=\pfpos, sloped, yshift=-\tshift pt] {$\approx 0$};
    \draw[\arrowthick, dashed] (s1prpr) -- (prprx1fx2f) 
        node[pos=\condpos, sloped, yshift=\tshift pt] {\textcolor{\failcolor}{$X_1~\&~X_2$ fail}}
        node[pos=\pfpos, sloped, yshift=-\tshift pt] {$\approx 0$};
        

    %
    
    \coordinate[above=\tdist cm of a0cord] (t0cord) {};
    \coordinate (t1cord) at ($(t0cord) + (2*\xdist, 0)$) {};
    \coordinate (t2cord) at ($(t1cord) + (3*\xdist, 0)$) {};
    \coordinate (tprevcord) at ($(t0cord) + (-3*\xdist, 0)$) {};

    \draw[decorate,decoration={brace,amplitude=10pt}] (t0cord) -- (t1cord) 
        node[midway,yshift=16pt] {\large Step $k$};
    \begin{scope}
        \clip (t1cord) circle (2.1 cm);
        \draw[decorate,decoration={brace,amplitude=10pt}] (t1cord) -- (t2cord) node[midway, xshift=-6.5cm, yshift=16pt] {\large Step $k+1$};
    \end{scope}
    \begin{scope}
        \clip (t0cord) circle (2.1 cm);
        \draw[decorate,decoration={brace,amplitude=10pt}] (tprevcord) -- (t0cord) node[midway, xshift=6.5cm, yshift=16pt] {\large Step $k-1$};
    \end{scope}

    \coordinate (dll) at ($(a0cord) + (0, -9.2)$){};
    \draw[thick, dashed, gray] ($(t0cord) + (0, 0.8)$) -- (a0) -- (dll){};
    \draw[thick, dashed, gray] 
    ($(t1cord) + (0, 0.8)$) -- (dndna1) -- (dndnx1f) -- 
    (dndnx2f) -- (dndnx1fx2f) -- (dnpra1) -- (dnprx1f) -- 
    (dnprx2f) -- (dnprx1fx2f) -- (prdna1) -- (prdnx1f) -- 
    (prdnx2f) -- (prdnx1fx2f) -- (prpra1) -- (prprx1f) -- 
    (prprx2f) -- (prprx1fx2f) --
    ($(prprx1fx2fcord)!1!(prprx1fx2fcord|-dll)$) {};
    
    \end{tikzpicture}
    \captionsetup{width=\textwidth}
    \caption{Decision tree for a series system of two deteriorating components that can be preventively replaced. System states are represented by round nodes, where white, light orange, and brown fillings indicate no, partial, and full system replacement, respectively. Maintenance decisions ($\bm{a}_k$) are represented as square cyan nodes. The action taken for each component (PR/DN) as well as the events (Fail/Safe) with the corresponding probabilities are included in the respective paths, where $p_{F,i}=\Pr(RUL_i\leq\Delta t)$, $c_{F,i}=c_c-c_{\infty,i}\cdot\EE{RUL_i \mid RUL_i \leq \Delta t}$ and $c_I^j$ depends on the specific policy, as described in detail in \Cref{subsubsec:general_structure}. System failure due to failure of at least one component is shown with filled black circles. Expected costs due to actions or failures are shown with lime boxes where they are incurred.}
    \label{fig:dt}
\end{figure}
}

\pgfmathdeclarefunction{lognormal}{3}{%
\pgfmathparse{%
    (exp(-(ln(#1+0.01)-#2)^2 / (2*#3^2)) / (#1*#3*sqrt(2)+0.01))%
    }%
}

\newcommand{\repgraphs}{
\begin{figure}[H]
\def\xdist{7.4}
\def\tlabelshift{-0.5}
\def\ytickheight{0.2}
\def\xtk{0.3}
\def\tickdist{4}
\def\ylabelshift{1}
\def\linethick{ultra thick}
\def\tickthick{very thick}
\def\rulcolor{red}
\def\domain{1.2:7.1}
\def\mu{1.4}
\def\sigma{0.4}
\def\Tfone{2.5}
\def\Tftwo{5.6}
\def\Tfonecolor{blue!80!black}
\def\Tftwocolor{green!50!black}
\def\plotyshift{0.1}
\def\plotfactor{2}
\def\textxshift{0}
\def\textyshift{1.3}
\def\ctickheight{0.1}
\def\yclabelshift{1}
\def\ycshift{1.2}
\def\cccolor{\Tfonecolor}
\def\cpcolor{\Tftwocolor}
\def\ectrcolor{orange!90!black}
\def\colorshift{0.04}
\def\yrshift{0.9}
\def\xrshift{1.0}
\def\rthick{thick}
\def\rtickheight{0.15}
\def\rlabelshift{0.8}


\centering
\begin{subfigure}{0.49\textwidth}
    \centering
    \begin{tikzpicture}
        \coordinate (0) at (0,0);
        \coordinate (end) at ($(0) + (\xdist, 0)$);
        \node[yshift=\tlabelshift cm] at (end) {$t$};
        \coordinate (tk) at ($(\xtk,0)$);
        \draw[\tickthick] ($(tk) + (0,\ytickheight)$) -- ($(tk) + 
        (0,-\ytickheight)$) node[anchor=north, yshift=-\ylabelshift]{\large $t_k$};
        \coordinate (tk1) at ($(tk) + (\tickdist,0)$);
        \draw[\tickthick] ($(tk1) + (0,\ytickheight)$) -- ($(tk1) + 
        (0,-\ytickheight)$) node[anchor=north, yshift=-\ylabelshift]{\large $t_{k+1}$};
    
        \draw[\rulcolor, very thick, domain=\domain, samples=200] 
        plot ({\x}, {\plotfactor*lognormal(\x,\mu,\sigma)+\plotyshift});
        \node[] at ($(tk1) + (\textxshift, \textyshift)$) {$\color{\rulcolor} \frulargs{i}{t}$};
    
        \draw[-latex, \linethick] (0) -- (end);
    
        \coordinate (cp) at ($(tk) + (0,-\ycshift)$);
        
        \draw[dotted, \tickthick, \cpcolor] (tk) -- ($(tk) + (0,-0.4)$);
        \draw[dotted, \tickthick, \cpcolor] ($(cp)+(0,0.4)$) -- (cp) node[pos=1, anchor=north, yshift=-\yclabelshift]{$\color{\cpcolor} c_f+c_{v,i}$};
        \draw[\tickthick, \cpcolor] ($(cp) + (0,\ctickheight)$) -- ($(cp) + (0,-\ctickheight)$);
        \draw[\tickthick, \cpcolor] ($(cp) + (0,\ctickheight)$) -- ($(cp) + (0,-\ctickheight)$);
        \coordinate (ectrl) at ($(cp) + (\colorshift,0)$);
        \coordinate (ectrr) at ($(end) + (-\colorshift,0) + (0,-\ycshift)$);
        \draw[\tickthick, \ectrcolor] (ectrl) -- (ectrr) node[midway, anchor=north, yshift=-\yclabelshift]{$\color{\ectrcolor} c_{\infty,i}$};
        \draw[\tickthick, \ectrcolor] ($(ectrl) + (0,\ctickheight)$) -- ($(ectrl) + (0,-\ctickheight)$);
    \end{tikzpicture}
\caption{Preventive replacement at $t_k$}
\label{subfig:rep_at_tk}
\end{subfigure}
\hfill
\begin{subfigure}{0.49\textwidth}
    \centering
    \begin{tikzpicture}
    \coordinate (0) at (0,0);
        \coordinate (end) at ($(0) + (\xdist, 0)$);
        \node[yshift=\tlabelshift cm] at (end) {$t$};
        \coordinate (tk) at ($(\xtk,0)$);
        \draw[\tickthick] ($(tk) + (0,\ytickheight)$) -- ($(tk) + 
        (0,-\ytickheight)$) node[anchor=north, yshift=-\ylabelshift]{\large $t_k$};
        \coordinate (tk1) at ($(tk) + (\tickdist,0)$);
        \draw[\tickthick] ($(tk1) + (0,\ytickheight)$) -- ($(tk1) + 
        (0,-\ytickheight)$) node[anchor=north, yshift=-\ylabelshift]{\large $t_{k+1}$};
    
        \draw[\rulcolor, very thick, domain=\domain, samples=200] 
        plot ({\x}, {\plotfactor*lognormal(\x,\mu,\sigma)+\plotyshift});
        \node[] at ($(tk1) + (\textxshift, \textyshift)$) {$\color{\rulcolor} \frulargs{i}{t}$};

        \pgfmathsetmacro\Tfoney{%
            \plotfactor*lognormal(\Tfone,\mu,\sigma)+\plotyshift %
        }
        \draw[\Tfonecolor, dashed, very thick] (\Tfone,0) -- (\Tfone,\Tfoney) 
        node[very near start, yshift=-0.55cm]{\textcolor{\Tfonecolor}{$T_{F,i}^{(1)}$}};
    
        \draw[-latex, \linethick] (0) -- (end);
    
        \coordinate (cc) at ($(\Tfone,0) + (0,-\ycshift)$);
        \draw[dotted, \tickthick, \cccolor]  ($(\Tfone,0)$) -- ($(\Tfone,0)+(0,-0.2)$);
        \draw[dotted, \tickthick, \cccolor]  ($(cc)+(0,0.4)$) -- (cc) node[pos=1, anchor=north, yshift=-\yclabelshift]{$\color{\cccolor} c_c$};
        \draw[\tickthick, \cccolor] ($(cc) + (0,\ctickheight)$) -- ($(cc) + (0,-\ctickheight)$);
        \draw[\tickthick, \cccolor] ($(cc) + (0,\ctickheight)$) -- ($(cc) + (0,-\ctickheight)$);
        \coordinate (ectrl) at ($(cc) + (\colorshift,0)$);
        \coordinate (ectrr) at ($(end) + (-\colorshift,0) + (0,-\ycshift)$);
        \draw[\tickthick, \ectrcolor] (ectrl) -- (ectrr) node[midway, anchor=north, yshift=-\yclabelshift]{$\color{\ectrcolor} c_{\infty,i}$};
        \draw[\tickthick, \ectrcolor] ($(ectrl) + (0,\ctickheight)$) -- ($(ectrl) + (0,-\ctickheight)$);
    \end{tikzpicture}
\caption{Do nothing at $t_k$, failure before $t_{k+1}$}
\label{subfig:rep_at_tk1_fail_bef_dt}
\end{subfigure}
\caption{Comparison of incurred costs preventive vs. corrective replacement of a component; taken from \cite{koutas2025leaf}.}
\label{fig:rep_graphs}
\end{figure}
}

\newcommand{\screpdt}{
\begin{figure}[H]
\def\xdist{2.5}
\def\sydist{1.75}
\def\aydist{1.1}
\def\xtextdist{1.2}
\def\ytextdist{1.0}
\def\tdist{3.5}
\def\width{0.4}
\def\radius{0.3}
\def\termradius{0.2}
\def\capsize{\normalsize}
\def\boundline{very thick}
\def\arrowthick{very thick}
\def\tshift{7}
\def\fixshift{0.4}
\def\safecolor{green!60!black}
\def\failcolor{red!80!black}
\def\condpos{0.6}
\def\pfpos{0.6}
\def\xenumshift{0.6}
\def\yenumshift{0.25}
\tikzset{circle_node/.style={circle, draw, \boundline, minimum size=\radius cm, inner sep=0pt}}
\tikzset{failed_node/.style={circle, draw, fill=black, \boundline, minimum size=\termradius cm, inner sep=0pt}}
\tikzset{replaced_node/.style={circle, draw, fill=orange!70!black, \boundline, minimum size=\radius cm, inner sep=0pt}}
\tikzset{act_node/.style={rectangle, draw, fill=cyan, \boundline, minimum size=\width cm, inner sep=0pt}}
\tikzset{cost_node/.style={rectangle, draw, fill=lime, \boundline, minimum size=1.2*\width cm, inner xsep=2pt, inner ysep=0pt}}
    \vspace{-1.2cm}
    \centering
    \begin{tikzpicture}
    \coordinate (a0cord) at (0,0);
    \node[act_node] (a0) at (a0cord) {};
    \coordinate[right=\xdist cm of a0cord, yshift=\sydist cm] (s1dncord) {};
    \node[circle_node] (s1dn) at (s1dncord) {};
    \coordinate[right=\xdist cm of a0cord, yshift=-\sydist cm] (s1prcord) {};
    \node[replaced_node] (s1pr) at (s1prcord) {};
    \coordinate[right=\xdist cm of s1dncord, yshift=\aydist cm] (dna1cord) {};
    \node[act_node] (dna1) at (dna1cord) {};
    \coordinate[right=\xdist cm of s1dncord, yshift=-\aydist cm] (dnfailedcord) {};
    \node[failed_node] (dnfailed) at (dnfailedcord) {};
    \coordinate[right=\xdist cm of s1prcord, yshift=\aydist cm] (pra1cord) {};
    \node[act_node] (pra1) at (pra1cord) {};
    \coordinate[right=\xdist cm of s1prcord, yshift=-\aydist cm] (prfailedcord) {};
    \node[failed_node] (prfailed) at (prfailedcord) {};

    \draw[\arrowthick] (a0) -- (s1dn) 
        node[midway, sloped, yshift=\tshift pt] {DN};
    \draw[\arrowthick] (a0) -- (s1pr)
        node[midway, sloped, yshift=\tshift pt] {PR};
    \draw[\arrowthick] (s1dn) -- (dna1)
        node[pos=\condpos, sloped, yshift=\tshift pt] {\textcolor{\safecolor}{Safe}}
        node[pos=\pfpos, sloped, yshift=-\tshift pt] {$1-p_F$};
    \draw[\arrowthick] (s1dn) -- (dnfailed)
        node[pos=\condpos, sloped, yshift=\tshift pt] {\textcolor{\failcolor}{Failed}}
        node[pos=\pfpos, sloped, yshift=-\tshift pt] {$p_F$};
    \draw[\arrowthick] (s1pr) -- (pra1)
        node[pos=\condpos, sloped, yshift=\tshift pt] {\textcolor{\safecolor}{Safe}}
        node[pos=\pfpos, sloped, yshift=-\tshift pt] {$\approx 1$};
    \draw[\arrowthick, dashed] (s1pr) -- (prfailed)
        node[pos=\condpos, sloped, yshift=\tshift pt] {\textcolor{\failcolor}{Failed}}
        node[pos=\pfpos, sloped, yshift=-\tshift pt] {$\approx 0$};

    \coordinate[below=\fixshift cm of s1dncord] (dnfixcord) {};
    \node[cost_node, anchor=north] (dnfix) at (dnfixcord) {$c_{\mathrm{DN}}$};
    \coordinate[below=\fixshift cm of s1prcord] (prfixcord) {};
    \node[cost_node, anchor=north] (prfix) at (prfixcord) {$c_{\mathrm{PR}}$};
    \coordinate[right=\fixshift cm of dna1cord] (dna1fixcord) {};
    \node[cost_node, anchor=west] (dna1fix) at (dna1fixcord) {$c_{\mathrm{DN,Safe}}$};
    \coordinate[right=\fixshift cm of dnfailedcord] (dnfailedfixcord) {};
    \node[cost_node, anchor=west] (dnfailedfix) at (dnfailedfixcord) {$c_{\mathrm{DN,Failed}}$};
    \coordinate[right=\fixshift cm of pra1cord] (pra1fixcord) {};
    \node[cost_node, anchor=west] (pra1fix) at (pra1fixcord) {$c_{\mathrm{PR,Safe}}$};
    \coordinate[right=\fixshift cm of prfailedcord] (prfailedfixcord) {};
    \node[cost_node, anchor=west] (prfailedfix) at (prfailedfixcord) {$c_{\mathrm{PR,Failed}}$};
    
    %
    \coordinate[above=\tdist cm of a0cord] (t0cord) {};
    \coordinate (t1cord) at ($(t0cord) + (2*\xdist, 0)$) {};
    \coordinate (t2cord) at ($(t1cord) + (3*\xdist, 0)$) {};
    \coordinate (tprevcord) at ($(t0cord) + (-3*\xdist, 0)$) {};

    \draw[decorate,decoration={brace,amplitude=10pt}] (t0cord) -- (t1cord) 
        node[midway,yshift=17pt] {\large Step $k$};
    \begin{scope}
        \clip (t1cord) circle (2.5 cm);
        \draw[decorate,decoration={brace,amplitude=10pt}] (t1cord) -- (t2cord) node[midway, xshift=-2.4cm, yshift=17pt] {\large Step $k+1$};
    \end{scope}
    \begin{scope}
        \clip (t0cord) circle (2.5 cm);
        \draw[decorate,decoration={brace,amplitude=10pt}] (tprevcord) -- (t0cord) node[midway, xshift=2.4cm, yshift=17pt] {\large Step $k-1$};
    \end{scope}

    \draw[thick, dashed, gray] ($(t0cord) + (0, 1)$) -- (a0) -- ($(a0cord) + (0, -0.5-\sydist-\aydist)$){};
        \node[gray] at ($(t0cord) + (0, 1.2)$) {$t=t_k$};
    \draw[thick, dashed, gray] ($(t1cord) + (0, 1)$) -- (dna1) -- (dnfailed) -- (pra1) -- (prfailed) --($(prfailedcord) + (0, -0.5)$) {};
        \node[gray] at ($(t1cord) + (0, 1.2)$) {$t=t_{k+1}$};
    
    \end{tikzpicture}
    \captionsetup{width=\textwidth}
    \caption{Decision tree for a deteriorating component that can be preventively replaced. Component states are represented by round nodes, where a white filling indicates no replacement, and a brown filling indicates component replacement. Maintenance decisions ($a_k$) are shown as square nodes. The action taken (PR/DN) as well as the condition of the component (Failed/Safe) with the corresponding event probabilities are included in the respective paths, where $p_F=\Pr(RUL\leq\Delta t)$. Component failure is shown with filled black circles. Action costs as well as consequence costs are shown with lime boxes where they are incurred. The Figure is adapted from \citep{koutas2025leaf}.}
    \label{fig:rep_decision_tree}
\end{figure}
}

\maketitle

\begin{abstract}
We propose predictive maintenance policies for managing replacements in multi-component series systems with economic dependence. The proposed policies 1) are guided by information available from component-specific prognostic models, 2) explicitly evaluate all available maintenance options with their potential consequences, and 3) consider the costs of the underlying renewal-reward process. Numerical investigations show that the derived policies achieve significantly lower long-running cost rates than benchmark policies, while being highly efficient and robust against overfitting.
\end{abstract}


%
%
%
%
\section{Introduction}
\label{sec:introduction}
%

Maintenance planning of multi-component systems has attracted growing interest over the past decades. System-level maintenance policies explicitly consider component interactions. These policies can be condition-based or predictive \cite[see, e.g.,][]{castanier2005condition,nguyen2015multi}. However, finding optimal maintenance policies for multi-component systems remains challenging \cite{dekker1997review,nicolai2008optimal}, which motivates this research.

\emph{Condition-based maintenance} (CBM) is a widely used approach for multi-component maintenance. It is a preventive strategy that leverages failure statistics and system-specific information from condition monitoring for decision-making \cite{ISO13306_2017,mobley2002introduction}. A prominent example of condition-based maintenance policies is the use of thresholds on observed indicators of degradation, whose exceedance triggers replacement or repair actions\cite{castanier2005condition}.  

\emph{Predictive maintenance} (PdM) utilizes a \emph{prognostic} model, which predicts the future health of the system based on the available condition-monitoring data, to inform maintenance policies \cite{mobley2002introduction,pecht2009prognostics,ran2019survey}. Such prognostic information can be obtained from a specific stochastic process assumption \cite[see, e.g.,][]{bouvard2011condition,van2013dynamic}, 
or by training a prognostic model from run-to-failure data \cite[see, e.g.,][]{dong2007segmental,li2018remaining}. PdM has become especially attractive with the increasing availability of training data and the improvement of data-driven machine learning (ML) models. The prognostic model commonly yields an estimate of the remaining useful life (RUL), i.e., the time until the system fails and ceases to perform its intended function \cite{kim2017prognostics}. To quantify the various sources of uncertainty about the current true state and future development, one typically resorts to a probability distribution over the RUL \cite{si2011remaining}.

In PdM, the prognostic model is the basis for making maintenance decisions. One defines policies as functions that take the RUL predictions as inputs and output maintenance actions. PdM policies can be grouped into two categories: black-box models or heuristic policies. Black-box decision models, such as Deep Reinforcement Learning (DRL) approaches  \cite{zhang2020deep,chen2023deep,najafi2023deep,zhang2023guided} often achieve the best performance. However, their training requires large amounts of data, which is not commonly available. In addition, the final policy is represented implicitly by the learned parameters of the neural network. Thus, the final decisions lack the transparency and interpretability required for safety-critical engineering applications. The preferred approach in practice is hence to use heuristic policies, which are simple and transparent maintenance selection rules \cite{bismut2021optimal}. A simple example of such a heuristic policy is to replace a component when the predicted RUL falls below a specified threshold. Typically, better policies can be obtained by considering the uncertainty in the RUL. One example is a policy that prescribes replacement when a specified quantile value of the RUL distribution is less than the time to the next possible replacement time \cite{kamariotis2024metric,he2023condition}.

Data-driven PdM policies are mostly formulated with little theoretical justification, relying instead on optimization to achieve good performance. This approach can lead to severe overfitting in applications with limited amounts of training data \cite{koutas2025leaf}. We aim to fill this gap by developing more robust, theory-guided heuristic PdM policies for multi-component systems. The proposed policies are based on decision trees and rely on renewal theory \cite{doob1948renewal, smith1958renewal} to directly quantify the value of a component's life extension, where the policy parameters do not require optimization to achieve good maintenance performance. By regarding the inter-arrival times as lifetimes of individual components and the associated maintenance/failure costs as the rewards, \emph{renewal-reward processes} find direct application in optimal maintenance problems. In our renewal theory framework, the ultimate goal is to find a policy that minimizes the system's long-running maintenance cost per unit time $c_{\infty}$ by accounting for the prognostic RUL predictions of the components.

%
%
%
%
\subsection{Related works on multi-component heuristics}
\label{subsec:related_work}
Renewal theory has been used since the 1960s to derive systematic preventive maintenance policies for multi-component systems \cite[e.g.,][]{radner1963opportunistic,vergin1977maintenance,okumoto1983optimum,sheu1991generalized}. These early works rely purely on statistical time-to-failure information, i.e., they do not use advanced condition information or prognostics\footnote{Some of these works do monitor the ``condition'' of the system. However, this relates to a mere check whether the system is still operational or has failed, effectively ensuring self-announcing failures. This is not to be confused with modern CBM, which yields more detailed information about the underlying state of the system.}. Overviews of these earlier works for various system dependencies and maintenance types are given in \cite{cho1991survey,dekker1997review}.

The next evolution of multi-component maintenance policies are condition-based maintenance policies, where (imperfect) information about the underlying degradation state of the components/system is obtained via continuous monitoring or scheduled inspections. Importantly, this information is restricted to the present state and does not include predictions about the future state of the system, i.e., no use of prognostics models. CBM approaches remain popular today due to ever-improving condition-monitoring systems, and are therefore the subject of ongoing research. However, because maintenance planning for multi-component systems is significantly more complex than that for single components, only a few studies \cite[e.g.,][]{castanier2005condition,zhu2010availability} have incorporated renewal theory. Reviews that discuss multi-component CBM in more detail can be found in \cite{alaswad2017review,de2020review}.

A substantial body of work has developed predictive maintenance policies for multi-component systems. 
\cite[e.g.,][]{camci2009system,van2013dynamic}. Only a few of these works incorporate renewal theory into the decision-making \cite[e.g., ][]{bouvard2011condition,nguyen2015multi,yang2024group}. By contrast to our contribution, these works include the policy's long-running cost rate, $c_\infty$, only for policy evaluation and/or as a target for parameter optimization \cite{nguyen2015multi}. Our proposed maintenance policies directly include $c_\infty$ in the cost expressions of the individual maintenance options. Furthermore, some works \cite{bouvard2011condition,yang2024group} define $c_\infty$ via the use of an individual component's predicted failure distribution rather than using the failure distribution of the \emph{population of components}, as required by  renewal theory. This point has already been raised in previous work \cite{kamariotis2024metric,koutas2025leaf}. 

%
%
%
%
\subsection{Contributions}
\label{subsec:contributions}
The objective of this paper is to obtain improved prognostics-based heuristic maintenance policies for multi-component systems. We achieve this by including theoretical considerations based on renewal theory into the derivation of the policies. Our specific contributions are summarized in the following.

We utilize the general method that we proposed in \cite{koutas2025leaf} to derive new preventive replacement heuristics for series systems with economic dependence. The derivation is performed by constructing a one-step decision tree that outlines the available maintenance options and the corresponding potential consequences. The probabilities of reaching each branch are obtained from the RUL distributions of the individual components' prognostic models. The associated branch costs are computed as the sum of the direct component costs and the costs associated with continuing the underlying renewal-reward process. We show that the policy parameters can be determined with suitable assumptions, without requiring optimization.

We compare the derived maintenance heuristics against a set of benchmark policies on a virtual RUL simulator. The proposed policies outperform optimized benchmark heuristics by as much as 35\% in terms of $c_{\infty}$ already for a 2-component system, while at the same time being highly efficient and robust against overfitting in low data regimes.
Based on the numerical investigations, we provide a detailed discussion of the strengths and weaknesses of the proposed policies and outline several possibilities for future extensions. 

The source code for this work can be found at \url{github.com/Dakout/branching_out}.

%
%
%
%
\section{Decision setting and notation}
\label{sec:decision_setting}
\subsection{Series system}
\label{subsec:selected_system}
We address the maintenance optimization of multi-component \emph{series systems}. Series system structures are common in numerous engineering applications, including power lines, oil/water/gas pipelines, manufacturing assembly lines, drivetrains, and aircraft flight control systems. They are of critical importance for reliability engineering and safety analysis, and thus of high interest to maintenance research \cite[see, e.g.,][]{chalabi2016optimisation}.

We classify the interaction of components in the system into three main categories, namely \emph{economic}, \emph{structural} and \emph{stochastic} dependence. Economic dependence offers cost-saving opportunities in case of joint maintenance of several components (or the opposite). If maintenance of a failed component implies maintenance of another component (i.e., if these components structurally form a part), this corresponds to structural dependency. Lastly, stochastic/probabilistic dependence signifies correlated lifetime distributions of components, e.g., due to common-cause failures \cite{thomas1986survey,dekker1997review}. In this paper, we consider only positive economic dependence, i.e., grouping maintenance of several components saves costs compared to individual maintenance. The applicability of our final policies to other dependencies is discussed in \Cref{sec:dicussion}.

%
%
%
%
\subsection{Preventive replacement of a series system}
\label{subsec:preventive_replacement}
We focus on preventive replacement as the maintenance action to decide on. We consider a series system of $M$ components, each equipped with a custom prognostic model, i.e., we assume that a probabilistic RUL prediction is available for each component. These predictions are fed into a heuristic policy, which outputs a joint action vector, assigning each component either a ``do nothing'' (DN) or ``preventively replace'' (PR) decision. The whole workflow is illustrated in \Cref{fig:decision_setting}. 

Let $t$ denote the time since the system has been installed. The set of discrete decision times for preventive replacement is denoted with $\{t_k=k\cdot \Delta t,k=1,2,\ldots\}$. A component $i$ can be replaced via two possible replacement types, namely (1) preventive replacement with preventive cost as a sum of fixed, $c_f$, and variable costs, $c_{v,i}$, and (2) corrective replacement upon failure with corrective cost $c_c$. We assume  $c_c \geq c_f+\sum_{i=1}^M c_{v,i}$ throughout this work, ensuring that corrective replacement is always more expensive than preventive replacement. Additional assumptions are that failures are self-announcing, replacements are perfect and are performed immediately and instantaneously upon failure of a component. Here, only components that have failed or for which a PR action has been chosen are replaced; other components are \emph{not} replaced, but continue from their current damaged state.\footnote{Another variant explored in some works \cite{murthy1985study,van2013dynamic} is that failure of a single component leads to destruction of the whole system, i.e., renewal of all components. Adjusting the proposed policies to this variant is straightforward and computationally more favourable, as the renewal cycles of the system can be more easily obtained.}
Moreover, the successive lifetimes of all components are independent and identically distributed with finite mean. Therefore, renewal-reward theory can be applied, where a component-wise renewal cycle is defined as the time interval between two successive replacements of a component. Finally, future costs are not discounted, as component lifetimes are assumed to be relatively short.

The failure time of a component $i$ is denoted with $T_{F,i}$. If reference to a specific sample is needed, we explicitly refer to it via a sample superscript $j$, e.g., $T_{F,i}^{(j)}$; otherwise, the sample-based notation is omitted wherever possible to improve readability. We further use the definition of the remaining useful life $RUL_i(t)=T_{F,i} - t$, with the shorthand $RUL_{i,k}\coloneqq RUL_i(t_k)$ for the discrete times $t_k$. In the considered setting, one obtains from each of the $M$ component-specific prognostic models at each $t_k$ a predicted RUL probability density function (PDF) denoted with $\frulargs{i,k}{\cdot}$. These predictions are based on all monitoring information available up to time $t_k$, and are subsequently used as input to the decision policy, as illustrated in \Cref{fig:decision_setting}.

Throughout this work, unless explicitly noted otherwise, the latest RUL-PDF of component $i$ is used as the reference distribution to compute cumulative probabilities and expectations, specifically, the probability of the RUL being in the interval $(a,b]$, \Cref{eq:prob_RUL_betweeen_a_b}, and the expected RUL given that it lies in the interval $(a,b]$, \Cref{eq:exp_RUL_betweeen_a_b}: 
\begin{align}
    \label{eq:prob_RUL_betweeen_a_b}
    \Pr(a<RUL_{i,k}\leq b) &\coloneqq \int_a^b \frulargs{i,k}{r} \dr \\
    \EE{RUL_{i,k} \mid a<RUL_{i,k}\leq b} & \coloneqq \frac{1}{\Pr(a<RUL_{i,k}\leq b)} \int_a^b r \cdot \frulargs{i,k}{r} \dr.
    \label{eq:exp_RUL_betweeen_a_b}
\end{align}

\workflow

%
%
%
%
\subsection{Policy evaluation}
\label{subsec:performance_evaluation}
In predictive maintenance, a policy $\Pi$ takes the prognostics information as an input, and provides a maintenance action as an output. A more rigorous definition of policies is provided in \Cref{sec:pdm_policies}. To numerically evaluate the average performance of $\Pi$, the policy is repeatedly applied to a population of systems that follow a certain time-to-failure (TTF) distribution. Upon failure/replacement, the system is assumed to be as good as new; therefore, we can use the elementary renewal-reward theory to compute the policy's long-running maintenance cost per unit time as \cite{grimmett2020probability}
\begin{equation}
\label{eq:ectr_general}
    c_{\infty} \mid \Pi = \lim_{t \rightarrow \infty} \frac{C(t) \mid \Pi}{t} = \frac{\EE{C \mid \Pi}}{\EE{T \mid \Pi}}.
\end{equation}
The expectations $\mathbb{E}[\cdot]$ in the numerator and denominator of \Cref{eq:ectr_general} are both taken with respect to the resulting renewal cycles of the system, and $C \mid \Pi$ and $T \mid \Pi$ are the life cycle costs and times that depend on the specific policy. For the sake of readability, we omit the explicit dependence on $\Pi$ in the following wherever possible. For our preventive replacement setting, $T$ denotes either the point in time when the system is preventively replaced or fails, whichever happens first, and $C$ denotes the corresponding preventive/corrective replacement costs. The overarching goal is to find policies which minimize $c_{\infty}$. In the following, explicit dependency of the performance on $\Pi$ is also omitted for the sake of readability. 

Renewal cycles of the series system only start upon the simultaneous failure or preventive replacement of all components. Since we anticipate this to be a rare event (especially with increasing number of components), we instead opt for a computation of $c_{\infty}$ based on the individual component's long-running cost rates $c_{\infty,i}$.

\begin{figure}[H]
    \centering
    \includegraphics[width=0.85\linewidth]{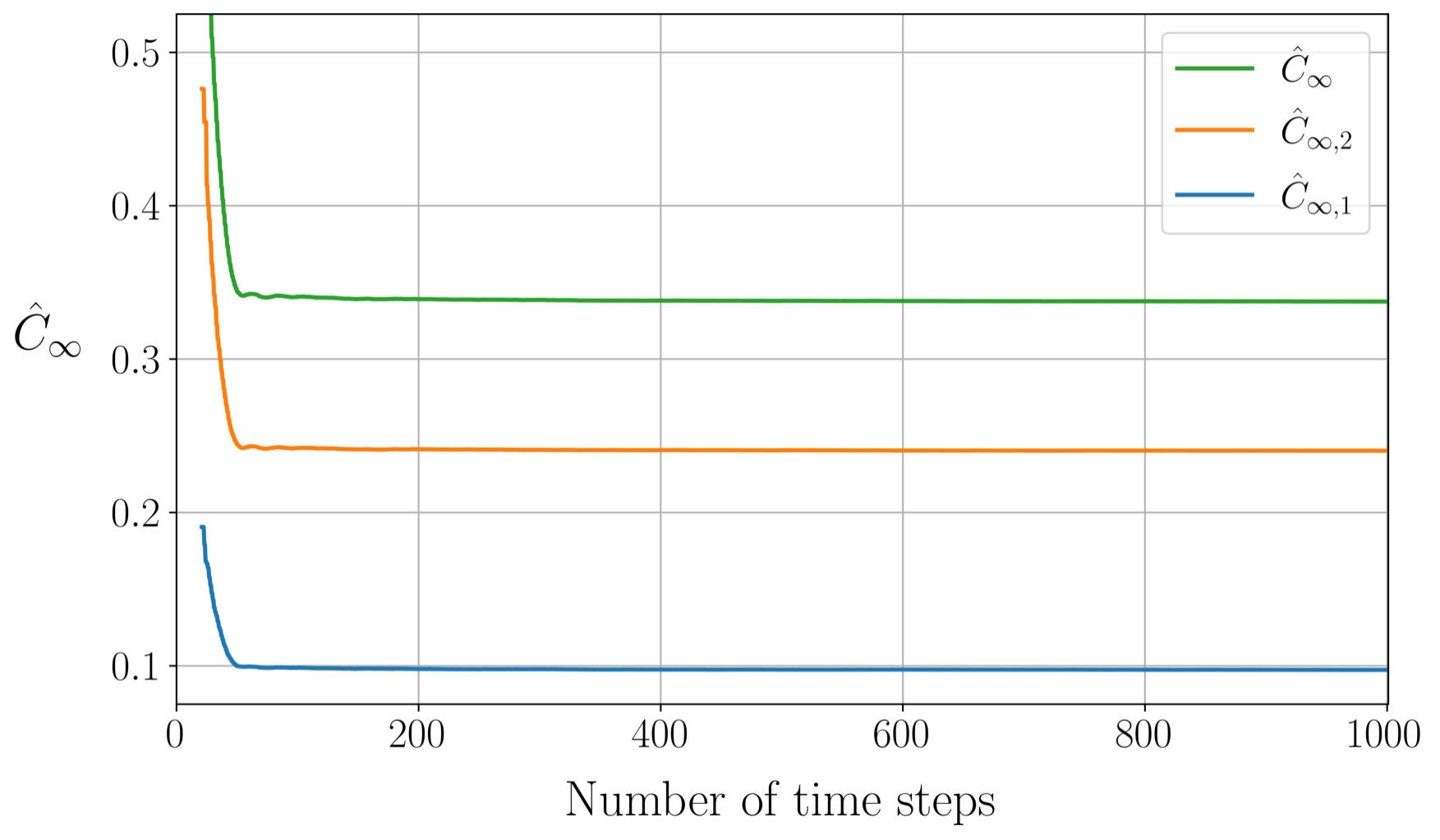}
    \caption{Long-running maintenance costs per unit time plotted over the number of time steps of length $\Delta t$. Shown are the long-running cost rates of two individual components as well as the overall series system. The $rh1^*$ policy (introduced in \Cref{subsubsec:probability_thresholds}) is applied to the case $c_f=10$, $c_{v,1}=10$, $c_{v,2}=40$. The results are obtained with Monte Carlo simulation with a batch size of $B=2,000$ until a maximum horizon of $H=10,000$ (i.e., $1,000$ steps).}
    \label{fig:convergence_plot}
\end{figure}

To numerically evaluate a heuristic policy for an $M$-component series system, we simulate a batch of trajectories with size $B$, each consisting of $M$ components until a horizon $H$. At the end of the horizon, every single component $i$ and trajectory $b$ will have a certain number of replacement events $N_{b,i}$, with corresponding life cycle costs and times $C_{b,i}^{(j)} ~\&~T_{b,i}^{(j)},~j=1,\ldots N_{b,i}$. Generally, upon reaching the horizon, each component has a unique residual lifetime without failure $R_{b,i}=H-\sum_{j=1}^{N_{b,i}}T_{b,i}^{(j)}$, which introduces an approximation. The resulting sample-based estimate of the long-running cost ratio is
\begin{equation}
    \label{eq:ectr_sys}
    \hat{C}_{\infty} \approx \frac{\sum_{i=1}^M \frac{
    \frac{1}{B}\sum_{b=1}^B 
    \frac{1}{N_{b,i}}\sum_{j=1}^{N_{b,i}}
    C_{b,i}^{(j)}}{
    \frac{1}{B}\sum_{b=1}^B
    \frac{1}{N_{b,i}}\sum_{j=1}^{N_{b,i}}
    T_{b,i}^{(j)}} \cdot H}{H} = 
    \sum_{i=1}^M \frac{
    \sum_{b=1}^B\sum_{j=1}^{N_{b,i}}C_{b,i}^{(j)}}{\sum_{b=1}^B\sum_{j=1}^{N_{b,i}}T_{b,i}^{(j)}}
    = \sum_{i=1}^M \hat{C}_{\infty,i},
\end{equation}
where $\hat{C}_{\infty,i}$ denotes the MC-approximation of component $i$'s long-running cost ratio.
According to \Cref{eq:ectr_sys}, the system's long-running cost rate estimate is obtained as a summation over the individual components' cost rates, while disregarding the residual lifetimes. We check the convergence of $\hat{C}_{\infty}$ with increasing horizon $H$, by recording the events and plotting the resulting cost rate over time. An example for such a convergence plot is given in \Cref{fig:convergence_plot}.
We perform this manual convergence verification for all investigated policies and system configurations. $B$ and $H$ are chosen accordingly to satisfy convergence.
For uncertainty quantification, we compute the variance of the maintenance performance based on the approximation in \Cref{eq:ectr_sys} and the independence assumption between the individual components\footnote{Due to the joint maintenance of components, there is correlation between the component lifetimes. However, since we choose $H$ and $B$ such that we have vanishing variance, and the computation of $\mathbb{V}\mathrm{ar}[{\hat{C}_{\infty,i}}]$ is already approximate, we deem this approximation as not critical.}:
\begin{equation}
    \label{eq:var_ectr}
    \Var{\hat{C}_{\infty}} \approx \sum_{i=1}^M \Var{\hat{C}_{\infty,i}}.
\end{equation}
$\Var{\hat{C}_{\infty,i}}$ involves computing the variance of a ratio of random variables, which does not have an analytical expression. Instead, given a sample of size $n$ for component $i$, one can use the first-order approximation \cite{van2000mean}:
\begin{equation}
    \label{eq:ectr_first_order}
    \Var{\hat{C}_{\infty,i}} \approx \frac{1}{n} \left[
    \frac{\hVar{C_i}}{\hEE{T_i}^2} + \frac{\hEE{C_i}^2 \cdot \hVar{T_i}}{\hEE{T_i}^4}
    - 2 \cdot \frac{\hEE{C_i} \cdot \widehat{\mathbb{C}\mathrm{ov}}\left[C_i, T_i\right]}{\hEE{T_i}^3}
    \right],
\end{equation}
where $\hEE{X}$, $\hVar{X}$, and $\widehat{\mathbb{C}\mathrm{ov}}\left[X,Y\right]$ denote unbiased estimators of $\EE{X}$, $\Var{X}$, and $\mathbb{C}\mathrm{ov}\left[X,Y\right]$, respectively.

%
%
%
%
\section{PdM policies}
\label{sec:pdm_policies}

%
%
%
%
\subsection{Definition of heuristic policies}
\label{subsec:definition_of_heuristic_policies}
In a general maintenance setting, individual decisions have an impact on the state of the underlying system as well as on the effect of future decision alternatives. Consequently, quantifying individual decision effects requires consideration of the full sequence of decisions throughout the structure's life. Due to uncertainty in past, present, and future system states as well as in the available information, this problem is formulated as a stochastic sequential decision problem (SDP)\cite{raiffa2000applied,kochenderfer2015decision,bismut2021optimal}. The corresponding full decision tree grows polynomially or even (super-)exponentially in the set of the number of system states, observations, actions, components and considered time horizon, which is referred to as \emph{curse of dimensionality} \cite{bellman1966dynamic} and \emph{curse of history} \cite{pineau2006anytime}.

The decision on which action to take at time step $t_k$ is determined by a function called \emph{policy} $\pi_{t_k}$. Since this paper generalizes the single-component heuristic policies introduced in \cite{koutas2025leaf} to a multi-component series system, we closely follow its notation, which is based on partially observable Markov decision processes (POMDPs) and reinforcement learning (RL) \cite{kaelbling1998planning,sutton1998reinforcement}. For the investigated multi-component system, we adopt a centralized control perspective. Thus, the policy definitions follow the notation of standard POMDPs rather than that of decentralized POMDPs (Dec-POMDPs) \cite{bernstein2002complexity}.

Given a set of system states $s \in \mathcal{S}$ and available actions $a \in \mathcal{A}$, a deterministic policy is a mapping $\pi: \mathcal{S}\rightarrow \mathcal{A}$, assigning an action $a=\pi_{t_k}(s)$ to each state $s$ \cite{kaelbling1996reinforcement}.\footnote{Note that this notation already assumes that the state $s$ is Markovian \cite{ethier2009markov}; otherwise, the policy depends on the full history of states and actions: $\pi_{t_k}(s_{0:t_k},a_{0:t_{k-1}})$ \cite{kochenderfer2015decision} Furthermore, policies are generally defined as probability distribution over $\mathcal{A}$, denoted as $\pi_{t_k}(a\mid s)=\Pr\left(a \mid S_{t_k}=s\right)$. However, stochastic policies are uncommon in PHM.}. Note that for our multi-component system, $s$ and $a$ can be factored as state and action vectors over the individual components: $s=[s_1,~\ldots,~s_M]$ and $a=[a_1,~\ldots,~a_M]$.

In practical cases, the true underlying state of the system is typically unknown. In POMDPs, this state uncertainty is represented by \emph{belief states} $b \in \mathcal{B}$, which are probability distributions over $\mathcal{S}$. The belief $b_{t}(s)=\Pr(S_{t}=s)$ denotes the probability of the system being in state $s$ at time $t$, with $0 \leq b_{t_k}(s) \leq 1, ~\forall s \in \mathcal{S}$ and $\sum_{s \in \mathcal{S}} b_{t_k}(s) = 1$. In POMDPs, the policy is then reformulated to be a function of the belief state: $\pi_{t_k}(b)$ \cite{kaelbling1998planning}, where $b$ is updated with Bayes' rule, taking into account the underlying system dynamics. 

For our predictive maintenance setting, we define the individual components' RULs as the state of the system, and the RUL-PDFs obtained by the prognostic models as the belief state. Also this system state and belief state can be represented as vectors over the component RULs and component RUL-PDFs, respectively:
\begin{align}
    \label{eq:state_as_rul}
    \bm{S}_{t_k}&:=[RUL_{1,k},~ \ldots,~ RUL_{M,k}]\\
    \bm{b}_{t_k}&:=[f_{RUL_{1,k}},~ \ldots,~f_{RUL_{M,k}}].
    \label{eq_belief_as_rul_pdf}
\end{align}
This set of RUL-PDFs is the input to all considered policies, which is why it is omitted in the following, and we simply write $\pi_{t_k}=\pi_{t_k}(b)$. 

A \emph{strategy} $ \Pi=\{\pi_{t_1}, \pi_{t_2}, \dots\} \in \mathcal{P}$ denotes the set of policies across all time steps, where $\mathcal{P}$ is the space of all admissible strategies. In the context of renewal theory, the globally optimal strategy $\Pi^*$ is the one minimizing the long-running maintenance cost per unit time:
\begin{equation}
    \label{eq:optimal_strategy}
    \Pi^*= \underset{\Pi \in \mathcal{P}}{\arg \min }~ \left[c_{\infty} \mid \Pi\right].
\end{equation}
Finding $\Pi^*$ in the entire strategy space is computationally unfeasible. Hence, to facilitate optimization, we reduce $\mathcal{P}$ by using a \emph{heuristic}, which is a parametrized strategy $\Pi_{\mathbf{w}}$ with parameters $\mathbf{w}=\{w_1, w_2, \dots, w_h\}$. Instead of \Cref{eq:optimal_strategy}, we then solve the following optimization problem instead \cite{bismut2021optimal}:
\begin{equation}
    \label{eq:optimal_heuristic_strategy}
    \mathbf{w}^*= \underset{\mathbf{w}}{\arg \min }~ \left[c_{\infty} \mid \Pi_{\mathbf{w}}\right], \qquad \Pi_{\mathbf{w}}\in \mathcal{P}_{\mathbf{w}}.
\end{equation}
Since $\mathcal{P}_{\mathbf{w}}\subseteq \mathcal{P}$, the solution of \Cref{eq:optimal_heuristic_strategy} does generally not coincide with the globally optimal solution of \Cref{eq:optimal_strategy}. 

We consider an infinite time horizon, under the assumption that components are continuously replaced. Although in practice the actual system lifetime is finite, the infinite-horizon assumption provides a good approximation when component lifetimes are significantly shorter than the lifetime of the overall system \cite{hinderer1999approximate}. Under an infinite-horizon formulation, the optimal strategies are \emph{stationary}, meaning the policies are identical at all time steps. Consequently, the notation can be simplified by omitting the time index $\Pi = \pi_{t_1} = \pi_{t_2} = \dots$ \cite{bismut2021optimal}. Accordingly, the terms \emph{policy} and \emph{strategy} are used interchangeably throughout this work.

%
%
%
%
\subsection{Our proposed multi-component doa heuristics}
\label{subsec:proposed_doa_heuristic}
%

%
%
%
%
\subsubsection{Derivation of hybrid PdM policies}
\label{subsubsec:derivation_hybrid_policies}
Heuristic policies can be viewed as a way to simplify the underlying complexity of the full decision tree. In this paper, we follow our general method proposed in \cite{koutas2025leaf} to derive predictive replacement policies for a multi-component series system. With this work, we thus extend the application of our heuristics from single components to multi-component systems. In short, at each time step, we construct a one-step simplified decision tree via the procedure described in \Cref{fig:general_method_overview}. 
\begin{figure}[H]
\centering
\begin{minipage}{0.56\textwidth}
\begin{tcolorbox}[title=General method for deriving PdM policies, colback=white, colframe=black]
    \begin{enumerate}[leftmargin=*, itemsep=5pt]
        \item[i.] Construct a simplified decision tree based on the available maintenance options and their potential consequences.
        \item[ii.] Assign probabilities to each branch based on the RUL-PDFs obtained by the prognostic models.
        \item[iii.] Assign to each branch the expected costs resulting from
        \begin{enumerate}[topsep=-5pt, itemsep=0pt]
            \item[a)] the monitored components,
            \item[b)] the underlying renewal-reward processes.
        \end{enumerate}
        \item[iv.] Simplify the cost expressions and derive the policy that chooses the action vector with minimal expected cost according to the simplified decision tree.
        \item[v.] Perform initialization of the selected parameters of the derived policy.
    \end{enumerate}
\end{tcolorbox}
\end{minipage}
\hfill
\begin{minipage}{0.42\textwidth}
\centering
    \begin{tikzpicture}[
        node distance=0.5cm,
        box/.style={rectangle, draw, rounded corners, align=center, minimum width=3.4cm, minimum height=0.9cm},
        arrow/.style={-latex, thick}
    ]
    \node[box] (tree) {Decision tree\\construction};
    \node[box, below=of tree] (prob) {Probability\\assignment};
    \node[box, below=of prob] (cost) {Cost\\assignment};
    \node[box, below=of cost] (policy) {Expected cost\\minimization};
    \node[box, below=of policy] (opt) 
    {Parameter\\initialization};
    
    \draw[arrow] (tree) -- (prob);
    \draw[arrow] (prob) -- (cost);
    \draw[arrow] (cost) -- (policy);
    \draw[arrow] (policy) -- (opt);

    \end{tikzpicture}
\end{minipage}
\caption{Detailed description of the general method for deriving PdM policies (left), with a graphical summary of the main stages (right). The figure is adapted from \cite{koutas2025leaf}.}
\label{fig:general_method_overview}
\end{figure}

To illustrate the procedure stage-by-stage, we show a one-step lookahead decision tree for a single component with the DN/PR decision options in \Cref{fig:rep_decision_tree}.

Stage i reflects the limited available maintenance options and the corresponding consequences in real-life maintenance settings. As an example, the potential inspection times of an aircraft are restricted by the aircraft's operation time and the availability of the crew. In this paper, this restriction is represented by the decision interval $\Delta t$. In \Cref{fig:rep_decision_tree}, the maintenance options are represented by the action branches (PR/DN) as well as the discrete times when these actions can be taken ($t_k$ vs. $t_{k+1}$). Regarding the consequences, the component can either fail or survive until the next decision opportunity upon a taken maintenance action.

In stage ii, the latest RUL distribution, provided by a prognostic model, is used to assign probabilities to each branch in the decision tree. Naturally, this RUL distribution changes from time step to time step, as new condition-monitoring data becomes available. 

We further assume that the decision intervals are small compared to the component lifetimes. Furthermore, we assume that a new component, which is installed upon failure/replacement of an old component, is unlikely to fail until the next decision opportunity; thus, its probability of failure is approximated as 0. This is represented by the dashed branch in \Cref{fig:rep_decision_tree}. Note that the hazard rate of a typical engineering component/system follows the bathtub curve. Hence, this assumption of vanishing failure probability within the first decision time step (which holds for our investigated dataset introduced in \Cref{subsec:data_set}) may not hold. However, the decision tree can be easily amended to include this failure probability.

In stage iii, costs are assigned to each branch. These costs concern both the direct costs following a maintenance action ($c_\mathrm{DN}$, $c_{\mathrm{PR}}$), as well as the costs that are incurred from each potential consequence following a certain maintenance action ($c_\mathrm{DN,Safe}$, $c_\mathrm{DN,Failed}$, $c_{\mathrm{PR, Safe}}$, $c_\mathrm{PR,Failed}$). These consequence costs include potential failure costs of the component as well as the expected benefits associated with keeping the component in service longer. Here, renewal theory is used to quantify these life extension benefits. Specifically, we use the component's long-running average cost per unit time of the respective maintenance policy ($c_{\infty,i}$). In \Cref{fig:rep_decision_tree}, the cost expressions are deliberately kept abstract, a more detailed breakdown of the individual costs is provided in the subsequent sections.

After the assignment of branch probabilities and costs, stage iv selects the action with minimal expected costs. 

Stage v concludes the policy derivation with an initialization of the policy's parameter ($c_{\infty,i}$). Note that ``parameter optimization'' strictly refers to solving  \Cref{eq:optimal_heuristic_strategy}. Here, the optimal values for the policy parameters are found by repeated application of the respective policy to a set of components, i.e., based on policy performance information. ``Parameter initialization'', on the other hand, refers to some kind of selection rule of the parameters, which is \emph{not} based on policy performance. Note that the initialization can still be data-driven and may involve an optimization stage, just with a different target function. We comparatively discuss parameter initialization/optimization for both the benchmark policies as well as the proposed $doa$ policies in later sections.

We refer to policies derived with the stages above as \emph{discrete option assessment} ($doa$) policies. In this work, the core simplification of the $doa$ policies consists of reducing the full sequential decision tree into a one-step lookahead decision tree, with the corresponding binary decision of \emph{replace now or replace later}. This one-step formulation was shown to be effective for one-component systems \cite{koutas2025leaf}. Since $doa$ policies combine an online one-step decision tree with a renewal theory-based offline value estimation of the future costs, they can be classified as \emph{hybrid} policies \cite{kochenderfer2022algorithms}. 

\screpdt
%

%
%
%
%
\newpage
\subsubsection{Decision tree construction for series systems}
\label{subsubsec:decision_tree_construction}
Following the decision tree for single components, we now construct the decision tree for a multi-component series system with economic dependence for the preventive replacement setting (stage i). Since two actions (PR/DN) are available for each component per time step, the action branches are action vectors (central control) and the number of action branches is growing as $2^M$. Following a certain action vector, each component can either fail or survive during $\Delta t$, leading to $2^M$ branches per action branch. Thus, in total, the one-step decision tree has a maximum of $2^{2M}$ branches. This decision tree is exemplarily depicted for 2 components in \Cref{fig:dt}.

Due to the stochastic independence of the components, each component is equipped with its prognostic model, yielding its marginal RUL distribution. The probabilities of reaching an individual branch (stage ii of the proposed method in \Cref{subsubsec:derivation_hybrid_policies}) are then computed with the RUL-PDFs of the individual components. Again, we assume that a component is unlikely to fail within the first decision interval, i.e., the probability of failure of a replaced component is approximated as 0. 

Regarding the expected costs of each branch (stage iii), the direct action costs for each action vector are the preventive replacement costs ($c_f~\&~c_{v,i})$, with the economic dependence taken into account. The consequence costs are the sum of potential component failure costs, $c_{F,i}$ (which are based on the corrective replacement cost $c_c$), as well as the costs assigned to the intact components, $c_I^j$. Both $c_{F,i}$ and $c_I^j$ include the costs of the underlying renewal-reward process (which are based on the cost rates $c_{\infty,i}$). The details of the individual expressions are outlined in \Cref{subsubsec:general_structure,subsubsec:doa1,subsubsec:doa2}. 

Therefore, to compute all final costs, one needs the cost structure of the underlying system ($c_f, ~c_{v,i},~c_c$) as well as the components' cost rates ($c_{\infty,i}$). The cost structure is given by the problem setting; hence, the individual cost rates are the set of parameters of the $doa$ heuristic, i.e., $\textbf{w}\coloneqq \bm{c_{\infty}}=[c_{\infty,1}, ~\ldots,~c_{\infty,M}]$. In this work, we only consider parameter initialization (stage v), and not optimization. The corresponding procedure is described in the following section.

%
%
%
%
\subsubsection{Parameter initialization}
\label{subsubsec:parameter_initialization}
\citet{kamariotis2024metric} present three different initialization methods to determine suitable values for the free parameters, denoted as $c_{\infty,i}^{0}$. We use the first initialization method, which considers maintenance without monitoring (NM). Here, renewal theory can be used to find the optimal deterministic replacement time (age-based replacement). The objective function for a given replacement time $\trep$ is formulated as:
\begin{equation}
\label{eq:obj_func_ectr_init_nm}
    g^{NM}_i(\trep) = \frac{p_{F,i}\cdot c_c + (1-p_{F_i})\cdot (c_f+c_{v,i})}{\int_0^{\trep}t\cdot f_{\Tf,i}(t) \dt + (1-p_{F,i})\cdot \trep}, \qquad \text{with~} p_{F,i}=\Pro{T_{F,i}\leq\trep}
\end{equation}
and the corresponding $c_{\infty,i}^0$ is found by solving the optimization problem 
\begin{equation}
    \label{eq:ectr_init}
    c_{\infty,i}^0 = \min_{\trep} g^{NM}_i(\trep), \quad \trep \in (0, \infty).
\end{equation}
The optimization can either be performed for analytical TTF distributions, $f_{T_F,i}(t)$, or for a given set of TTF samples in which case one obtains a numerical approximation $\hat{C}_{\infty,i}^0\approx c_{\infty,i}^0$. 

As can be seen in \Cref{eq:obj_func_ectr_init_nm}, the optimization is univariate, i.e., one can perform $M$ independent individual optimizations to obtain $\bm{c_{\infty}}^0$, which is easier than a joint optimization over all parameters. The reason for the univariate optimization is that, generally, each component has a custom TTF distribution, cost structure and starting age. Thus, the deterministic replacement times of each component are unique and thus there are no system effects to take into consideration.

Note that although finding $c_{\infty,i}^{0}$ with the presented method involves an optimization step, the target function is different from the final policy's performance described in \Cref{eq:optimal_heuristic_strategy}. Hence, the method presented herein is still referred to as parameter initialization instead of optimization.

\decisiontree
%
%
%
%
%
\subsubsection{General structure}
\label{subsubsec:general_structure}
In this section, we describe the general structure and the detailed procedure of obtaining the final heuristic $doa$ policy. For that, we find a general formulation of the consequences (costs and system states) that arise from different actions, i.e., the first and second layer of branches in the decision tree in \Cref{fig:dt}. We then get the final policy as the selection of the action branch with minimal expected cost. Since we show the derivation of the final policy for a specific example time step $t_k$, we drop the subindex $k$ wherever possible for ease of notation.

At decision time $t_k$, if one decides to preventively replace component $i$, this is indicated with the binary variable $z_i=1$; otherwise $z_i=0$. The full action vector is hence given by $\bm{z}=[z_1,~\ldots,~z_M]$, where $\bm{z} \in \mathcal{Z}\coloneqq\{0,1\}^M$. For a given action vector $\bm{z}$, the total consequences are the immediate preventive replacement costs as well as expected future costs due to 1) failure of components and 2) survival of components. For ease of notation, we omit the explicit dependence of all following costs on $\bm{z}$.

The immediate preventive replacement costs, $c_P$, are given as
\begin{equation}
\label{eq:cP}
    c_{P} \coloneqq c_P \mid \bm{z} 
    = \left(\sum_{i=1}^M z_i \cdot (c_f + c_{v,i})\right) - \max\left(0,~-1+\sum_{i=1}^M z_i\right) \cdot c_f.
\end{equation}
The second part of the above equation represents the economic dependence of the components, where for multiple preventive replacements $c_f$ still only has to be payed once. \Cref{eq:cP} has been used to compute the costs that are displayed below the white, light orange, and brown circle nodes in \Cref{fig:dt}.

Next, to obtain the expected failure and survival costs for a certain action vector, one has to consider all possible system states that could arise from $\bm{z}$, with their corresponding probabilities. The system state $s^j$ depends on the individual component vector states $\bm{x}^j = [x_1^j,~ \ldots,~ x_M^j]$ as well as the underlying system logic, where 
\begin{equation}
    s^j = 
    \begin{cases}
    1 & \text{if system failed} \\
    0 & \text{if system is working}    
    \end{cases}
    \qquad \qquad \& \qquad \qquad 
    x_i^j = 
    \begin{cases}
    1 & \text{if component }i \text{ failed} \\
    0 & \text{if component }i \text{ is working}   
    \end{cases}
\end{equation}
For an $M$-component series system, $j=1,\ldots,2^M$ and $s^j=1,~\forall j>1$. The collection of all possible system states depending on all possible combinations of component states can be represented with the truth table shown in \Cref{tab:system_truth_table}.\footnote{For a series system, the truth table is trivial. We still show the truth table because the implementation relies on it, to facilitate the adaptation of the proposed $doa$ policies to other multi-component systems.}
\begin{table}[H]
    \renewcommand{\arraystretch}{1.2}
    \setlength{\abovecaptionskip}{8pt}
    \centering
    \begin{tabular}{c|cccc|c}
    \hline
    $j$ & $x_1^j$ & $x_2^j$ & \ldots & $x_M^j$ & $s^j$ \\
    \hline
    1 & 0 & 0 & 0 & 0 & 0 \\
    2 & 1 & 0 & 0 & 0 & 1 \\
    3 & 0 & 1 & 0 & 0 & 1 \\
    4 & 1 & 1 & 0 & 0 & 1 \\
    \vdots & \vdots & \vdots & \vdots & \vdots & \vdots \\
    $2^M$ & 1 & 1 & 1 & 1 & 1 \\
    \hline
    \end{tabular}
    \caption{Truth table for an $M$-component series system.}
    \label{tab:system_truth_table}
\end{table}
To compute the probabilities for each component vector state, we first introduce the short-hand notation for the one-step failure probability of a component as $p_{F,i}=\Pr(RUL_i \leq \Delta t)$. 
For a system with $M$ statistically independent components, the probability of reaching a specific component vector state, $\bm{x}^j$, is
\begin{equation}
    \label{eq:pr_xj}
    \Pr(\bm{x}^j) = 
    \prod_{i=1}^M 
        x_i^j \cdot p_{F,i} + (1-x_i^j)\cdot(1-p_{F,i}).  
\end{equation}
Since we are interested in the future failure and survival costs for a given action vector $\bm{z}$, we have to incorporate the change of failure probabilities due to preventive replacement. In \Cref{subsubsec:decision_tree_construction}, we assumed that the lifetime of components is much larger than the decision interval $\Delta t$, which leads to $\Pr(T_{F,i}\leq \Delta t)\approx 0$. The effect of this assumption is that components cannot fail after they have been replaced, i.e., the corresponding component vector states are assigned an occurrence probability of 0. This constraint can be included in \Cref{eq:pr_xj} with the help of $\bm{z}$ as
\begin{equation}
    \label{eq:pr_xj_given_z}
    p^j \coloneqq\Pr(\bm{x}^j \mid \bm{z}) = 
    \prod_{i=1}^M (1-z_i) \cdot \left[ 
        x_i^j \cdot p_{F,i} + (1-x_i^j)\cdot(1-p_{F,i})
    \right] + z_i\cdot(1-x_i^{j}).
\end{equation}
The expression given in \Cref{eq:pr_xj_given_z} has been used to obtain the probabilities for the consequence branches in \Cref{fig:dt}.

Corresponding to the branch probabilities, the individual branch consequences are the sum of the expected failure costs $c_{F,i}$ and the cost attributed to the survived components $c_I^j$. Both of these incorporate the future costs associated with the continuation of the underlying renewal-reward process. 

We start with $c_{F,i}$. To show how the renewal-reward process costs are incorporated into $c_{F,i}$, we first consider again a single component. Upon replacement (preventive or corrective) of a component, a new component is installed. Since there is no information on this new component, its maintenance costs are represented with the average cost rate $c_{\infty,i}$. Since preventive replacement and failure occur at different times, this ``switch'' to the average process occurs at different times, ultimately leading to different future cost estimations. \Cref{fig:rep_graphs} illustrates this for two cases for a single component: 1) PR at $t_k$ and 2) DN at $t_k$, followed by failure of the component before $t_{k+1}$. 

\repgraphs

One can observe in \Cref{fig:rep_graphs} that failure of a component incurs higher direct costs ($c_c>c_f+c_{v,i}$), but the time frame for paying the long-running costs, $c_{\infty,i}$, is reduced by $RUL_i$ in case of $RUL_i\leq \Delta t$. Since for $t>T_{F,i}^{(1)}$ the long-running costs are identical for preventive and corrective replacement, we can cancel them. Taking the expectation over $T_{F,i}$, the expected costs for corrective replacement relative to preventive replacement are
\begin{equation}
    \label{eq:expected_failure_costs_comp_i}
    c_{F,i} = c_{c} - c_{\infty,i}\cdot\EE{RUL_i \mid RUL_i \leq \Delta t}. 
\end{equation}
The subtrahend in \Cref{eq:expected_failure_costs_comp_i} represents the value of keeping component $i$ in service longer, i.e., the value of life extension \cite{koutas2025leaf}.

Including the effects of maintenance actions (i.e., upon PR, a component cannot fail until the next decision time step), the expected failure costs for a certain component vector state are obtained as
\begin{equation}
    \label{eq:cF}
    c_F^j \coloneqq\EE{C_{F}^j \mid \mathbf{z}} =
    \sum_{i=1}^M (1-z_i) \cdot x_i^j \cdot c_{F,i}.
\end{equation}
The estimate of the future costs of the intact components for a certain component vector state, $C_I^j$, depends on the specific $doa$ version. For now, we leave it with the general formulation
\begin{equation}
    \label{eq:cIj_general}
    c_I^j \coloneqq \EE{C_I^j \mid \bm{z}}
\end{equation}
and introduce two candidate formulations for $C_I^j$ in \Cref{subsubsec:doa1,subsubsec:doa2}.

Given the expected preventive, failure, and survival costs in \Cref{eq:cP,eq:cF,eq:cIj_general}, respectively, the expected total cost for a given action vector is obtained as
\begin{equation}
    \label{eq:cTot}
    c_{tot} \coloneqq \EE{C_{tot} \mid \bm{z}} = c_P + \sum_{j=1}^{2^M} p^j \cdot ( \mathrm{c}_F^j + \mathrm{c}_I^j) = c_P + \bm{p} \cdot (\bm{c}_F + \bm{c}_I).
\end{equation}
Finally, the $doa$ heuristic policy is formulated as a selection of the action vector with minimum expected total cost:
\begin{equation}
    \label{eq:doa_heur_general}
    \bm{a}_k = \Pi^{doa}_{\bm{c_{\infty}}} = \underset{\bm{z}}{\arg \min} ~\EE{C_{tot} \mid \bm{z}}.
\end{equation}
%
%
%
%
%
\subsubsection{Replacement at next decision opportunity}
\label{subsubsec:doa1}
How to quantify the costs of components for which DN was chosen at $t_k$, and then survive until the next decision time $t_{k+1}$? A first option that we investigate is to assume that components surviving until $t_{k+1}$ are then preventively replaced at $t_{k+1}$. Thus, ultimately, we consider the decision problem of replacing immediately at $t_k$ versus at $t_{k+1}$. The validity of this assumption is discussed in more detail in \Cref{sec:dicussion}. Nevertheless, this simple cost assumption proved to be quite effective in various single-component decision settings \cite{koutas2025leaf}. For this reason, it is also employed in this work.

If the policy yields that it is better to preventively replace component $i$ at $t_k$, then component $i$ is indeed replaced at $t_k$. On the other hand, if the policy favours PR of component $i$ at $t_{k+1}$, then component $i$ is not necessarily replaced at $t_{k+1}$ upon its survival. Instead, one awaits the survival of the component, collects new prognostic information, and reconsiders the decision between the two options of replacement at $t_{k+1}$ versus at $t_{k+2}$.

To obtain the cost expression for postponing the preventive replacement by one time step, we revisit the arguments from \Cref{fig:rep_graphs}. If one gets to preventively replace component $i$ at $t_{k+1}$ compared to $t_k$, then the switch to the average maintenance process (and its associated cost quantified with rate $c_{\infty,i}$) is postponed by $\Delta t$. Thus, for a single component, one saves costs $c_{\infty,i} \cdot \Delta t$, which leads to the single-component cost expression for a surviving component:
\begin{equation}
    \label{eq:cIi_doa1}
    c_{I,i} = c_f + c_{v,i} - c_{\infty,i} \cdot \Delta t.
\end{equation}
Again, the subtrahend in \Cref{eq:cIi_doa2} represents the value of extending component $i$'s life until the next decision opportunity. To consider the economic dependency of replacing multiple components at $t_{k+1}$, we include a correction term for the fixed costs, similar to \Cref{eq:cP}. To further consider the constrained number possible system states resulting from a specific action vector, we introduce the dependence on $\bm{z}$ and $\bm{x}^{(j)}$. These two amendments then lead to the multi-component cost expression for a set of surviving components:
\begin{equation}
    \label{eq:cIj_doa1}
    \mathrm{c}_I^j = \left(\sum_{i=1}^M (1-z_i) \cdot (1 - x_i^j) \cdot c_{I,i}\right) - c_f\cdot\max\left(0,~-1 +\sum_{i\in\mathcal{M}} (1 - z_i)\cdot (1 - x_i^j) \right).
\end{equation}
We call the heuristic that assumes replacement of survived components at the next decision opportunity as $doa1$.

%
%
%
%
\subsubsection{Individual perfect replacement}
\label{subsubsec:doa2}
The $doa1$ policy introduced in \Cref{subsubsec:doa1} uses preventive replacement at the next decision opportunity as a proxy for ``replace later'' of a surviving component. This renders $doa1$ a conservative policy that tends to perform a replacement too early \cite{koutas2025leaf}.

As a less conservative proxy for ``replace later'', and thus also as a less conservative policy, we propose and investigate a second option. Here, the assumption is that a surviving component can be preventively replaced right before its expected failure time conditional on $RUL_i>\Delta t$. The life extension of the component, and hence also the delay of the switch to the average maintenance process, is equal to the conditional expected failure time:
\begin{equation}
    \label{eq:cIi_doa2}
    c_{I,i} = c_f + c_{v,i} - c_{\infty,i}\cdot \EE{RUL_i \mid RUL_i > \Delta t}.
\end{equation}
The assigned cost for a surviving component in \Cref{eq:cIi_doa2} is strictly lower than the assigned cost of $doa1$ in \Cref{eq:cIi_doa1}, leading to a less conservative policy. Taking into account the dependence on $\bm{z}$ and $\bm{x}^j$, the multi-component cost expression for a set of surviving components takes the form of
\begin{equation}
    \label{eq:cIj_doa2}
    \mathrm{c}_{I}^j = \sum_{i=1}^M (1-z_i) \cdot (1 - x_i^j) \cdot c_{I,i}.
\end{equation}
Here, the economic dependence of the individual components is \emph{not} considered, as the assumed replacement time of component $i$ is continuous, and hence is unlikely to coincide with the replacement times of other components.

We call the heuristic that assumes preventive replacement of survived components at their conditional expected failure time $doa2$.

%
%
%
%
\subsection{Benchmark policies}
\label{subsec:benchmark_policies}
For thorough evaluation of the proposed $doa$ heuristics, we compare their performance against a set of existing replacement heuristics. Here, plenty of options are available, as outlined in more detail in the introduction. The benchmark heuristic replacement policies are termed $rhi$ and are introduced in the following \Cref{subsubsec:probability_thresholds,subsubsec:reliability_threshold}.

%
%
%
%
\subsubsection{Component-wise probability thresholds}
\label{subsubsec:probability_thresholds}
The first benchmark heuristic replacement policy, termed $rh1$, is a standard strategy used frequently in the literature for risk-based inspection planning \cite[e.g.,][]{luque2019risk,nielsen2018computational}) and PHM \cite[e.g.,][]{kamariotis2024metric,nguyen2019new,he2023condition}). $rh1$ is a component-wise PdM policy (i.e., no consideration of system effects) that preventively replaces component $i$ when its predicted probability of failure until the next decision opportunity exceeds a threshold value $p_{\mathrm{thres},i}$:
\begin{equation}
    \label{eq:rep_threshold}
    a_{k,i} = \Pi^{rh1}_{p_{\mathrm{thres},i}} = 
    \begin{cases}
        \mathrm{PR}, & \text { if } 
        p_{\mathrm{thres},i} < \Pr\left(RUL_i \leq \Delta t\right) \\ 
        \mathrm{DN} & \text { else}.
    \end{cases}
\end{equation}
Since the individual actions are chosen independently of each other, the final joint action vector is obtained as
\begin{equation}
    \label{eq:rh1}
    \bm{a}_k = [a_{k,1},~\ldots,~a_{k,M}].
\end{equation}
The parameters of $rh1$ are the component-specific probability thresholds, i.e., $\textbf{w}=\bm{p}_\mathrm{thres}=[p_{\mathrm{thres},1},~\ldots,~p_{\mathrm{thres},M}]$. A default choice for the probability thresholds is $p_{\mathrm{thres},i}=(c_f+c_{v,i})/c_c$\footnote{This threshold can be derived by comparing the preventive maintenance costs $c_f+c_{v,i}$ with the direct expected failure costs until the next decision time $\Pr\left(RUL_i \leq \Delta t\right)\cdot c_c$.}. However, this choice is suboptimal and can perform poorly \cite{kamariotis2024metric,koutas2025leaf}. Thus, if sufficient training data is available, it is preferable to find $\bm{p}_{\mathrm{thres}}^*$ by solving the optimization problem in \Cref{eq:optimal_heuristic_strategy}. In this work, we solely use the policy with the optimized thresholds, denoted as $rh1^*$.
%
%
%
%
\subsubsection{System reliability threshold}
\label{subsubsec:reliability_threshold}
The second investigated heuristic benchmark policy, denoted $rh2$, was introduced by \citet{nguyen2015multi}\footnote{\citet{nguyen2015multi} consider a much more general replacement setting, i.e., a general system consisting of series and parallel subsystems, the inclusion of inspections, as well as the distinction of critical and non-critical components. Here, we introduce their policy adapted to our specific decision setting.}. It consists of a two-stage process, where first the series system's reliability $R=\prod_{i=1}^M R_i = \prod_{i=1}^M \Pr(RUL_i>\Delta t)$ is compared against a threshold $r_\mathrm{thres}$.

If the system reliability falls below this threshold, then in the second stage, all possible action vectors $\bm{z}$ are checked. Out of all actions that improve the system's reliability above the reliability threshold, the action with the highest reliability improvement to action cost ratio is chosen. Mathematically, this policy is expressed as: 
\begin{equation}
    \label{eq:rh2}
    \bm{a}_{k} = \Pi^{rh2}_{r_{\mathrm{thres}}} = 
    \begin{cases}
        \bm{a}_{\max}, & \text { if } 
        R < r_{\mathrm{thres}} \\ 
        [\mathrm{DN},~\ldots,~\mathrm{DN}] & \text { else},
    \end{cases}
\end{equation}
where $\bm{a}_{\max}$ is found via
\begin{equation}
    \label{eq:amax}
    \bm{a}_{\max} = \underset{\bm{z} \in \mathcal{Z}_r}{\arg \max} \frac{R(\bm{z}) - R}{\EE{C_P \mid \bm{z}}}, \qquad \text{where} ~ \mathcal{Z}_r=\{\bm{z} \mid R(\bm{z}) \geq r_{\mathrm{thres}}\}.
\end{equation}
$\EE{C_P \mid \bm{z}}$ and $R(\bm{z})=1-\Pr(\bm{x}^0 \mid \bm{z})$ can be computed from \Cref{eq:cP,eq:pr_xj_given_z}, respectively.

The sole free parameter of $rh2$ is the system reliability threshold, $r_{\mathrm{thres}}$, which can either be defined a-priori, or found by solving the optimization problem in \Cref{eq:optimal_heuristic_strategy}. In this work, we exclusively use the policy with the optimized system reliability threshold, denoted as $rh2^*$.

%
%
%
%
\subsection{Limited training data}
\label{subsec:low_data}
In many engineering applications, especially for safety-critical systems, run-to-failure data is usually sparse. Hence, parameter initialization/optimization of a chosen heuristic has to be performed with a limited set of training data. In this case, both the initialization function in \Cref{eq:ectr_init} and the performance function in \Cref{eq:optimal_heuristic_strategy} are not smooth over the whole parameter domain, which poses challenges. In the following, we outline the differences of $doa$ parameter initialization to the benchmarks' parameter optimization in the low data regime.

%
%
%
%
\subsubsection{Doa parameter initilization with few training data}
\label{subsubsec:doa_init_low_data}
For finding $c_{\infty,i}^0$, merely a set of time-to-failure samples is needed, based on which the optimal deterministic replacement time for the population of components is computed. Information from a prognostic model is not required, as $c_{\infty,i}^0$ is found under the assumption of simple age-based replacement. For $M$ components, $M$ separate and independent single-parameter optimizations are performed. Furthermore, this approach is agnostic to the underlying system logic and, hence, can be used for any type of system. 

\begin{figure}[H]
    \centering
    \begin{subfigure}{0.49\textwidth}
        \includegraphics[width=\textwidth]{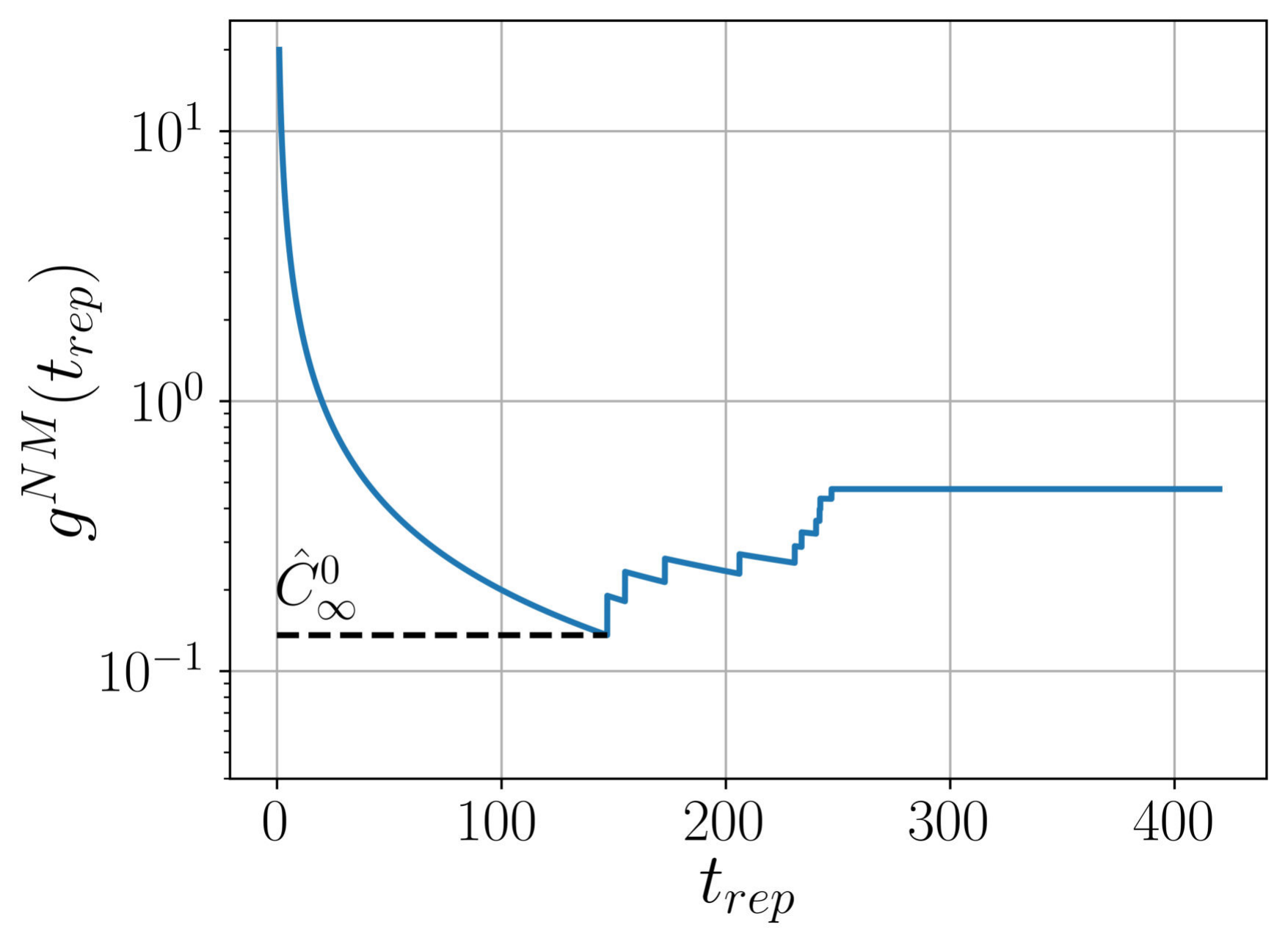}
        \label{subfig:doa_noisy_init_10}
    \end{subfigure}
    \hfill
    \begin{subfigure}{0.49\textwidth}
        \includegraphics[width=\textwidth]{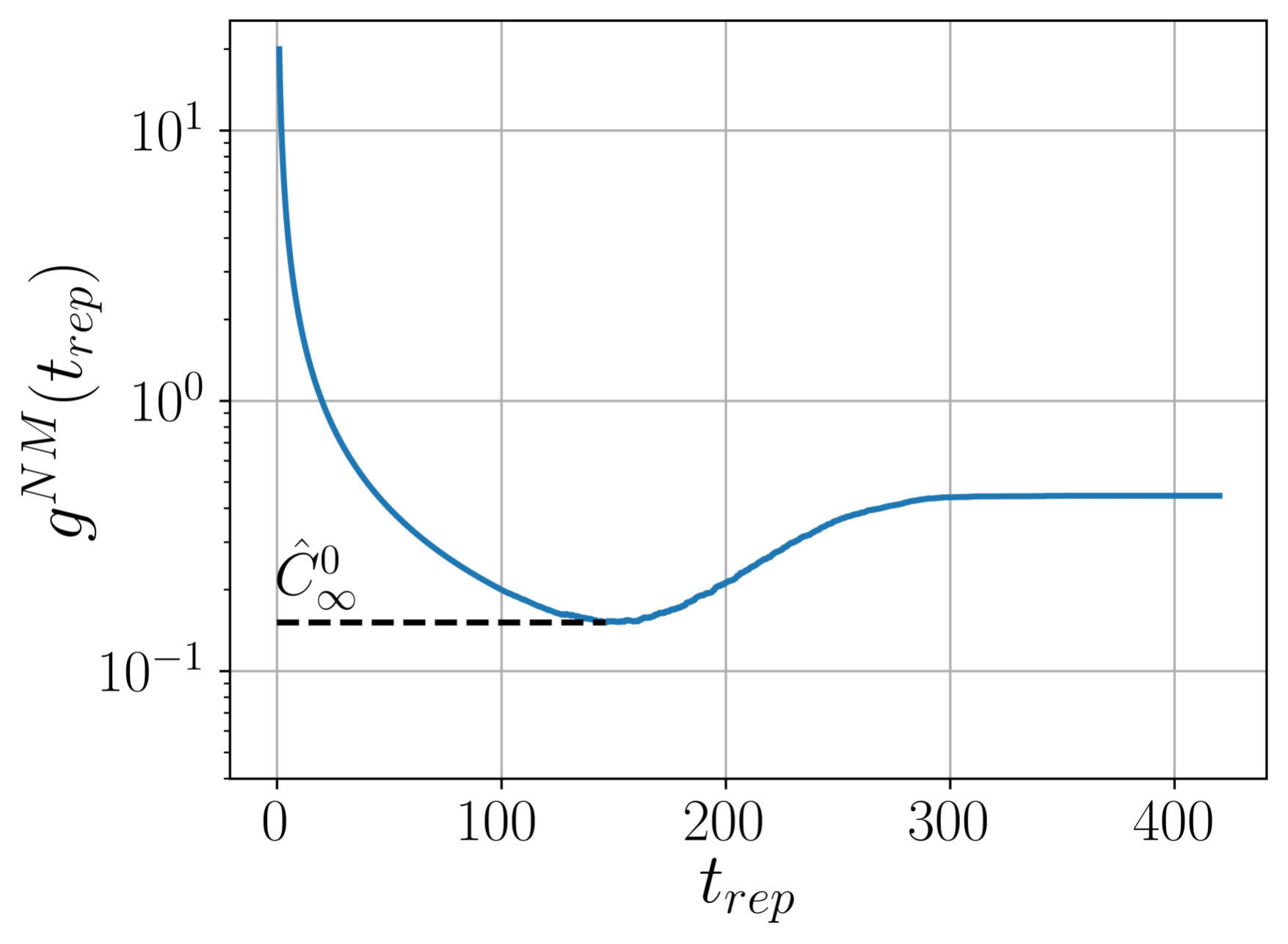}
        \label{subfig:doa_noisy_init_1000}
    \end{subfigure}
    \vspace{-1em}
    \caption{Depiction of $doa$'s initialization target function for a component with $c_f=10$, $c_v=10$, $c_c=100$ computed with 10 (left) and $1,000$ (right) failure samples.}
    \label{fig:doa_noisy_init}
\end{figure}

The shape of the objective function for different sample sizes is shown in \Cref{fig:doa_noisy_init}. For $n$ samples, the target $g^{NM}(\trep)$ has $n$ jumps, between which the function is monotonically decreasing. Finding the global minimum of this function is trivial and can be performed in negligible time via fine discretization of $t$, with corresponding entry-wise evaluation of \Cref{eq:obj_func_ectr_init_nm}. With increasing sample size, the target function gets increasingly smoother.
%

%
%
%
%
\subsubsection{Parameter optimization with few training data}
\label{subsubsec:benchmark_parameter_optimization_few_data}
Finding a heuristic's optimal parameter values with respect to the final policy performance constitutes a more difficult task than the parameter initialization outlined in \Cref{subsubsec:doa_init_low_data}. Firstly, one needs a set of time-to-failure samples \emph{as well as} corresponding prognostic information over the individual components' life cycles. In addition, the optimization is generally performed jointly over all parameters of the heuristic policy (instead of separately and independently).

Still, optimization of the benchmark policies' parameters is necessary because either the default parameter values lead to poor performance \cite[see, e.g.,][]{kamariotis2024metric,koutas2025leaf}, or default values for the parameters are not available. By contrast, for the $doa$ policies, reasonable default values for $c_{\infty,i}$ can be derived, which is why optimization with respect to the policy's performance is not necessary. If one nevertheless wants to perform optimization to find $c_{\infty,i}^*$, the procedure outlined in the following can be applied as well. 
\begin{figure}[H]
    \centering
    \begin{subfigure}{0.49\textwidth}
        \includegraphics[width=\textwidth]{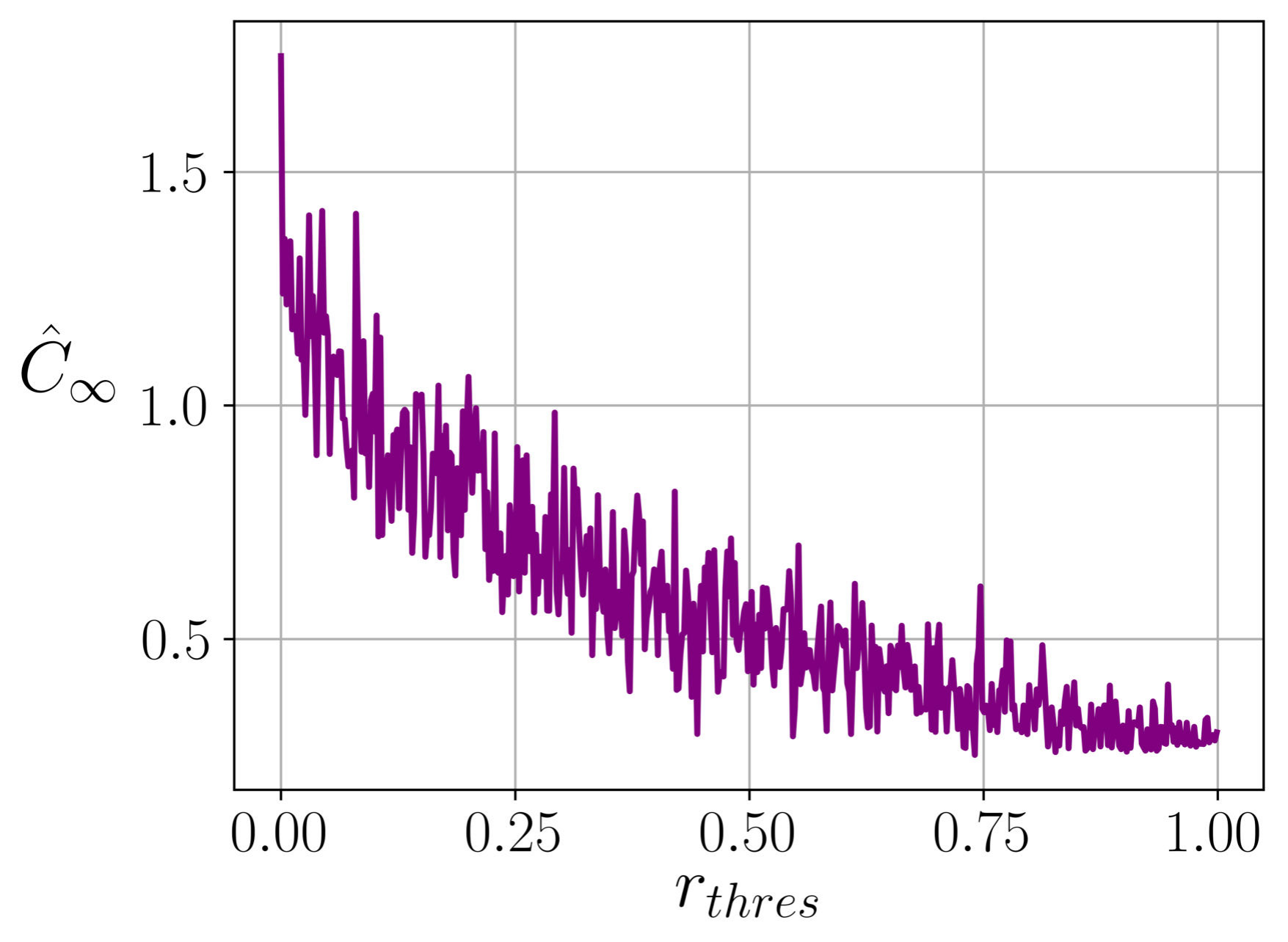}
        \label{subfig:rh2_noisy_target_10}
    \end{subfigure}
    \hfill
    \begin{subfigure}{0.49\textwidth}
        \includegraphics[width=\textwidth]{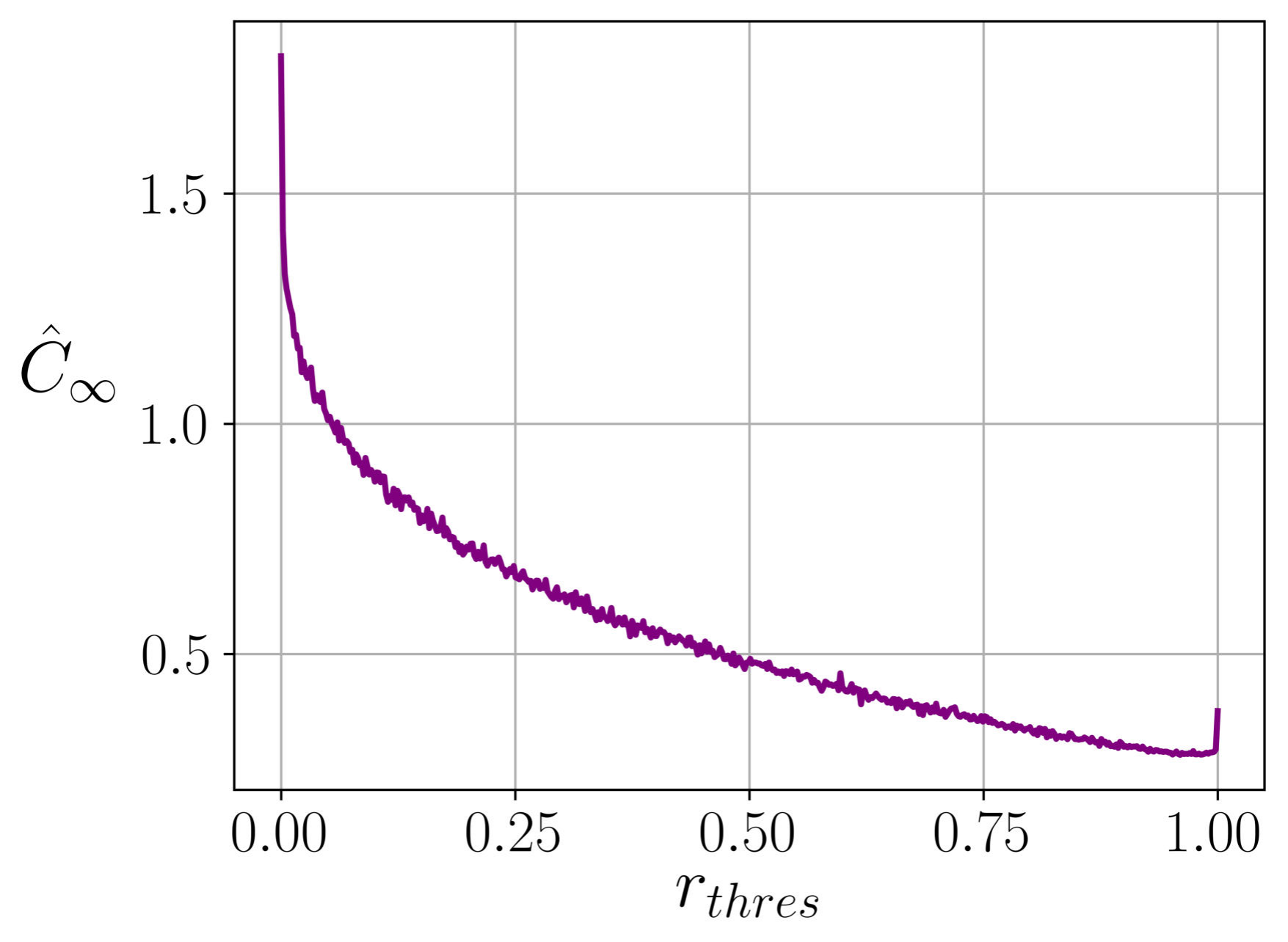}
        \label{subfig:rh2_noisy_target_1000}
    \end{subfigure}
    \vspace{-1em}
    \caption{$rh2$ performance ($\hat{C}_{\infty}$) over the domain $r_{\mathrm{thres}}\in[0,1]$ shown for a set of 10 (left) and $1,000$ (right) training data samples per component. The specific setup is the same as in \Cref{subsec:limited_training_data}.}
    \label{fig:noisy_ectr}
\end{figure}
The shape of the performance function for different sample sizes is shown in \Cref{fig:noisy_ectr}. Again, with increasing amounts of training data, the underlying target function becomes increasingly smooth; hence, solving the corresponding optimization problem becomes easier. However, in case of low training data availability, optimization is not straightforward since the underlying $\hat{C}_{\infty}$ is noisy. This is because even small parameter perturbations can have large effects on the maintenance performance due to the components' economic dependence.\footnote{For a single-component system, $\hat{C}_{\infty}$ can be piecewise constant over the parameter space \cite{koutas2025leaf}.} Since the failure costs are high compared to the preventive replacement costs, the optimal system reliability threshold is also high with $r_{thresh}\approx0.98$ for $1,000$ samples in \Cref{fig:noisy_ectr}.

Optimization of this noisy target function (especially in higher dimensions) is challenging. Due to the ubiquity of these problems in many scientific areas, various approaches for handling this noise have been developed over the past decades, such as Bayesian optimization or evolutionary algorithms. In this work, we opt for a stochastic optimizer that is inspired by natural selection, namely a genetic algorithm (GA). 

\Cref{sec:app_optimization_via_ga} provides an overview of the workings of GAs along with the specific implementation adopted in this study.

%
%
%
%
\section{Numerical investigations on a virtual RUL simulator}
\label{sec:results}
A good policy should perform well over a broad range of application settings. For this reason, we test the performance of the proposed $doa$ heuristics with parameter initialization and compare it to the optimized benchmark policies $rhi^*$ over a wide range of settings. These include varying amounts of components to show scaling properties, numerous cost settings by varying fixed to variable costs $c_f/c_v$ as well as varying component importances $c_{v,1}/c_{v,2}$. We also assess the robustness of the heuristics in applications with few training data.

For all investigations, $\Delta t=10$ is fixed as the decision interval, and $c_c=100$ is fixed for the corrective costs. Furthermore, we choose the investigation settings so that $c_f + \sum_{i=1}^M c_{v,i} \leq c_c$ to reflect higher consequences of corrective replacement.

%
%
%
%
\subsection{Dataset}
\label{subsec:data_set}
We investigate the performance of the proposed heuristic policies in a controlled setting that enables accurate testing and verification. To this end, we make use of the RUL simulator introduced in \cite{kamariotis2024metric}. The simulator provides \emph{calibrated}\footnote{A well-calibrated probabilistic model outputs predictions whose credible intervals correspond to the observed frequencies in the long run \cite{gneiting2007probabilistic}, i.e., the true RUL lies in the $(1-\alpha)\%$ credible interval in $(1-\alpha)\%$ of cases, with sufficient test data.} RUL predictions. Moreover, the use of the simulator enables the generation of large amounts of test data, which facilitates the exact evaluation of the different policies' performances.

Without loss of generality, we assume the same underlying TTF distribution for all components, which is why the component index $i$ is neglected in the following. The TTF distribution of the components takes the form of the Gaussian distribution: $T_{F} \sim \mathcal{N}(\mu=225, \sigma=40)$.\footnote{For the $doa$ policies, we assume that no component fails until the next decision time. In the implementation, this is enforced by rejecting a new component if its failure time is smaller than $\Delta t$; hence slightly altering the assumed normal TTF distribution.} The RUL simulator models prediction errors $\epsilon_k$ at time step $t_k$ by a lognormal process $RUL_{k}=\epsilon_k\cdot RUL_{\mathrm{true},k}$, where the logarithms of the prediction errors $\ln (\epsilon_k)$ follow a zero-mean multivariate normal distribution. The resulting marginal predicted RUL distributions, $f_{RUL_k}$, are lognormal with mean $\mu_{RUL_k}\approx1.083\epsilon_k\cdot RUL_{\mathrm{true},k}$ and standard deviation $\sigma_{RUL_k}\approx0.451\epsilon_k\cdot RUL_{\mathrm{true},k}$. Thus, the simulator yields slightly biased RUL predictions, which is, however, not critical for performance. The RUL predictions at different $t_k$ are correlated, with correlation length equal to $l_{\mathrm{corr}}=50$. 
A detailed description of the RUL simulator can be found in \Cref{sec:virtual_RUL_simulator}. An example realization of RUL predictions obtained with this simulator is shown in \Cref{fig:decision_setting}.
%
%
%
%
\subsection{Performance with abundant training data}
\label{subsec:abundant_training_data}
We start the investigation by assuming abundant training data is available, i.e., by setting the batch size $B$ and simulation horizon $H$ high enough such that the system's computed $\hat{C}_{\infty}$ has converged (see \Cref{fig:convergence_plot}). Here, parameter initialization/optimization for each policy is not an issue, as the respective target functions are smooth and a variety of (gradient-based) optimizers can be used to find the initialized/optimized parameters (see, e.g., \Cref{subfig:doa_noisy_init_1000,subfig:rh2_noisy_target_1000}).

\begin{figure}[H]
    \centering
    \includegraphics[width=0.9\linewidth]{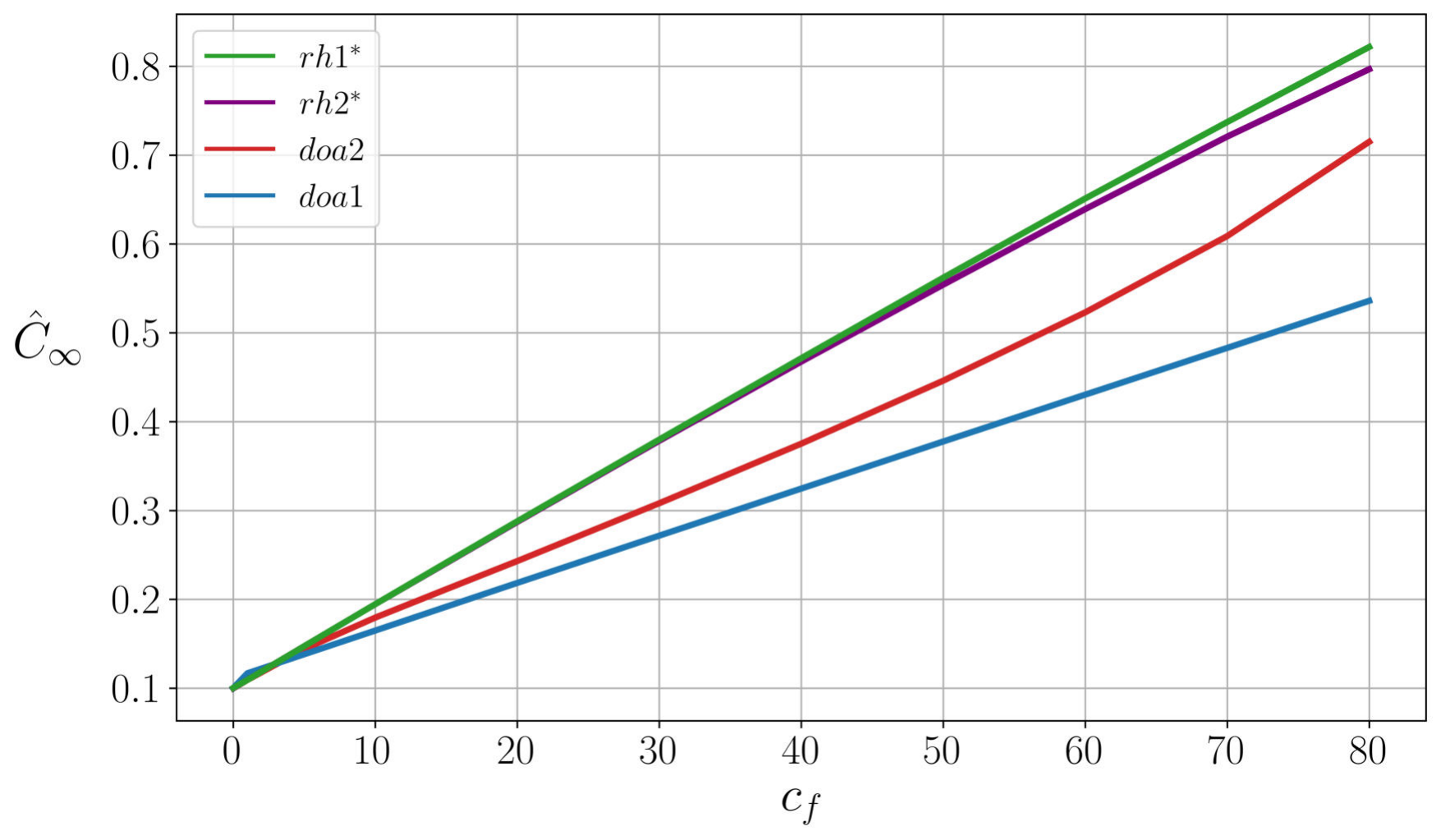}
    \caption{Comparison of the proposed $doa$ policies with the benchmarks $rh1^*$ and $rh2^*$ for a 2-component system over a range of fixed costs $c_f$. The following settings were used for training the optimized policies as well as for computing the system's $c_{\infty}$ via \Cref{eq:ectr_sys} with Monte Carlo simulation: $B=2,000$ \& $H=500$ for $doa2$, and $B=2,000$ \& $H=1,000$ for $rh1^*$, $rh2^*$, and $doa1$. The corrective and all variable costs were fixed to $c_c=100$ \& $c_{v,i}=10$. The figure shows the mean of each policy as a solid line, as well as the corresponding 95\% credible intervals.}
    \label{fig:2comps_cf}
\end{figure}

We first investigate the performance over a wide range of $c_f$ values, while keeping the variable costs, $c_{v,i}$, fixed. This covers industries with minimal set-up cost, e.g., due to low technician travel times or automated replacements, up to industries with high set-up costs, e.g., due to long technician travel times or the use of specialized equipment. The corresponding results for a 2-component system are shown in \Cref{fig:2comps_cf}. 

\begin{figure}[H]
    \centering
    \includegraphics[width=0.9\linewidth]{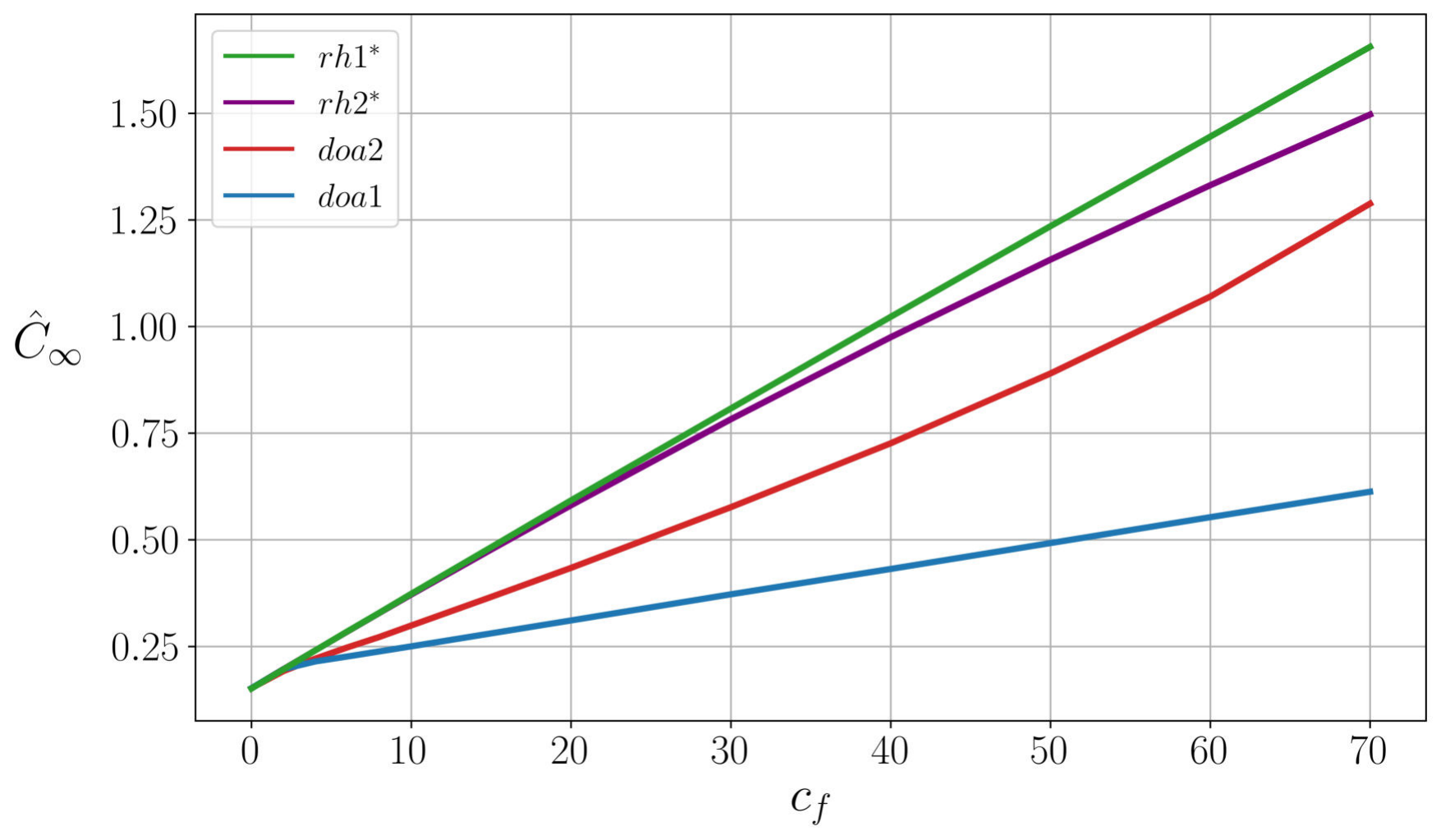}
    \caption{Comparison of the proposed $doa$ policies with the benchmarks $rh1^*$ and $rh2^*$ for a 5-component system over a range of fixed costs $c_f$. The following settings were used for training the optimized policies as well as for computing the system's $c_{\infty}$ via \Cref{eq:ectr_sys} with Monte Carlo simulation: $B=500$ \& $H=400$ for $doa2$, and $B=2,000$ \& $H=1,000$ for $rh1^*$, $rh2^*$, and $doa1$. The corrective and all variable costs were fixed to $c_c=100$ \& $c_{v,i}=6$. The figure shows the mean of each policy as a solid line, as well as the corresponding 95\% credible intervals.}
    \label{fig:5comps_cf}
\end{figure}
One can see that for low set-up costs $c_f < 4$, where coordination and planning of multiple component replacements is less crucial, the performances of all policies are quite similar. With increasing $c_f$, where coordination is increasingly more important, the proposed $doa$ policies significantly outperform the benchmarks. 
$doa1$ emerges as the clear winner, with cost reductions of up to $35\%$ at $c_f/c_{v,i}=8$. Regarding the benchmarks, $rh2^*$ slightly outperforms $rh1^*$ over the whole investigated parameter range, with increasing difference for increasing $c_f$. This is attributable to $rh2$ explicitly considering all action combinations and the economic dependence of the components, whereas $rh1$ independently chooses replacement actions for each component.

In summary, the long-running maintenance costs per unit time of all policies increases linearly with increasing $c_f$, but the slope of the $doa$ heuristics is significantly smaller compared to the benchmark policies. To confirm the scalability of the $doa$ heuristics, this investigation is performed also for 5 components; the results of which are shown in \Cref{fig:5comps_cf}. The pattern of increasing advantage of the $doa$ policies with increasing $c_f$ can again be observed, with achieved cost reductions of $doa1$ of up to $65\%$ for $c_f/c_{v,i}=11.7$.

\begin{figure}[H]
    \centering
    \includegraphics[width=0.88\linewidth]{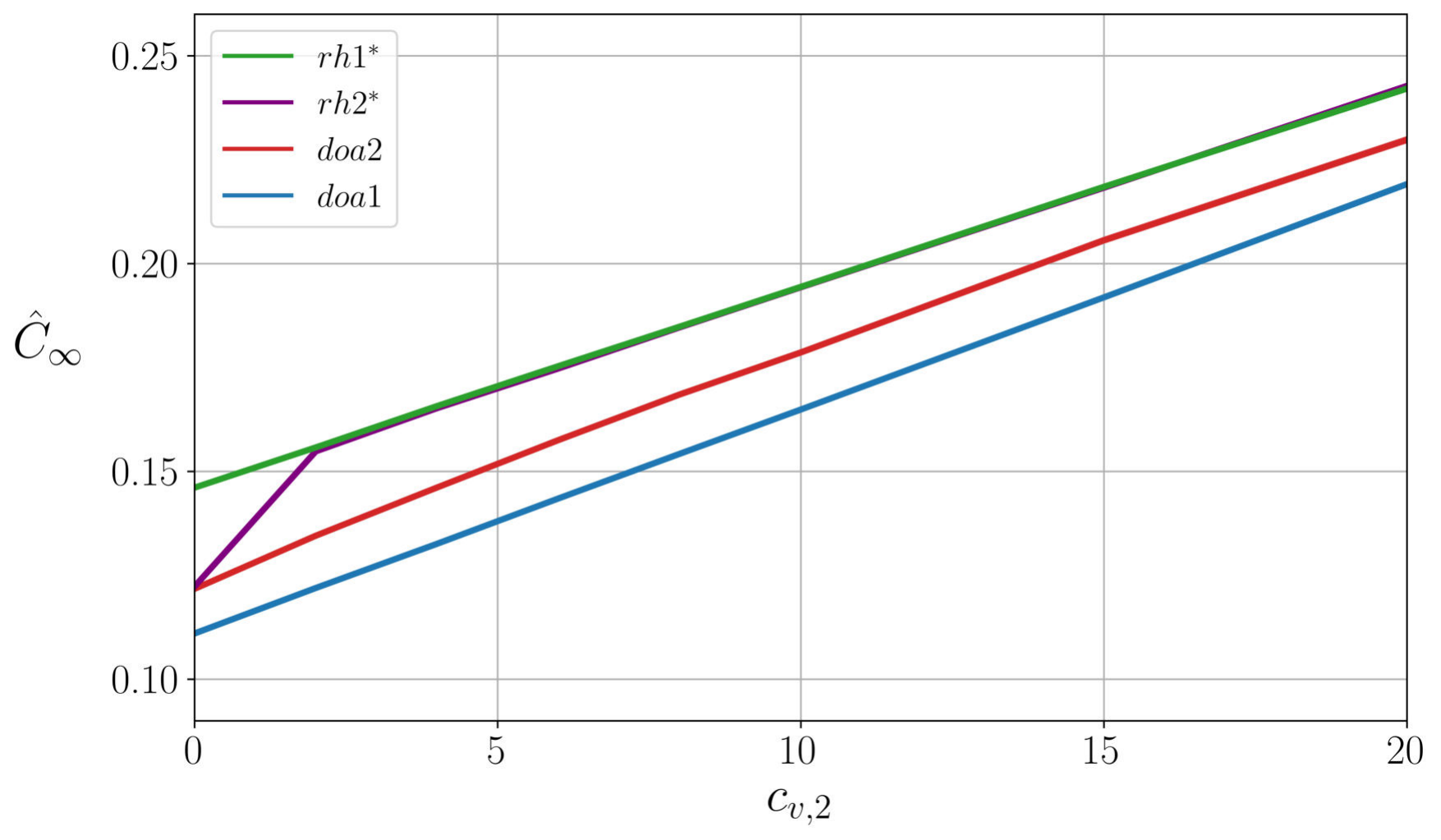}
     \caption{Comparison of the proposed $doa$ policies with the benchmarks $rh1^*$ and $rh2^*$ for a 2-component system over a range of variable costs $c_{v,2}$. The following settings were used for training the optimized policies as well as for computing the system's $c_{\infty}$ via \Cref{eq:ectr_sys} with Monte Carlo simulation: $B=1,000$ \& $H=500$ for $doa2$, and $B=2,000$ \& $H=1,000$ for $rh1^*$, $rh2^*$, and $doa1$. The remaining costs are fixed to $c_c=100$, $c_{f}=10$ \& $c_{v,1}=10$. The figure shows the mean of each policy as a solid line, as well as the corresponding 95\% credible intervals.}
    \label{fig:2comps_cv}
\end{figure}

For the next investigation, we plot the performance of all policies over a wide range of relative component importances, $c_{v,2}/c_{v,1}$, while keeping $c_f$ fixed. This investigation covers series systems of components with equal cost, as well as combinations of cheap and expensive components. The corresponding results, again for a 2-component system, are shown in \Cref{fig:2comps_cv}.
The benchmarks again show similar performance over the majority of the investigated range. The noticeable difference occurs for $c_{v,2}=0$. Here, component $2$ should always be preventively replaced whenever PR is chosen for component 1, as no additional costs are incurred. The policies that explicitly consider the individual action costs in the selection process ($rh2$, $doa1$, $doa2$) can capture this behaviour and hence significantly outperform $rh1$.

Regarding the proposed policies, both $doa1$ and $doa2$ consistently outperform the benchmarks, where $doa1$ yields the best performance over the whole investigated range. If one would extend the performance curves until $c_{v,2}/c_{v,1}=80$, the performance of $doa1$ gets gradually worse until it becomes the worst-performing policy. The reason for this result lies in the conservatism of $doa1$. For low variable costs, being conservative leads to earlier replacements; hence, more preventive replacements and fewer failures. Due to the low PR cost, it is preferable to choose more frequent PR in exchange for fewer costly failures. This is why the conservative $doa1$ outperforms the less conservative $doa2$ (and $rhi^*$). By contrast, for high $c_{v,2}$, the opposite is the case, i.e., it is preferable to accept more failures. Here, $doa2$ achieves the best performance. Nevertheless, the case of $c_{v,2}/c_{v,1}=20$ already represents an extreme cost imbalance of components; even higher cost ratios are unlikely to occur in practice. Thus, these cases are of minor importance compared to the shown range, where $doa1$ achieves the best performance. 

In summary, the long-running maintenance costs of all policies increase linearly with increasing $c_{v,2}$. The slopes are similar for all policies, but $doa1$ stands out as the best-performing policy and with the steepest performance slope. Based on the results of \Cref{fig:2comps_cf,fig:5comps_cf}, smaller/larger $c_f$ would lead to smaller/larger performance gaps, but with sustained slope differences.

%
%
%
%
\subsection{Performance with limited training data}
\label{subsec:limited_training_data}
We investigate how the policies perform in low data regimes. Here, the respective target functions for parameter initialization/optimization are noisy (see, \Cref{subfig:doa_noisy_init_10,subfig:rh2_noisy_target_10}). Parameter initialization for the $doa$ policies is not an issue, but parameter optimization for the benchmark heuristics becomes more challenging. We use the genetic algorithm described in \Cref{sec:app_optimization_via_ga} to find the optimal set of parameters. 

\begin{figure}[H]
    \centering
    \includegraphics[width=0.9\linewidth]{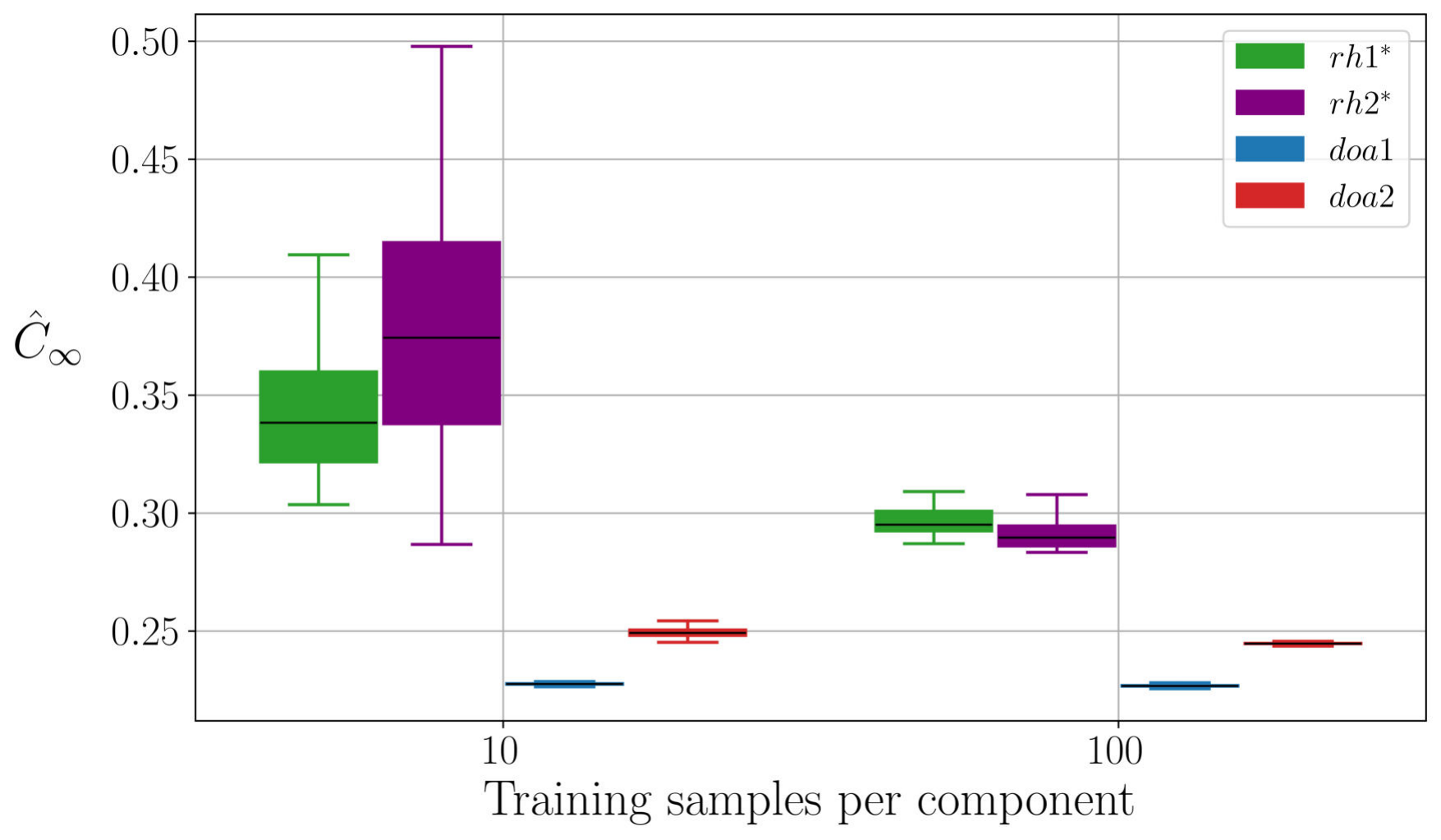}
    \caption{Boxplot comparison of the proposed $doa$ policies with initialized parameters against the benchmarks $rhi^*$ with optimized parameters. A 5-component system is used with $c_c=100$, $c_f=6$, $\bm{c}_{v}=[2,4,6,8,10]$. The initialization/optimization was carried out for 10 and 100 TTF training samples per component; the subsequent evaluation of the policies with the initialized/optimized parameters was conducted with $B=500$, $H=400$. The investigation was repeated 100 times for all policies but $doa2$, for which it was repeated 50 times; the respective medians of the runs are highlighted in black.}
    \label{fig:low_data}
\end{figure}

In these low data settings, the final maintenance performance can highly depend on the specific set of TTF samples. This heuristic-specific sample dependence is shown in \Cref{fig:low_data} for the repeated evaluation of a 5-component system with only 10 and 100 TTF samples per component, respectively. One can clearly see that the variance in final performance is much higher for the benchmarks compared to the proposed $doa$ policies.

Evidently, parameter optimization on low amounts of data leads to heavy overfitting. By contrast, the proposed doa policies show excellent robustness against overfitting, even in very low data domains. $Doa1$ again emerges as the best policy, showing both the best median performance as well as the lowest performance variance. 
%

%
%
%
%
\section{Discussion}
\label{sec:dicussion}

%
%
%
%
\subsection{Strengths and weaknesses of the proposed policies}
\label{subsec:strengths_and_weaknesses}
The proposed $doa$ policies were shown to outperform the benchmarks in the numerical investigations. Even more importantly, they possess essential practical advantages over existing data-driven policies for replacement decisions. 

A key advantage is that the policies can be deployed in low-data regimes. 
No (joint) parameter optimization over the free parameters, $c_{\infty,i}^0$, is necessary, as initialization options are available that already lead to superior performance. The initialization of the component-specific cost rates can be performed sequentially and independently, which constitutes a simple task. Furthermore, finding the optimum of the underlying target function is trivial via fine discretization, even when only few training samples are available. For the initialization option used in this work, merely a set of TTF samples is needed; information from a prognostic model is not even required. As evidenced by the results in \Cref{subsec:low_data}, this initialization is highly robust against overfitting. This effect has already been observed in \cite{koutas2025leaf,baker1995can} in the case of maintenance of single-component systems. Hence, the $doa$ policies can be reliably used also in low-data domains and when no past prognostic data is available. This is contrasted by the benchmark policies, which, in case of low data availability, require an additional step of noisy optimization. This constitutes a hard task, leads to heavy overfitting\footnote{The genetic algorithm used in this work resulted in heavy overfitting to the respective TTF samples and, thus, in large variance of the final policy performance. It is possible that better results (i.e., lower variance and better mean) of the benchmark policies can be achieved for different optimization procedures, e.g., via the usage of Bayesian optimization \cite[see, e.g., ][]{frazier2018tutorial}. However, we suspect that this would merely lead to a slight mitigation of the overfitting issue, as this effect was already observed for single component systems \cite{koutas2025leaf}. In addition, a more detailed investigation of different optimization procedures is not the scope of this work.}, and still results in worse performance.

Another advantage of the $doa$ policies is that they are based on an online evaluation of the available decision options at each decision time step $t_k$. Thus, one can perform various changes to the underlying system during production without an effect on the performance of the heuristics. One such possible amendment is the change of the cost structure, e.g., when a new replacement method is developed that reduces the current variable cost of a component, $c_{v,i}$. Here, one simply needs to recompute the $c_{\infty,i}^0$ (with the already available TTF samples) and adjust the cost formulations for $c_P$ in \Cref{eq:cP} and $c_I^j$ in, e.g., \Cref{eq:cIj_doa1} with this newly computed cost ratio. Another possible amendment is the addition/removal of components, where the corresponding decision tree can be simply expanded/reduced. By contrast, for the benchmark policies, a change of the cost structure necessitates a recomputation of \emph{all} policy parameters. If the system is expanded by new components, one would first need to collect run-to-failure prognostics data for these new components, before parameter optimization can even be performed. 

The main drawback of the $doa$ policies is that their computational cost scales badly with the number of components. With two available actions, the one-step decision tree has $2^{2M}$ branches, which even for moderate number of components $M$ renders evaluation of longer horizons $H$ difficult at best and unfeasible at worst. 
In \Cref{subsubsec:reduction_of_computational_effort}, we propose a few approaches for reducing this computational effort in future work. 

Moreover, the decision tree for the $doa$ policies is constructed for a specific type of system. While changes of the cost structure or number of components can be trivially implemented, a change to the system structure, e.g., via the introduction of a redundant component, necessitates a larger reworking of the decision tree. Thus, the analysis of the individual branches would need to be redone for a custom system. In \Cref{subsubsec:amendment_to_other_systems}, we outline a few steps on how this analysis can be performed also for other systems. 

\subsection{Suggestions for future work}
\label{subsec:future_work}
%

%
%
%
%
\subsubsection{Reduction of computational effort}
\label{subsubsec:reduction_of_computational_effort}
The main limitation of the proposed $doa$ policies is that its computational cost scales poorly with the number of components. In general, the most expensive operations are the conditional expectations, $\EE{RUL_i \mid \cdot}$ in \Cref{eq:expected_failure_costs_comp_i,eq:cIi_doa2}. Computing these conditional RUL expectations is unnecessary for components that are far from failure. Therefore, a potential way to reduce the computational effort by at least an order of magnitude is to not evaluate selected action branches if they fulfill some criteria. One option is to formulate appropriate upper/lower bounds based on the cost structure ($c_f,~c_{v,i},~c_c$) and failure probabilities $p_{F,i}$ (which are comparably fast to evaluate) that serve to exclude branches. Such criteria should be derived separately for each of the presented policies.
 
A more advanced approach would be to use dynamic upper/lower bounds and to prune suboptimal branches based on these computed bounds. Numerous search strategies (i.e., the order in which branches are evaluated) and pruning rules (i.e., rules that prevent the evaluation of suboptimal branches) have been proposed in the discrete optimization community. For an overview, the reader is referred to \cite{morrison2016branch}.

Moreover, in the current implementation, the simulation horizon $H$ is chosen generously to ensure $\hat{C}_{\infty}$ convergence with a corresponding visual check. This can be improved by the inclusion of automatic convergence checks in the simulation process, leading to a dynamic adjustment of $H$ and thus to a reduction of computational overhead.

%
%
%
%
\subsubsection{Further policies}
\label{subsubsec:further_policies_and_benchmarks}
%
Future work could explore other $doa$ policy variants than the two investigated ones, such as $doa2$, but with the inclusion of economic dependence in case multiple components are scheduled to be replaced at similar times. 
Another option is to use the proposed versions, but to perform optimization to find $c_{\infty,i}^*$. This would likely lead to further performance gains, but with the drawbacks of higher computational cost and overfitting in low data domains \cite[see also][]{koutas2025leaf}.

%
%
%
%
%
\subsubsection{Adaptation of the doa policies to other systems}
\label{subsubsec:amendment_to_other_systems}

The series system is frequently encountered in engineering applications and therefore has substantial practical relevance. Nevertheless, future research should extend the formulation of the proposed $doa$ policies to $k$-out-of-$n$ systems and general systems.

Adapting the $doa$ policies to other systems with economic dependence is not trivial, but we provide a few suggestions. Firstly, one can use the general decision tree in \Cref{fig:dt} as a starting point to outline all the possible action and component vector states, with the corresponding branch probabilities derived from the individual prognostic models. Only the cost estimations of the individual branches have to be adjusted to the specific system. For general serial-parallel systems, not all component failures lead to system failure. Hence, one needs to specify further characteristics, such as self-announcing vs. hidden component failures. Such a general system has a different structure function, i.e., a different $s^j$ column in the truth \Cref{tab:system_truth_table}, which is trivial to determine. The updated truth table can then be used with the amended cost functions to find the expected costs for each branch. The main challenge lies in the formulation of suitable renewal theory-based cost estimations for failed and surviving components, i.e., $c_{F,i},~c_I^j$ in \Cref{eq:cF,eq:cIj_general}. In the current implementation, a component is immediately replaced upon PR/failure.\footnote{A convenient consequence of this assumption is that all components are working at all decision times, enabling a stateless implementation.} For a general system with parallel parts, this is not necessarily the case, and one can have failed components that are not immediately replaced. Therefore, the cost expressions for individual components depend on the states of other components, which introduces significant complications.

We have considered systems with economic dependence, but without structural or stochastic dependence, to limit the complexity of the numerical investigations. Treating systems with structurally dependent components within the framework is straightforward, since only the cost estimations of the individual branches have to be adjusted to incorporate the additional costs of forced joint replacements. Furthermore, the consideration of systems consisting of dependent components constitutes a prognostic challenge rather than a maintenance planning challenge. With stochastic dependence, the state of one component can influence the lifetime distribution of other components. This necessitates the switch from $M$ marginal RUL distributions to an $M-$dimensional joint RUL distribution. Training such a joint prognostic model is significantly more challenging than training a set of independent prognostic models, but this is not the focus of this paper. If a joint prognostic model is available, then we hypothesize that our proposed policies can be applied to systems with stochastically dependent components as well.

%
%
%
%
\section{Concluding remarks}
\label{sec:concluding_remarks}
We present novel policies for preventive replacement of an $M$-component series system with economic dependence. The policies are derived by constructing a reduced decision tree to model the available decision options and their consequences. The probability of reaching a certain branch is computed with the set of available RUL-PDFs obtained from the component-specific prognostic models. The expected costs of each action are
then computed as the sum of the direct action and failure costs, as well as the costs of the assumed underlying renewal-reward process. The latter accounts for the costs associated with the components beyond the limited time horizon of the decision tree. In this way, the derived $doa$ policies account for the benefits of potential life extensions of individual components by replacing them later rather than early.

We compared our proposed $doa$ policies against two benchmark policies on a simulated dataset for various cost structures and different data regimes. Our results show that the $doa$ policies can significantly outperform the presented benchmark policies and achieve a reduction of the long-running cost rate of up to 35\%  reduction for a 2-component system. This performance improvement was achieved using only the initialized $doa$ parameters against the optimized benchmark policies. If one optimizes $doa$ parameters, $c_{\infty,i}^*$, the performance gap to the benchmarks will likely increase.
In addition to their superior performance, the main advantage of the proposed policies is that they perform well even without a data-driven optimization. Consequently, they are less susceptible to the overfitting issues that affect the benchmark policies when only limited data are available. Under these conditions, the proposed $doa$ policies once again outperform the benchmarks significantly, with better median performance as well as excellent robustness against overfitting.
Among the different variants of the proposed policy,  we generally recommend the use of $doa1$ for obtaining replacement decisions of an $M$-component series system. In our investigations, this policy obtained the best performances for all cases, coupled with the best robustness against overfitting.

This work highlights the importance of further theory-driven investigations for the development of multi-component maintenance policies. The over-reliance of current state-of-the-art approaches on data is a major problem for applications with a limited amount of run-to-failure data, which is the case for many application domains. As we demonstrated, theory-driven approaches can help identify better and more robust policies, and thus constitute a step towards improving the efficiency and safety of engineering systems.

%
%
%
%
\section*{Acknowledgement}
\label{sec:acknowledgement}
This work was supported by the German Federal Ministry for Economic Affairs and Climate Action (BMWK) through
the aviation research program LUFO VI-3 in the project BIG-ROHU.

\bibliographystyle{unsrtnat}
\bibliography{references}

@article{koutas2025leaf,
  title={{Leaf it to renewal: Improved predictive maintenance policies via renewal theory and decision trees}},
  author={Koutas, Daniel and Straub, Daniel},
  journal={arXiv preprint arXiv:2509.20145},
  year={2025}
}

@misc{ISO13306_2017,
  author       = {{International Organization for Standardization}},
  title        = {{ISO 13306:2017 --- Maintenance terminology}},
  year         = {2017},
  howpublished = {Standard},
  address      = {Geneva, Switzerland},
  organization = {International Organization for Standardization},
  note         = {ISO 13306:2017},
}

@article{thomas1986survey,
  title={{A Survey of Maintenance and Replacement Models for Maintainability and Reliability of Multi-Item Systems}},
  author={Thomas, LC},
  journal={Reliability Engineering},
  volume={16},
  number={4},
  pages={297--309},
  year={1986},
  publisher={Elsevier}
}

@article{nicolai2008optimal,
  title={{Optimal Maintenance of Multi-Component Systems: a Review}},
  author={Nicolai, Robin P and Dekker, Rommert},
  journal={Complex System Maintenance Handbook},
  pages={263--286},
  year={2008},
  publisher={Springer}
}

@book{mobley2002introduction,
  title={{An Introduction to Predictive Maintenance}},
  author={Mobley, R Keith},
  year={2002},
  publisher={Elsevier}
}

@article{cho1991survey,
  title={A survey of maintenance models for multi-unit systems},
  author={Cho, Danny I and Parlar, Mahmut},
  journal={European Journal of Operational Research},
  volume={51},
  number={1},
  pages={1--23},
  year={1991},
  publisher={Elsevier}
}

@article{dekker1997review,
  title={{A Review of Multi-Component Maintenance Models with Economic Dependence}},
  author={Dekker, Rommert and Wildeman, Ralph E and Van der Duyn Schouten, Frank A},
  journal={Mathematical Methods of Operations Research},
  volume={45},
  number={3},
  pages={411--435},
  year={1997},
  publisher={Springer}
}

@article{alaswad2017review,
  title={A review on condition-based maintenance optimization models for stochastically deteriorating system},
  author={Alaswad, Suzan and Xiang, Yisha},
  journal={Reliability Engineering \& System Safety},
  volume={157},
  pages={54--63},
  year={2017},
  publisher={Elsevier}
}

@article{de2020review,
  title={{A review on maintenance optimization}},
  author={De Jonge, Bram and Scarf, Philip A},
  journal={{European Journal of Operational Research}},
  volume={285},
  number={3},
  pages={805--824},
  year={2020},
  publisher={Elsevier}
}

@article{radner1963opportunistic,
  title={Opportunistic replacement of a single part in the presence of several monitored parts},
  author={Radner, Roy and Jorgenson, Dale W},
  journal={Management Science},
  volume={10},
  number={1},
  pages={70--84},
  year={1963},
  publisher={INFORMS}
}

@article{vergin1977maintenance,
  title={{Maintenance Scheduling for Multicomponent Equipment}},
  author={Vergin, Roger C and Scriabin, Michael},
  journal={AIIE Transactions},
  volume={9},
  number={3},
  pages={297--305},
  year={1977},
  publisher={Taylor \& Francis}
}

@article{okumoto1983optimum,
  title={{An Optimum Group Maintenance Policy}},
  author={Okumoto, K and Elsayed, EA},
  journal={Naval Research Logistics Quarterly},
  volume={30},
  number={4},
  pages={667--674},
  year={1983},
  publisher={Wiley Online Library}
}

@article{sheu1991generalized,
  title={{A Generalized Block Replacement Policy with Minimal Repair and General Random Repair Costs for a Multi-unit System}},
  author={Sheu, Shey-Huei},
  journal={Journal of the Operational Research Society},
  volume={42},
  number={4},
  pages={331--341},
  year={1991},
  publisher={Taylor \& Francis}
}

@article{baker1995can,
  title={{Can models fitted to small data samples lead to maintenance policies with near-optimum cost?}},
  author={Baker, RD and Scarf, PA},
  journal={IMA Journal of Management Mathematics},
  volume={6},
  number={1},
  pages={3--12},
  year={1995},
  publisher={Oxford University Press}
}

@article{castanier2005condition,
  title={{A condition-based maintenance policy with non-periodic inspections for a two-unit series system}},
  author={Castanier, Bruno and Grall, Antoine and B{\'e}renguer, Christophe},
  journal={Reliability Engineering \& System Safety},
  volume={87},
  number={1},
  pages={109--120},
  year={2005},
  publisher={Elsevier}
}

@article{zhu2010availability,
  title={{Availability optimization of systems subject to competing risk}},
  author={Zhu, Yada and Elsayed, Elsayed A and Liao, Haitao and Chan, Ling-Yau},
  journal={European Journal of Operational Research},
  volume={202},
  number={3},
  pages={781--788},
  year={2010},
  publisher={Elsevier}
}

@article{bouvard2011condition,
  title={{Condition-based dynamic maintenance operations planning \& grouping. Application to commercial heavy vehicles}},
  author={Bouvard, Keomany and Artus, Samuel and B{\'e}renguer, Christophe and Cocquempot, Vincent},
  journal={Reliability Engineering \& System Safety},
  volume={96},
  number={6},
  pages={601--610},
  year={2011},
  publisher={Elsevier}
}

@article{camci2009system,
  title={{System Maintenance Scheduling With Prognostics Information Using Genetic Algorithm}},
  author={Camci, Fatih},
  journal={IEEE Transactions on Reliability},
  volume={58},
  number={3},
  pages={539--552},
  year={2009},
  publisher={IEEE}
}

@article{van2013dynamic,
  title={{A dynamic predictive maintenance policy for complex multi-component systems}},
  author={Van Horenbeek, Adriaan and Pintelon, Liliane},
  journal={Reliability Engineering \& System Safety},
  volume={120},
  pages={39--50},
  year={2013},
  publisher={Elsevier}
}

@article{nguyen2015multi,
  title={Multi-level predictive maintenance for multi-component systems},
  author={Nguyen, Kim-Anh and Do, Phuc and Grall, Antoine},
  journal={Reliability Engineering \& System Safety},
  volume={144},
  pages={83--94},
  year={2015},
  publisher={Elsevier}
}

@article{murthy1985study,
  title={Study of a multi-component system with failure interaction},
  author={Murthy, DNP and Nguyen, DG},
  journal={European Journal of Operational Research},
  volume={21},
  number={3},
  pages={330--338},
  year={1985},
  publisher={Elsevier}
}

@article{dong2007segmental,
  title={{A segmental hidden semi-Markov model (HSMM)-based diagnostics and prognostics framework and methodology}},
  author={Dong, Ming and He, David},
  journal={Mechanical systems and signal processing},
  volume={21},
  number={5},
  pages={2248--2266},
  year={2007},
  publisher={Elsevier}
}

@article{li2018remaining,
  title={{Remaining useful life estimation in prognostics using deep convolution neural networks}},
  author={Li, Xiang and Ding, Qian and Sun, Jian-Qiao},
  journal={Reliability Engineering \& System Safety},
  volume={172},
  pages={1--11},
  year={2018},
  publisher={Elsevier}
}

@article{yang2024group,
  title={{Group machinery intelligent maintenance: Adaptive health prediction and global dynamic maintenance decision-making}},
  author={Yang, Li and Zhou, Shihan and Ma, Xiaobing and Chen, Yi and Jia, Heping and Dai, Wei},
  journal={Reliability Engineering \& System Safety},
  volume={252},
  pages={110426},
  year={2024},
  publisher={Elsevier}
}

@article{chalabi2016optimisation,
  title={{Optimisation of preventive maintenance grouping strategy for multi-component series systems: Particle swarm based approach}},
  author={Chalabi, N and Dahane, Mohammed and Beldjilali, Bouziane and Neki, Abdelkader},
  journal={Computers \& Industrial Engineering},
  volume={102},
  pages={440--451},
  year={2016},
  publisher={Elsevier}
}

@book{kochenderfer2022algorithms,
  title={{Algorithms for Decision Making}},
  author={Kochenderfer, Mykel J and Wheeler, Tim A and Wray, Kyle H},
  year={2022},
  publisher={MIT press}
}

@article{he2023condition,
  title={Condition-based maintenance optimization for multi-component systems considering prognostic information and degraded working efficiency},
  author={He, Rui and Tian, Zhigang and Wang, Yifei and Zuo, Mingjian and Guo, Ziwei},
  journal={Reliability Engineering \& System Safety},
  volume={234},
  pages={109167},
  year={2023},
  publisher={Elsevier}
}

@article{luque2019risk,
  title={{Risk-based optimal inspection strategies for structural systems using dynamic Bayesian networks}},
  author={Luque, Jesus and Straub, Daniel},
  journal={Structural Safety},
  volume={76},
  pages={68--80},
  year={2019},
  publisher={Elsevier}
}

@article{nielsen2018computational,
  title={{Computational framework for risk-based planning of inspections, maintenance and condition monitoring using discrete Bayesian networks}},
  author={Nielsen, Jannie S{\o}nderk{\ae}r and S{\o}rensen, John Dalsgaard},
  journal={Structure and Infrastructure Engineering},
  volume={14},
  number={8},
  pages={1082--1094},
  year={2018},
  publisher={Taylor \& Francis}
}

@article{hinderer1999approximate,
  title={{Approximate solution of Markov renewal programs with finite time horizon}},
  author={Hinderer, Karl and Waldmann, Karl-Heinz},
  journal={SIAM Journal on Control and Optimization},
  volume={37},
  number={2},
  pages={502--520},
  year={1999},
  publisher={SIAM}
}

@book{grimmett2020probability,
  title={{Probability and Random Processes}},
  author={Grimmett, Geoffrey and Stirzaker, David},
  year={2020},
  publisher={{Oxford University Press}}
}

@book{ethier2009markov,
  title={{Markov Processes: Characterization and Convergence}},
  author={Ethier, Stewart N and Kurtz, Thomas G},
  year={2009},
  publisher={John Wiley \& Sons}
}

@article{kaelbling1996reinforcement,
  title={{Reinforcement Learning: A Survey}},
  author={Kaelbling, Leslie Pack and Littman, Michael L and Moore, Andrew W},
  journal={{Journal of Artificial Intelligence Research}},
  volume={4},
  pages={237--285},
  year={1996}
}

@book{sutton1998reinforcement,
  title={{Reinforcement Learning: An Introduction}},
  author={Sutton, Richard S and Barto, Andrew G and others},
  volume={1(1)},
  number={1},
  year={1998},
  publisher={MIT press Cambridge}
}

@article{kaelbling1998planning,
  title={Planning and acting in partially observable stochastic domains},
  author={Kaelbling, Leslie Pack and Littman, Michael L and Cassandra, Anthony R},
  journal={{Artificial Intelligence}},
  volume={101},
  number={1-2},
  pages={99--134},
  year={1998},
  publisher={Elsevier}
}

@article{bernstein2002complexity,
  title={{The Complexity of Decentralized Control of Markov Decision Processes}},
  author={Bernstein, Daniel S and Givan, Robert and Immerman, Neil and Zilberstein, Shlomo},
  journal={Mathematics of Operations Research},
  volume={27},
  number={4},
  pages={819--840},
  year={2002},
  publisher={INFORMS}
}

@article{pecht2009prognostics,
  title={{Prognostics and Health Management of Electronics}},
  author={Pecht, Michael},
  journal={{Encyclopedia of Structural Health Monitoring}},
  year={2009},
  publisher={Wiley Online Library}
}

@article{kim2017prognostics,
  title={{Prognostics and Health Management of Engineering Systems}},
  author={Kim, Nam-Ho and An, Dawn and Choi, Joo-Ho},
  journal={Switzerland: Springer International Publishing},
  year={2017},
  publisher={Springer}
}

@article{si2011remaining,
  title={{Remaining useful life estimation--A review on the statistical data driven approaches}},
  author={Si, Xiao-Sheng and Wang, Wenbin and Hu, Chang-Hua and Zhou, Dong-Hua},
  journal={{European Journal of Operational Research}},
  volume={213},
  number={1},
  pages={1--14},
  year={2011},
  publisher={Elsevier}
}

@article{ran2019survey,
  title={{A Survey of Predictive Maintenance: Systems, Purposes and Approaches}},
  author={Ran, Yongyi and Zhou, Xin and Lin, Pengfeng and Wen, Yonggang and Deng, Ruilong},
  journal={arXiv preprint arXiv:1912.07383},
  pages={1--36},
  year={2019},
  publisher={tech. rep., 12 2019}
}

@article{kamariotis2024metric,
  title={A metric for assessing and optimizing data-driven prognostic algorithms for predictive maintenance},
  author={Kamariotis, Antonios and Tatsis, Konstantinos and Chatzi, Eleni and Goebel, Kai and Straub, Daniel},
  journal={Reliability Engineering \& System Safety},
  volume={242},
  pages={109723},
  year={2024},
  publisher={Elsevier}
}

@article{bismut2021optimal,
  title={{Optimal Adaptive Inspection and Maintenance Planning for Deteriorating Structural Systems}},
  author={Bismut, Elizabeth and Straub, Daniel},
  journal={Reliability Engineering \& System Safety},
  volume={215},
  pages={107891},
  year={2021},
  publisher={Elsevier}
}

@book{kochenderfer2015decision,
  title={{Decision Making Under Uncertainty: Theory and Application}},
  author={Kochenderfer, Mykel J},
  year={2015},
  publisher={MIT press}
}

@article{zhang2020deep,
  title={{Deep reinforcement learning for condition-based maintenance planning of multi-component systems under dependent competing risks}},
  author={Zhang, Nailong and Si, Wujun},
  journal={Reliability Engineering \& System Safety},
  volume={203},
  pages={107094},
  year={2020},
  publisher={Elsevier}
}

@article{chen2023deep,
  title={{A deep reinforcement learning approach for maintenance planning of multi-component systems with complex structure}},
  author={Chen, Jiahao and Wang, Yu},
  journal={Neural Computing and Applications},
  volume={35},
  number={21},
  pages={15549--15562},
  year={2023},
  publisher={Springer}
}

@article{najafi2023deep,
  title={{A deep reinforcement learning approach for repair-based maintenance of multi-unit systems using proportional hazards model}},
  author={Najafi, Seyedvahid and Lee, Chi-Guhn},
  journal={Reliability Engineering \& System Safety},
  volume={234},
  pages={109179},
  year={2023},
  publisher={Elsevier}
}

@article{zhang2023guided,
  title={{Guided probabilistic reinforcement learning for sampling-efficient maintenance scheduling of multi-component system}},
  author={Zhang, Yiming and Zhang, Dingyang and Zhang, Xiaoge and Qiu, Lemiao and Chan, Felix TS and Wang, Zili and Zhang, Shuyou},
  journal={Applied Mathematical Modelling},
  volume={119},
  pages={677--697},
  year={2023},
  publisher={Elsevier}
}

@article{smith1958renewal,
  title={{Renewal Theory and Its Ramifications}},
  author={Smith, Walter L},
  journal={Journal of the Royal Statistical Society Series B: Statistical Methodology},
  volume={20},
  number={2},
  pages={243--284},
  year={1958},
  publisher={{Oxford University Press}}
}

@article{doob1948renewal,
  title={{Renewal Theory From the Point of View of the Theory of Probability}},
  author={Doob, Joseph L},
  journal={Transactions of the American Mathematical Society},
  volume={63},
  number={3},
  pages={422--438},
  year={1948},
  publisher={JSTOR}
}

@book{raiffa2000applied,
  title={{Applied Statistical Decision Theory}},
  author={Raiffa, Howard and Schlaifer, Robert},
  volume={78},
  year={2000},
  publisher={John Wiley \& Sons}
}

@article{bellman1966dynamic,
  title={{Dynamic Programming}},
  author={Bellman, Richard},
  journal={Science},
  volume={153},
  number={3731},
  pages={34--37},
  year={1966},
  publisher={American Association for the Advancement of Science}
}

@article{pineau2006anytime,
  title={{Anytime Point-Based Approximations for Large POMDPs}},
  author={Pineau, Joelle and Gordon, Geoffrey and Thrun, Sebastian},
  journal={{Journal of Artificial Intelligence Research}},
  volume={27},
  pages={335--380},
  year={2006}
}

@article{nguyen2019new,
  title={A new dynamic predictive maintenance framework using deep learning for failure prognostics},
  author={Nguyen, Khanh TP and Medjaher, Kamal},
  journal={Reliability Engineering \& System Safety},
  volume={188},
  pages={251--262},
  year={2019},
  publisher={Elsevier}
}

@article{gneiting2007probabilistic,
  title={Probabilistic forecasts, calibration and sharpness},
  author={Gneiting, Tilmann and Balabdaoui, Fadoua and Raftery, Adrian E},
  journal={Journal of the Royal Statistical Society Series B: Statistical Methodology},
  volume={69},
  number={2},
  pages={243--268},
  year={2007},
  publisher={{Oxford University Press}}
}

@article{whitley1994genetic,
  title={{A Genetic Algorithm Tutorial}},
  author={Whitley, Darrell},
  journal={Statistics and Computing},
  volume={4},
  number={2},
  pages={65--85},
  year={1994},
  publisher={Springer}
}

@article{karafotias2014parameter,
  title={{Parameter Control in Evolutionary Algorithms: Trends and Challenges}},
  author={Karafotias, Giorgos and Hoogendoorn, Mark and Eiben, {\'A}goston E},
  journal={IEEE Transactions on Evolutionary Computation},
  volume={19},
  number={2},
  pages={167--187},
  year={2014},
  publisher={IEEE}
}

@incollection{lobo2007adaptive,
  title={{Adaptive Population Sizing Schemes in Genetic Algorithms}},
  author={Lobo, Fernando G and Lima, Claudio F},
  booktitle={Parameter Setting in Evolutionary Algorithms},
  pages={185--204},
  year={2007},
  publisher={Springer}
}

@incollection{goldberg1991comparative,
  title={{A Comparative Analysis of Selection Schemes Used in Genetic Algorithms}},
  author={Goldberg, David E and Deb, Kalyanmoy},
  booktitle={Foundations of Genetic Algorithms},
  volume={1},
  pages={69--93},
  year={1991},
  publisher={Elsevier}
}

@book{back2000evolutionary,
  title={{Evolutionary Computation 1: Basic Algorithms and Operators}},
  author={B{\"a}ck, Thomas and Fogel, David B and Michalewicz, Zbigniew},
  year={2000},
  publisher={IOP Publishing}
}

@article{deb1995simulated,
  title={{Simulated Binary Crossover for Continuous Search Space}},
  author={Deb, Kalyanmoy and Agrawal, Ram Bhushan and others},
  journal={Complex Systems},
  volume={9},
  number={2},
  pages={115--148},
  year={1995},
  publisher={[Champaign, IL, USA: Complex Systems Publications, Inc., c1987-}
}

@book{mitchell1998introduction,
  title={{An Introduction to Genetic Algorithms}},
  author={Mitchell, Melanie},
  year={1998},
  publisher={MIT press}
}

@article{deb2002fast,
  title={{A Fast and Elitist Multiobjective Genetic Algorithm: NSGA-II}},
  author={Deb, Kalyanmoy and Pratap, Amrit and Agarwal, Sameer and Meyarivan, TAMT},
  journal={IEEE Transactions on Evolutionary Computation},
  volume={6},
  number={2},
  pages={182--197},
  year={2002},
  publisher={Ieee}
}

@article{golberg1989genetic,
  title={{Genetic Algorithms in Search, Optimization \& Machine Learning}},
  author={Golberg, David E},
  journal={Addison-Wesley},
  volume={1989},
  number={102},
  pages={36},
  year={1989}
}

@article{mitchell1994evolving,
  title={{Evolving Cellular Automata to Perform Computations: Mechanisms and Impediments}},
  author={Mitchell, Melanie and Crutchfield, James P and Hraber, Peter T},
  journal={Physica D: Nonlinear Phenomena},
  volume={75},
  number={1-3},
  pages={361--391},
  year={1994},
  publisher={Elsevier}
}

@article{van2000mean,
  title={{Mean and Variance of Ratio Estimators Used in Fluorescence Ratio Imaging}},
  author={Van Kempen, GMP and Van Vliet, LJ},
  journal={Cytometry: The Journal of the International Society for Analytical Cytology},
  volume={39},
  number={4},
  pages={300--305},
  year={2000},
  publisher={Wiley Online Library}
}

@article{morrison2016branch,
  title={{Branch-and-bound algorithms: A survey of recent advances in searching, branching, and pruning}},
  author={Morrison, David R and Jacobson, Sheldon H and Sauppe, Jason J and Sewell, Edward C},
  journal={Discrete Optimization},
  volume={19},
  pages={79--102},
  year={2016},
  publisher={Elsevier}
}

@article{frazier2018tutorial,
  title={{A Tutorial on Bayesian Optimization}},
  author={Frazier, Peter I},
  journal={arXiv preprint arXiv:1807.02811},
  year={2018}
}

\newpage
\appendix
\counterwithin{figure}{section}
\counterwithin{table}{section}

\label{appendix}

%
%
%
%
\section{Parameter optimization via genetic algorithms}
\label{sec:app_optimization_via_ga}

The basic functionality of GAs is to iteratively evolve a population of candidate solutions via a set of selection and recombination operations to generate new samples in the search space. The goal is to find good solutions by preserving favourable features and eliminating undesirable traits \cite{whitley1994genetic}. The following paragraphs outline the steps of the specific GA used in this work, along with a more detailed explanation of the workings and their application to our particular decision setting.

\subsection{Initialization} 
    
The first step of GAs consists of generating a population of candidate solutions, also called individuals. Each individual is represented by a set of parameters (genes). For an $M$-component series system, an individual of $rh1$ is represented by the parameter vector $[p_{\mathrm{thres},1}^{(i)},~\ldots,~p_{\mathrm{thres},M}^{(i)}]$.

The population size $n_p$ is a crucial parameter: a value that is too low may hinder the GA from finding high-quality solutions, whereas a value that is too large wastes computational resources. However, finding a good a-priori estimate of the optimal population size is difficult, as it is also problem-dependent. Thus, a common approach, and also the one used in this work, is to do some preliminary experiments with a number of different population sizes to see which one works best. For a more detailed discussion on choosing the optimal population size, the reader is referred to \cite{lobo2007adaptive}.   

Given $n_p$, to generate an initial population for the $M$-dimensional parameter vector of $rh1$, we use Monte Carlo sampling. In our decision setting, corrective replacements have higher consequences than preventive replacements. To incorporate this domain knowledge, we therefore bias the initial samples $p_{thres,i}^{(i)}$ towards lower failure probabilities. Specifically, we choose the Beta distribution, whose PDF is described by:
\begin{equation}
    f_X(x; a,b) = \frac{\Gamma(a+b)}{\Gamma(a)\Gamma(b)} x^{\alpha-1} (1-x)^{\beta-1}, \qquad x \in (0,1),~0<\alpha,\beta
\end{equation}
with shape parameters $\alpha,\beta$ and gamma function $\Gamma(\cdot)$. Setting the mode of the distribution to be at $c_{r,i}=\frac{c_f+c_{v,i}}{c_c}$ (i.e., the default probability threshold of $rh1$) and choosing a large variance to still cover the whole parameter space, leads to
\begin{equation}
    \label{eq:beta_params}
    \beta=3 \quad \& \quad \alpha_i=\frac{1+(\beta-2)c_{r,i}}{1-c_{r,i}}.
\end{equation}

\subsection{Candidate Evaluation}

Each individual of the population is evaluated by a \emph{fitness} function and, thus, assigned a fitness. In our replacement setting, the fitness of each individual corresponds to the heuristic's performance achieved with the individual's parameters. As visualized in \Cref{subfig:rh2_noisy_target_10}, for very low sample domains, this performance can vary heavily, even for small parameter changes.

\subsection{Parent selection}

The next GA step involves the selection of ``parent'' candidates from which ``offspring'' candidates are generated. Generally, the selection procedure favours individuals with better fitness value \cite{back2000evolutionary}. In this work, we choose a tournament-based selection, which draws $n_t$ individuals randomly (with uniform distribution) from the population, selects the individual with the highest fitness from this group for further genetic processing, and repeats this procedure until the desired number of offspring is reached. We choose this tournament-based selection for its simplicity, fast selection, and easy parallelization possibility. A more detailed comparison of various selection schemes for GAs can be found in \cite{goldberg1991comparative}.

\subsection{Genetic processing}

Once parents are selected, the two main options for further genetic processing are \emph{crossover} (recombination) and \emph{mutation}, with many different variations, combinations, and importance weights \cite{karafotias2014parameter,back2000evolutionary}.

Crossover defines a certain combination operation of parents to generate offspring, where, ideally, the parents' beneficial traits are preserved and merged. Hence, crossover corresponds to exploitation of high-quality genetic material \cite{mitchell1998introduction}. Mutation, on the other hand, represents random variations of an individual's traits and hence corresponds to exploration of the search space \cite{back2000evolutionary,mitchell1998introduction}. Given a set of parents, the occurrence probability of crossover and mutation is set by $p_c$ and $p_m$, respectively. Generally, mutation can be seen as a lower-probability event compared to crossover, to reintroduce lost genetic material \cite{golberg1989genetic}.\footnote{Note that in this work, first the crossover and then mutation operator are applied to an individual independently; hence an individual could undergo both operations. In other evolutionary algorithm works, an individual can only be subject to one operation.}

In this work, we employ simulated binary crossover (SBX) introduced in \cite{deb1995simulated} for the recombination of our continuous and bounded threshold parameters ($p_{thres,i}\in [0,1]$). Here, offspring are sampled somewhere between (or slightly beyond) the parents' values. The probability distribution over the child's distance to its parents is controlled by the distribution index $\eta_c$. The details of SBX are explained in \cite{deb1995simulated}.

As a complementary operator to SBX, we use polynomial mutation to introduce random perturbations to individual variables. The probability distribution over the magnitude of the changes is again controlled by a distribution index $\eta_m$. Once an individual is randomly selected for mutation, the genes are randomly and independently selected for mutation with probability $p_{m,i}$ \cite{deb2002fast}.

\subsection{Generation selection}

For the selection of the next generation, we choose an \emph{elitist} approach described in \cite{mitchell1994evolving}. Here, the $n_e$ individuals with the best fitness of the current population, i.e., the elites, are directly copied into the next generation without modification. The remaining $n_p-n_e$ offspring are generated by repeating the process described in steps 3 and 4. This elitist approach ensures the survival of the highest-quality genetic material and upkeeps selection pressure \cite{mitchell1994evolving,mitchell1998introduction}.

A summary of all GA hyperparameter used in this work is provided in \Cref{tab:general_ga_params}. For a more detailed overview on parameter selection of genetic algorithms, the reader is referred to \cite{karafotias2014parameter}.
\begin{table}[H]
    \renewcommand{\arraystretch}{1.2}
    \setlength{\abovecaptionskip}{8pt}
    \centering    
    \begin{tabular}{ll}
        \toprule
        Parameter & Value \\
        \midrule
        Population size $n_p$ & 250 \\
        Tournament size $n_t$ & 3 \\
        Cross-over probability $p_{c}$ & 0.9 \cite{deb2002fast}\\
        Cross-over distribution index $\eta_c$ & 20 \cite{deb2002fast} \\
        Mutation selection probability $p_{m}$ & $0.2$ \\
        Per-gene mutation probability $p_{m,i}$ & $1/M$ \cite{deb2002fast} \\
        Mutation distribution index $n_m$ & 20 \cite{deb2002fast} \\
        \# Elites $n_e$ & 50 \\
        \# Generations $n_g$ & 25 \\
        \bottomrule
    \end{tabular}
    \caption{GA parameters for $rh2$.}
    \label{tab:general_ga_params}
\end{table}
%

%
%
%
%
%
%
\section{Virtual RUL simulator}
\label{sec:virtual_RUL_simulator}
This section is more or less an iteration of the RUL simulator description in \cite{kamariotis2024metric}.

We assume an underlying normal distribution of the time-to-failure of the components, i.e., $T_\mathrm{F}\sim\mathcal{N}(\mu=225,\sigma=40)$. At each time step $t_k$, we make a prediction about the remaining useful life of the component, denoted as $RUL_k$. This prediction deviates from the true RUL, denoted $RUL_{\mathrm{true}, k}=T_\mathrm{F}-t_k$, by a certain error $\epsilon_k$. We define this prediction in a logarithmic scale as:
\begin{equation}
    \label{eq:rul_pred_error}
    \ln\left(RUL_k \right) = \ln\left(RUL_{\mathrm{true},k}\right) + \ln\left( \epsilon_k \right).
\end{equation}
In this work, we assume a multivariate normal (MVN) distribution for the logarithmic prediction errors over time. Specifically, 
\begin{equation}
    \label{eq:mvn_of_errors}
    [\ln(\epsilon_1), \dots, \ln(\epsilon_n)] \sim \mathrm{MVN}(\mathbf{0}, \mathbf{\Sigma}).
\end{equation}
The mean of the logarithmic prediction errors is 0. The covariance matrix $\mathbf{\Sigma}$ (not to be confused with the strategy $\Sigma$ in \Cref{sec:pdm_policies}!) is constructed as
\begin{equation}
    \label{eq:cov_mat}
    \mathbf{\Sigma} = \mathbf{D} \cdot \mathbf{R} \cdot \mathbf{D},
\end{equation}
where $\mathbf{D}$ is a diagonal matrix containing the standard deviation of the prediction errors $\sigma_{\ln(\epsilon_k)}$ on the main diagonal. For this work, we fix $\sigma_{\ln(\epsilon_k)}=0.4~\forall k$. On the other hand, $\mathbf{R}$ is a correlation coefficient matrix. If the predictions of the prognostic model are assumed to be independent at $t_k$, then $\mathbf{R}$ reduces simply to a diagonal matrix with $1$ on the main diagonal. In this work, we infuse a certain structure into the prognostic model's predictions by assuming:
\begin{equation}
    \label{eq:corr_mat}
    \mathbf{R} = [\rho_{ij}], \qquad \rho_{ij} = \exp \left(-\frac{|t_i-t_j|}{l} \right),
\end{equation}
where $l$ is the hyperparameter defining the correlation length.

The procedure to generate the RUL-PDFs at each time step is then:
\begin{enumerate}
    \item Draw a failure sample $\Tfi$ from the underlying time-to-failure distribution.
    \item Compute the true RUL values at each time step via $RUL_k^{(i)}=\Tfi-t_k$.
    \item Sample an error vector $[\ln(\epsilon_1^{(i)}), \dots, \ln(\epsilon_n^{(i)})]$ from the MVN distribution in \Cref{eq:mvn_of_errors}.
    \item Compute the mean values of the prediction via \Cref{eq:rul_pred_error}: $\mu_{\ln\left(RUL_{k} \right)}^{(i)}=\ln(RUL_{\mathrm{true},k}^{(i)}) + \ln(\epsilon_k^{(i)})$.
    \item The predicted RUL-PDF is then a lognormal distribution, i.e., $\ln\left(RUL_k^{(i)} \right) \sim \mathcal{N}\left(\mu_{\ln\left(RUL_k \right)}^{(i)}, \sigma_{\ln(\epsilon_k)}\right)$.
\end{enumerate}

\end{document}